\documentclass[10pt]{article}

\usepackage{float} 
\usepackage{wrapfig}

\usepackage{amsmath,amsfonts,amssymb,amsthm,epsf}
\input epsf 
\usepackage{calc}
\usepackage{epsfig}
\usepackage{epstopdf}
 
\numberwithin{equation}{section}
\newtheorem{theorem}{Theorem}[section]

\newtheorem{proposition}[theorem]{Proposition}

\newtheorem{assumption}{Assumption}

\DeclareMathOperator{\Span}{span}

\newcommand{\RR}{{\rm I\kern -1.6pt{\rm R}}}
\newcommand{\beq}{\begin{equation}} 
\newcommand{\eeq}{\end{equation}}
\newcommand{\bes}{\begin{split}} 
\newcommand{\es}{\end{split}}
\newcommand{\bef}{\begin{equation}} 
\newcommand{\ef}{\end{equation}}

\newcommand{\beginsupplement}{%
        \setcounter{table}{0}
        \renewcommand{\thetable}{S\arabic{table}}%
        \setcounter{figure}{0}
        \renewcommand{\thefigure}{S\arabic{figure}}%
     }

\begin{document}
\voffset= -0.5in
\title{Error bounds and analysis of proper orthogonal decomposition model reduction methods using snapshots from the solution and the time derivatives}
\author{Tanya Kostova, Geoffrey Oxberry, Kyle Chand and William Arrighi }
\maketitle
\begin{abstract}
Proper orthogonal decomposition methods for model reduction utilize information about the solution at certain time and parameter points to generate a reduced space basis. In  this paper, we compare two proper orthogonal decomposition methods for reducing large systems of ODEs.  The first method is based on collecting snapshots from the solutions only; the second method uses snapshots from both the solutions and their time derivatives. To compare the methods, we derive new bounds for the 2-norm of the approximation error induced by the each of the methods. The bounds are represented as a sum of two terms: the first depends on the size of the first neglected singular value while the second depends only on the spacings between the snapshots. We performed numerical experiments to compare the errors from the two model reduction  methods applied to the semidiscretized FitzHugh-Nagumo system and  investigated the relation between the behavior of the numerically observed  error and the error bounds. We find that the error bounds, though not tight, provide insights and justification for using time derivative snapshots in POD model reduction for dynamical systems.
\end{abstract}

\bigskip

{\bf Authors' Address:}

Lawrence Livermore National Laboratory, 7000 East Avenue, Livermore, CA 94551, USA; 

\bigskip

{\bf E-mail of corresponding author:} kostova@llnl.gov

\vskip 2in
{\bf \Large IM release number: LLNL-JRNL-663838}

\vfill\eject
\section{Introduction}

In many areas of science and technology, complex multi-physics  time-dependent 
problems are modeled by large systems of differential 
equations. Their analysis often poses huge computational challenges as it requires
multiple,  prohibitively expensive simulations, in terms of  time and
memory.  As known from the theory of dynamical systems, seemingly complex high-dimensional dynamical systems can have low - dimensional global attractors, or low-dimensional invariant manifolds containing the physically relevant solutions of the system \cite{guck, rega, foias}.  In these cases, reduced order models (ROMs)
exploiting the attractor or center manifold structure can dramatically reduce the time
and memory needed to execute the corresponding full-order model (FOM).

One of the most studied approaches to building reduced order models is based on 
Proper Orthogonal Decomposition (POD) \cite{sirovich,  antoulas}. Normally, POD
 uses  the
(vector) solution of the FOM at selected time, space and parameter values, typically 
called "snapshots" \cite{sirovich, burkardt}, to calculate a set of singular (orthonormal) vectors spanning the full space. In the context of POD, if the dimension of the FOM is $n$ (a large number), for a given $l<n$, the span of the $l$ vectors corresponding to the $l$ largest singular values defines the reduced basis space (RBS).   The ROM is defined and solved in the RBS and the prolongation of its solution back in $n$-dimensional space serves as an approximation to the solution of the FOM. The ultimate goal is to find an optimal (ideally small)  $l$ so that the approximate solution is sufficiently close to the solution of the FOM. Two typical problems are important for POD ROMs: (a) given a fixed error tolerance, find the the snapshots generating the smallest basis
  satisfying that tolerance; (b) given a fixed basis size, find the snapshots generating a basis that minimizes some error indicator. The specific choice of snapshots used affects
greatly the size of the ROM basis and the ROM solution error.

The availability of {\it a priori error estimates} or   {\it error bounds} that include characteristics of the RBS (e.g., dimension, exact locations of the snapshots, type of snapshots, etc.)  is important for understanding the origin of the approximation error and controlling it. {\it A priori} estimates could be used directly to determine the optimal reduced dimension of the model or to develop {\it a posteriori} error estimators to approximate it computationally. Although several papers deriving  error bounds for  POD ROMs have
been published, e.g.  \cite{hinze, iliescu,
 kv1, kv2,  rathinam} the bounds rarely contain information about the spacing between snapshots and none explicitly includes the exact locations of the snapshots. Including the latter information can be potentially useful for for informing snapshot selection and lower the computational resource requirements of POD. Notably, Kunisch and Volkwein \cite{kv1}  developed error bounds for POD ROMs for linear parabolic problems which are in the form of a sum of two terms: one depending on the singular values corresponding to the orthogonal to the RBS basis vectors  and the other depending on  the (uniform) distance between the time points where the snapshots were taken.  This result was generalized and applied to other problems in subsequent publications, {\it e.g.}  \cite{kv2, hinze}. In the latter two publications the authors used a POD ROM based on two types of snapshots: a) collected from the solution and b) difference quotients (DQs) calculated from these snapshots. Although the DQs are in the span of the solution snapshots, their use was justified by a need to avoid the appearance of a blow-up term in the derivation of the error bound. 

A few other authors have used DQs or  values of the time derivative as snapshots. In a variant of POD ROM, the Discrete Empirical Interpolation Method, applied to reduce nonlinear dynamical systems, Chaturantabut and Sorensen \cite{chatur1, chatur2}  used (time-derivative)  snapshots from  the nonlinear part of the system. The incentive for using derivative snapshots in their method was to  reduce the computational complexity of the POD ROM. However, so far, there is no clear answer whether using DQs or time derivative snapshots provides an advantage in terms of accuracy of approximation of the FOM solution.  In a recent paper,  
Iliescu and Wang \cite{iliescu} posed the question whether there is advantage (and if so, what) to include the DQs in the POD calculation. The question was prompted by results of two groups of authors  \cite{chapelle, singler} who developed convergence analysis for POD ROMS that did not use DQs ({\it i.e.} providing evidence that using DQs is not necessary for avoiding blow-up of the POD error), as well as by conflicting numerical results on the benefit of using DQs in POD ROM formulations \cite{iliescu}. Iliescu and Wang defined a POD optimality criterion and explored the optimality, in several norms, of two POD methods - the first utilizing only solution snapshots, and the second - utilizing both solution and DQ snapshots.  They found that the method with DQ snapshots was optimal in all norms while the one using only solution snapshots was optimal only in some of the norms.  

In this paper, we compare two POD methods for reducing large systems of ODEs.  The first POD ROM method is based on snapshots from the solutions only; the second method uses snapshots from both the solutions and their time derivatives. The time derivative snapshots in general do not belong to the span of the solution snapshots, therefore the number of singular vectors corresponding to non-zero singular values  is normally twice larger ($2l$) in the second case compared to the first ($l$).  The additional computational cost to  calculate the time derivative snapshots is negligible and thus, this is a cheap way to augment an existing basis; yet using derivative snapshots brings in additional information about the dynamics of the solutions. Thus, a method utilizing derivative information could produce more accurate ROM approximations which would be valuable when $l$ is relatively small. These arguments provide an  incentive for comparing these  two methods. 

Our initial  incentive, however, came from the expression for the error (\ref{err}), where we noticed that if the vector ${\bf f}({\bf x}, t)$ was a snapshot, it would be in the span of the reduced basis vectors and the second term in the right-hand side of (\ref{err}) would be zero. Intuitively, this property could make the error smaller compared to the case when this term is not zero. We explore this conjecture by  deriving error bounds for the two methods. We have derived bounds in a form consisting of two terms: one depending on the largest neglected singular value and one depending on the distance between snapshots but not on the singular value. Assuming that the largest neglected singular value is small, we show that the bound of the  2-norm of the error has a second-order dependence on the length of the time intervals between snapshots when only solution snapshots are used, and a fourth order dependence on same length when also time derivative snapshots are used. 

To the best of our knowledge, this paper introduces two main new results.
The first one is the actual form of the error bound which involves the time moments of the snapshots (via the time intervals between the snapshots) and the largest neglected singular value (instead of all neglected singular values). Error bounds and error estimates involving the time moments of the snapshots contain information which is potentially significant for rational snapshot selection.  

The second, and more significant,  innovation is that the bounds allow for comparison between the two POD ROM methods (with and without derivative snapshots) and specifically suggest that if the first neglected singular value is  sufficiently small, the method with time derivative information could be more accurate. Thus, these bounds give insights that we test numerically. Interestingly, we find that the behavior of the errors in numerical experiments with cases of the discretized FitzHugh-Nagumo system can be explained  by the form of the bounds, although these bounds are not tight.

The paper is organized as follows. Section 2 includes preliminary information and derivations. Section 3 includes derivations of the error bounds for the two methods. Section 4 includes three numerical experiments illustrating the the validity of the insights obtained from the error bounds.

\section{\label{s2}{Notations and relevant background}} 
Throughout the paper we will indicate vector- and  matrix- valued variables by bold upper-
and lower-case letters and scalar values by normal typesetting.
We consider a dynamical system of the form
\begin{equation}
\label{system1}
\dot {\bf x}= {\bf f}({\bf x}(t),  t), {\bf x}(0)={\bf x}^0, t\in[0, T], 
\end{equation}
where, respectively, boldface letters denote vectors, ${\bf x}(t)=(x_1(t), . . ., x_n(t))^T \in {\mathbb R}^n,  {\bf f}:  {\mathbb R}^{n+1}\to  {\mathbb R}^n$. 
\begin{assumption}\label{ass1}
It is assumed that $\frac{\partial f_i}{\partial x_j}\in {\mathcal C}( {\mathbb R}^n)$  and $\frac{\partial f_i}{\partial t} \in {\mathcal C}(0,T)$. Further smoothness assumption on ${\bf f}$ will be imposed in Proposition \ref{s2}.
\end{assumption}
\bigskip
 
Let  ${\bf y}(t)$ be a solution of (\ref{system1}).  In the following theoretical treatment it is assumed that snapshots of the solution and its derivatives are collected at certain time points from ${\bf y}(t)$.  We will explore the approximation error incurred by using a reduced order model (ROM) constructed from a POD basis using combinations of these snapshots.
We  first introduce some definitions and derive useful relationships.

\subsection {\label{s1} Using only  ${\bf y}(t_i)$  as snapshots.} 
For the sake of self-consistency and setting up notations, we will revisit some well known definitions and results related to POD. 

\bigskip 

We consider a set of time points, not necessarily equispaced, $t_i \in [0, T], i=1,. . . , m\leq n$ where the latter inequality has been set mostly for clarity of presentation.
We consider the matrix $ Y=\Bigl[ {\bf y}(t_1) \vdots  . . . \vdots {\bf y}(t_m) \Bigr] \in {\mathbb R}^{n\times m}$, with $rank({\bf Y})=m^Y$. Then,

(i)  ${\bf Y}{\bf Y}^T$ is $n\times n$ matrix, ${\bf Y}^T{\bf Y}$ is an $m\times m$ matrix  and $rank ({\bf Y}{\bf Y}^T)= rank({\bf Y}^T{\bf Y})=m^Y\leq m$.  The squares of the $m$ largest singular values of ${\bf Y}$, $(\sigma^Y_i)^2, i=1, . . ., m$,  are the eigenvalues of ${\bf Y}^T{\bf Y}$, while the rest are equal to zero.  

\medskip

(ii) Let ${\bf u}^Y_i=(u^Y_{1i}, u^Y_{2i}, . . ., u^Y_{ni})^T\in {\mathbb R}^n$ and  ${\bf v}^Y_j =(v^Y_{1j}, v^Y_{2j}, . . ., v^Y_{mj})^T\in {\mathbb R}^m$ be orthonormal  sets of eigenvectors of  ${\bf Y}{\bf Y}^T$ and ${\bf Y}^T{\bf Y}$.  The vectors ${\bf u}^Y_i\in {\mathbb R}^n$ and  ${\bf v}^Y_j \in {\mathbb R}^m$  are called {\it left} and {\it right singular vectors} of ${\bf Y}$, respectively. Let ${\bf U}^Y$ and ${\bf V}^Y$ be the matrices whose columns consist of these eigenvectors, {\it i.e.} ${\bf U}^Y({\bf U}^Y)^T={\bf I}_{n}, {\bf V}^Y({\bf V}^Y)^T={\bf I}_m$. 
It is known that
\beq
\label{Yu}
{\bf Y}^T {\bf u}^Y_i=\sigma^Y_i {\bf v}^Y_i
\eeq 
and 
\beq
\label{Yv}
{\bf Y} {\bf v}^Y_i=\sigma^Y_ i{\bf u}^Y_i,  i=1,. . ., m,
\eeq
where $\sigma^Y_i\geq 0$. 

Thus, ${\bf Y}{\bf V}^{Y}={\bf U}^{ Y}{\bf \Sigma}^Y$
where ${\bf \Sigma}^Y$ is $n\times m$ matrix with diagonal entries $\sigma^Y_i>0, i=1, . . . , m^Y$  and 0s otherwise, and it is assumed that $\sigma^Y_1\geq\sigma^Y_2\geq . . . \geq \sigma^Y_{m^Y}$.
This leads to the well known singular decomposition of ${\bf Y}$ ({\it e.g.} \cite{antoulas} ):
\beq
\label{Y}
{\bf Y}={\bf U}^Y{\bf \Sigma}^Y ({\bf V}^Y)^T,
\eeq
and the {\it elementwise dyadic decomposition} 

\beq
\label{repr}
{\bf y}(t_i)=\sum_{k=1}^{m^Y} \sigma^Y_k v^Y_{ik} {\bf u}^Y_{k} .
\eeq

(iii) In practice, once $m$ snapshots are calculated, from them  $l\leq m^Y\leq m$ left singular vectors corresponding to the $l$ greatest singular values of ${\bf Y}$ are used to construct the basis for the ROM. The rationale is based on the Schmidt-Eckart-Young-Mirsky theorem \cite{antoulas} which asserts that the truncated dyadic decomposition  $${\bf X}=\sum^l_{i=1} \sigma_i {\bf u}_i {\bf v}_i$$ minimizes $||{\bf Y}-{\bf X}||_2$ over all ${\bf X}$ such that $\operatorname{rank} {\bf X}=l$. The  value of $l$ is typically chosen so that either the reduced basis approximates the snapshots with sufficient accuracy, or to satisfy constraints on the ROM (memory usage and computational time are both related to $l$). From (\ref{repr}) the following 2-norm  bound of the difference between the snapshot ${\bf y}(t_i)$ and its $l$-truncated representation is derived:

\beq
\label{erest1}
\bes
&||{\bf y}(t_i)-\sum_{k=1}^{l} \sigma^Y_k v^Y_{ik} {\bf u}^Y_k||_2= ||\sum_{k=l+1}^{m^Y} \sigma^Y_k v^Y_{ik} {\bf u}^Y_k||_2 \\
&= \sqrt{\langle \sum_{k=l+1}^{m^Y} \sigma^Y_k v^Y_{ik} {\bf u}^Y_k, \sum_{k=l+1}^{m^Y} \sigma^Y_k v^Y_{ik} {\bf u}^Y_k\rangle} = \sqrt{\sum_{k=l+1}^{m^Y}( \sigma^Y_k v^Y_{ik} {\bf u}^Y_k)^2}\\
&=\sqrt{\sum_{k=l+1}^{m^Y}( \sigma^Y_k v^Y_{ik})^2} \leq \sigma^Y_{l+1}\sqrt{\sum_{k=l+1}^{m^Y}(v^Y_{ik})^2}\leq \sigma^Y_{l+1}, 
\end{split}
\eeq
where $\langle , \rangle$ denotes the Euclidean dot product. 

\bigskip
In the sections that follow, we will use the following partitionings and notations. Given a matrix ${\bf Y} = {\bf U}^Y{\bf \Sigma}^Y ({\bf V}^Y)^T$, we partition the columns of ${\bf U}^Y$ as follows.
\beq
{\bf U}^Y=\Bigl[\tilde {\bf U}^Y \vdots \breve {\bf U}^Y \vdots  \bar {\bf U}^Y \Bigr]
\eeq
where $\tilde {\bf U}^Y=[{\bf u}_1 \vdots . . . \vdots {\bf u}_l]$ is the {\it reduced basis matrix} ,   $\breve {\bf U}^Y=[{\bf u}_{l+1} \vdots  . . .  \vdots  {\bf u}_{m^Y}]$  contains the truncated basis vectors with non-zero singular values and $\bar {\bf U}^Y=[{\bf u}_{m^Y+1}   \vdots  . . . \vdots {\bf u}_n]$ is the ${\bf Y}{\bf Y}^T$ {\it null space matrix}: span$(\bar{\bf U}^Y)=\ker({\bf YY}^T)$. Note that  $\breve {\bf U}^Y$ and $\bar {\bf U}^Y$ may be empty matrices. 

\bigskip 

Respectively,  ${\bf u}^Y_{i}, i=1, . . ., l\leq m^Y$ are denoted as $\tilde {\bf u}^Y_i$; if $l<m^Y$, ${\bf u}_{l+1}, . . . , {\bf u}^Y_{m^Y}$ are denoted as $\breve {\bf u}^Y_i, i=1, . . ., m^Y-l$; if $m^Y<n$, ${\bf u}^Y_{m^Y+1},..., {\bf u}^Y_n$ are denoted as $\bar {\bf u}^Y_i, i=1, .  . ., n-m^Y$.

\bigskip

Since rank$({\bf Y}{\bf Y}^T)=m^Y, {\bf Y}{\bf Y}^T\bar{\bf U}=0_{n\times (n-m^Y)}$ (the null $n\times (n-m^Y)$ matrix). Further, from (\ref{repr}) and the orthonormality of ${\bf U}^Y$ it follows   $ {\bf Y}^T \bar {\bf U}^Y=0_{m\times (n-m)}. $
\medskip

\subsection{Using both  ${\bf y}(t_i)$ and ${\bf f}({\bf y}(t_i), t_i)$ as snapshots.} 
Consider a collection of $2m$  snapshots from both  ${\bf y}(t_i)$ and  ${\bf f}({\bf y}(t_i),  t_i), i=1, . . ., m$ with $2m\leq n$. (Again, the latter inequality is not a necessary condition, but rather imposed for clarity of exposition.)
Let ${\bf Z}$ be the matrix made of these vectors.
\beq
{\bf Z}=[{\bf y}(t_1) \, \vdots \,  . . . \, \vdots \, {\bf y}(t_m) \, \vdots \,  {\bf f}({\bf y}(t_1), t_0) \, \vdots \, . . . \, \vdots \, {\bf f}({\bf y}(t_m), t_m)]
\in {\mathbb R}^{n\times 2m}.  
\eeq
${\bf Z}{\bf Z}^T\in {\mathbb R}^{n\times n}$ and ${\bf Z}^T{\bf Z} \in {\mathbb R}^{2m\times 2m}.$

Let ${\bf u}^Z_i\in {\mathbb R}^{n}, i=1, . . ., n$ and ${\bf v}^Z_i \in {\mathbb R}^{2m}, i=1, . . ., 2m$ be the orthonormal left and right column singular vectors of $Z$ correspondingly, let ${\bf U}^Z=\Bigl[{\bf u}^Z_1 \, \vdots  \, . . . \, \vdots \,{\bf u}^Z_{n}\Bigr], {\bf V}^Z=\Bigl[{\bf v}^Z_1 \, \vdots  \,  . . . \, \vdots  \, {\bf v}^Z_{2m}\Bigr]$  and ${\bf \Sigma}^Z \in {\mathbb R}^{n \times 2m}$ be the matrix of singular values $\sigma_1\geq . . .  \sigma_{2m}$. 

Then as in (\ref{Y}), 
\beq
\bes
&{\bf Z}={\bf U}^Z{\bf \Sigma}^Z ({\bf V}^Z)^T,  \text{ with }\operatorname{rank}(Z)= m^Z\\
\end{split}
\eeq
and for $l\leq m^Z$, as in (\ref{erest1}),
\beq
||{\bf y}(t_i)-\sum_{k=1}^{l} \sigma^Z_k v^Z_{ik} {\bf u}^Z_k||_2\leq \sigma^Z_{l+1} \text{  and }
||{\bf f}({\bf y}(t_i),  t_i)-\sum_{k=1}^{l} \sigma^Z_k v^Z_{i+m \, k} {\bf u}^Z_k||_2\leq \sigma^Z_{l+1}.
\eeq

Similarly to the previous section, we denote by $\tilde {\bf U}^Z, \breve {\bf U}^Z, \bar {\bf U}^Z$ the  matrices corresponding to the eigenvectors associated with the first $l$ singular values, the $l+1, . . ., m^Z$ singular values and the nullspace of ${\bf Z}{\bf Z}^T$.

\subsection{Projection estimates}
Let the column vectors ${\bf u}_k\in {\mathbb R}^n, k=1, . . ., n$ be an orthonormal basis in ${\mathbb R}^n$ consisting of the singular vectors of a snapshot matrix.  
In light of the above considerations, let $l \leq m \leq n$ be given integers and let us define the matrices ${\bf \tilde U} = [\tilde{\bf u}_i]= [{\bf u}_i], \text{ for } i=1, . . ., l; \, {\bf \breve U} =[\breve{\bf u}_i]= [{\bf u}_i],  \text{ for } i=l+1, . . ., m; \, {\bf \bar U} = [\bar{\bf u}_i]=[{\bf u}_i],  \text{ for }  i=m+1, . . ., n$ (one or both of the latter may be empty). 

Let us consider the $n\times n$ matrices (again, one or both of the latter may be empty) $$\tilde {\bf P}={\bf \tilde U}\tilde {\bf U}^T, \, \breve {\bf P}={\bf \breve U}\breve {\bf U}^T, \, \bar {\bf P}={\bf \bar U}\bar {\bf U}^T.$$
$\tilde {\bf P}$ is the projection onto $\Span \{\tilde {\bf u}_i\}$ and $\breve {\bf P}$ and $ \bar {\bf P}$ are projection matrices onto $\Span\{\breve {\bf u}_j\}$ and  $\Span\{\bar {\bf u}_i\}$.

We will be using of several occasions the following result which is easy to prove 
\begin{proposition} \label{prop1} 
$\tilde {\bf P} +\breve {\bf P}+ \bar {\bf P}={\bf I}_n$ and $||\tilde {\bf P}||_2=||\breve {\bf P}||_2=||\bar {\bf P}||_2= 1$.
\end{proposition}

\bigskip 

Let us consider the specific projection matrices derived from ${\bf Y}$ and ${\bf Z}$, $\tilde {\bf P}^Y= \tilde {\bf U}^Y (\tilde {\bf U}^Y)^T, \breve {\bf P}^Y= \breve {\bf U}^Y (\breve {\bf U}^Y)^T, \bar {\bf P}^Y= \bar {\bf U}^Y (\bar {\bf U}^Y)^T.$ 
The following estimates of the projections of the snapshots in Section \ref{s1} are derived as in (\ref{erest1}).

\beq
\label{corol1}
\bes
&||\tilde {\bf P}^Y {\bf y}(t_i)||_2= ||\sum_{k=1}^{l} \sigma^Y_{k} v^Y_{ik}\tilde {\bf u}^Y_k||_2 =\sqrt{\sum_{k=1}^l (\sigma^Y_{k} v^Y_{ik})^2 }\leq \sigma^Y_1;\\
&||\breve {\bf P}^Y {\bf y}(t_i)||_2= ||\sum_{k=1}^{m^Y-l} \sigma^Y_{k+l} v^Y_{i \, k+l}\breve {\bf u}^Y_k||_2 = \sqrt{\sum_{k=1}^{m^Y-l}(\sigma^Y_{k+l} v^Y_{i \, k+l})^2} \leq \sigma^Y_{l+1}\\
\end{split}
\eeq
and, 
\beq
\label{corol11}
||\bar {\bf P}^Y {\bf y}(t_i)||_2= 0.
\eeq

Similarly,  if $\tilde {\bf P}^Z= \tilde {\bf U}^Z (\tilde {\bf U}^Z)^T, \breve {\bf P}^Z= \breve {\bf U}^Z (\breve {\bf U}^Z)^T, \bar {\bf P}^Z= \bar {\bf U}^Z (\bar {\bf U}^Z)^T$, 
the following estimates of the projections of the snapshots from the solution and its time derivatives in Section \ref{s2} are derived as in (\ref{erest1}).

\beq
\label{rr}
\bes
&||\tilde {\bf P}^Z {\bf y}(t_i)||_2 \leq \sigma^Z_1 ;\\
&||\tilde {\bf P}^Z {\bf f}({\bf y}(t_i))||_2\leq \sigma^Z_1;\\
&||\breve {\bf P}^Z {\bf y}(t_i)||_2 \leq \sigma^Z_{l+1};\\
&||\breve {\bf P}^Z  {\bf f}({\bf y}(t_i))||_2\leq \sigma^Z_{l+1};\\
&||\bar {\bf P}^Z {\bf y}(t_i)||_2= 0;\\
&||\bar {\bf P}^Z {\bf f}({\bf y}(t_i))||_2= 0,
\end{split}
\eeq

Having derived these estimates, in what follows, we will denote 
\beq \label{not} \breve {\bf P}+\bar {\bf P}=\tilde {\bf P}^{\perp}={\bf I}_n-\tilde {\bf P}.\eeq 

Therefore, 
\beq
\label{rrq1}
||(\tilde {\bf P}^X)^\perp {\bf y}(t_i)||_2 = ||(\breve {\bf P}^X+\bar {\bf P}^X)  {\bf y}(t_i)||_2 = ||\breve {\bf P}^X {\bf y}(t_i)||_2 \leq \sigma^X_{l+1}, \text{ for } X=Y, Z,\\
\eeq
and
\beq
\label{rrq2}
||(\tilde {\bf P}^Z)^\perp {\bf f}({\bf y}(t_i))||_2 = ||(\breve {\bf P}^Z+\bar {\bf P}^Z)  { {\bf f}(\bf y}(t_i))||_2 = ||\breve {\bf P}^Z { {\bf f}(\bf y}(t_i))||_2 \leq \sigma^Z_{l+1}. 
\eeq

\section{Error bounds for the two ROM methods}
\subsection{{\label{deriv}} Derivation of an upper bound in the general case} 
Let ${\bf x}(t)$  be any solution of (\ref{system1}) for some parameter set and initial conditions possibly different from those of $\bf y$. Suppose that a  number $m$ of snapshots from the solution ${\bf y}$ have been calculated as described in the previous sections and a POD basis  ${\bf u}_i\in {\mathbb R}^n, i=1, . . . , n$  has been constructed. Since ${\bf u}_i$ is a basis in ${\mathbb R}^n$, ${\bf x}(t)$ can be decomposed as:
\beq
\label{rep}
{\bf x}(t)={\bf \tilde U} \tilde {\bf x} (t)+{\bf \breve U} \breve{\bf x}(t) +{\bf \bar U} \bar{\bf x}(t),
\eeq
where 
\beq
\tilde {\bf x} =\tilde {\bf U}^T {\bf x} \in \Span\{\tilde {\bf u}_i\} \in {\mathbb R}^l,  \, 
\breve{\bf x}= \breve {\bf U}^T {\bf x} \in \Span\{\breve {\bf u}_i\}\in {\mathbb R}^{m-l}, \,
\bar{\bf x}=\bar {\bf U}^T {\bf x} \in \Span\{\bar {\bf u}_i\}\in {\mathbb R}^{n-m}
\eeq
are solutions of the dynamical system
\beq
\label{abc}
\begin{split}
&{\dot{\tilde{\bf x}}} = \tilde {\bf U}^T {\bf f}(\tilde U \tilde {\bf x} (t)+{\bf \breve U} \breve{\bf x}(t) +{\bf \bar U} \bar{\bf x}(t),  t), \quad  \tilde {\bf x} (0)=\tilde {\bf U}^T {\bf x}(0)\\
&{\dot {\breve{\bf x}}}= \breve {\bf U}^T {\bf f}(\tilde U \tilde {\bf x} (t)+{\bf \breve U} \breve{\bf x}(t) +{\bf \bar U} \bar{\bf x}(t), t), \quad \breve{\bf x}(0)=\breve {\bf U}^T {\bf x}(0)\\
&{\dot {\bar{\bf x}}}=\bar {\bf U}^T {\bf f}(\tilde U \tilde {\bf x} (t)+{\bf \breve U} \breve{\bf x}(t) +{\bf \bar U} \bar{\bf x}(t),  t), \quad  \bar{\bf x}(0)=\bar {\bf U}^T {\bf x}(0).
\end{split}
\eeq

In the POD ROM approach the large $n$-dimensional ODE system (\ref{abc}) is replaced with an $l$-dimensional ODE system of the form
\beq
\label{lower}
{\dot{\bf z} }=\tilde {\bf U}^T {\bf f}({\bf \tilde U} {\bf z} (t), t), {{\bf z} }(0)=\tilde {\bf U}^T {\bf x}(0),
\eeq
The prolongation of ${\bf z} $ in ${\mathbb R}^n$,  ${\bf x}_{\tilde {\bf P}}={\bf \tilde U} {\bf z} $  is the POD  approximation of  the solution ${\bf x}$ of the full system and satisfies the equation (as also stated in \cite{serban}):
\beq
\label{apprsol}
\bes
&\dot{\bf x}_{\tilde {\bf P}}(t)= \tilde {\bf P} {\bf f}({\bf x}_{\tilde {\bf P}}, t)\\ &{\bf x}_{\tilde {\bf P}}(0)= \tilde {\bf P} {\bf x}(0), 
\end{split}
\eeq
while the exact solution of the full system satisfies 
\beq
\bes
\label{np}
&\dot {\bf x}=(\tilde {\bf P}+\tilde {\bf P}^\perp) {\bf f}({{\bf x}}, t), \\
&{\bf x}(0)=(\tilde {\bf P}+\tilde {\bf P}^\perp) {\bf x}(0),.\end{split}
\eeq

Subtracting (\ref{apprsol}) from (\ref{np}), the equation for the error ${\bf e}(t) ={\bf x}(t) - {\bf x}_{\tilde {\bf P}}(t)$ is
\beq
\label{err}
\begin{split}
&\dot {\bf e} = \tilde {\bf P}\bigl[ {\bf f}({\bf x}, t)-{\bf f}({\bf x}_{\tilde {\bf P}},  t)\bigr] + \tilde {\bf P}^{\perp}{\bf f}({\bf x},  t)\\
&{\bf e}(0) =  \tilde {\bf P}^{\perp}{\bf x}(0).
\end{split}
\eeq

We use the  latter system to evaluate the error.  As $\bf f$ is assumed to be continuously differentiable, ${\bf f}({\bf x}, t)-{\bf f}({\bf x}_{\tilde {\bf P}},  t)=\frac{\partial {\bf f}}{\partial {\bf x}}({\bf x_*}(t), t){\bf e}(t)$, where for each $t$, ${\bf x_*}(t)$ is a value defined by Taylor's theorem such that $x_{*,i}(t)$ is in the open  interval with ends $ x_i(t)$ and  $x_{\tilde {\bf P},i}(t)$. 
We denote 
$${\bf A}^*_{\tilde {\bf P}}(t)=\tilde {\bf P}\frac{\partial {\bf f}}{\partial {\bf x}}({\bf x_*}(t), t)$$ 
and note that it is dependent on the projection $\tilde {\bf P}$ also via ${\bf x}_{\tilde {\bf P}}$ through ${\bf x_*}(t)$.  

We note that ${\bf A}^*_{\bar {\bf P}}(t)$ is bounded on $[0,T]$ and thus its 2-norm is bounded. Let $\Lambda_{\tilde {\bf P}}$, be a Lipschitz constant such that
\beq\label{lambda} ||{\bf A}_{\tilde {\bf P}}^*(t)||_2\leq \Lambda_{\tilde {\bf P}}, \forall t\in [0,T].\eeq

We next proceed to  integrate  system (\ref{err}) (the integrals exist because of the boundedness):
\beq
{\bf e}(t) = \tilde {\bf P}^{\perp} {\bf x}(0) +\int_0^t \tilde {\bf P}^{\perp}  {\bf f}({\bf x}, s) ds +\int_0^t  {\bf A}^*(s){\bf e}(s)  ds.
\eeq

Since $\int_0^t \tilde {\bf P}^{\perp}  {\bf f}({\bf x}, s) ds=\int_0^t \tilde {\bf P}^{\perp} \dot {\bf x}(s) ds=  \tilde {\bf P}^{\perp} ( {\bf x}(t)- {\bf x}(0))$, it is easily established that
\beq
\label{err1}
{\bf e}(t) = \tilde {\bf P}^{\perp} {\bf x}(t) +\int_0^t \tilde {\bf P} {\bf A}^*(s){\bf e}(s)  ds.
\eeq

Applying 2-norms to both sides of (\ref{err1}) and  the triangle inequality and noting that $||\int_a^b \phi(s) ds||_2\leq \int_a^b ||\phi(s)||_2 ds $ (\cite{rudin})
we obtain
\beq
\label{errnorm}
||{\bf e}(t)||_2 \leq ||\tilde {\bf P}^{\perp} {\bf x}(t)||_2 +\int_0^t || {\bf A}^*(s) {\bf e}(s)||_2  ds
\eeq

Taking into consideration Proposition \ref{prop1} and (\ref{lambda}), the latter results finally in 
\beq
\label{fin}
||{\bf e}(t)||_2 \leq ||\tilde {\bf P}^{\perp} {\bf x}(t)||_2 + \Lambda_{\tilde {\bf P}} \int_0^t || {\bf e}(s)||_2  ds
\eeq

We will use this last inequality and a formulation of a Gronwall-type theorem from \cite{brunner}, Theorem 1.5.1, stating the following.

\begin{theorem}{If $g, A\geq 0$ are real-valued continuous functions on $[0,T]$ and if the continuous function $\eta$ satisfies $\eta(t)\leq g(t)+\int_0^t A(s)\eta(s) ds, t\in[0,T]$, then 
$$\eta(t)\leq g(t) + \int_0^t A(s)  g(s) \exp\Bigl(\int_s^tA(u) \, du\Bigr) \,ds, \forall t\in [0,T].$$}
\end{theorem}

Let  $$g(t)=||\tilde {\bf P}^{\perp} {\bf x}(t)||_2.$$ 
Clearly $g(t)$ is continuous. Applying Gronwall's inequality, we get the estimate

\beq
\label{errnorm2}
||{\bf e}(t)||_2 \leq g(t) +\int_0^t \Lambda_{\tilde {P}} e^{\Lambda_{\tilde {P}} (t-s)} g(s)  ds \leq  e^{\Lambda_{\tilde P} t}\max_{t\in [0,T]} g(t).
\eeq

The above is a bound for the difference between a solution of system (\ref{system1}) and a solution of the modified system (\ref{apprsol}). In deriving this bound, we have not used so far how the projection matrix $\tilde {\bf P}$ was constructed.  Obviously, the error will be zero if the solution ${\bf x}(t)$ is orthogonal to $\tilde {\bf P}^\perp$, {\it i.e.}  if ${\bf x}(t)$ is in the span of the basis vectors comprising $\tilde {\bf U}$. Clearly, if the snapshots used to define $\tilde {\bf P}$ are chosen so that  ${\bf x}(t)$ is  $\tilde {\bf P}^\perp {\bf x}(t)\approx 0$  over the whole interval $[0,T]$  ({\it e.g.} the solution snapshots are sufficiently  dense and the solution does not change rapidly in time) it is obvious that the error will be small. It also makes sense to consider snapshots from the time derivatives, or, possibly, other characteristics of the solution that contribute to reducing the term $g(t)$. We explore this idea in the next sections. 

\subsection{Error bounds for the approximation of the solution from which the snapshots were collected.}
One could use the reduced basis derived from a particular solution of the FOM (\ref{system1}). It is important to derive error bounds in these cases; however, here we start with deriving an {\it a priori} bound of the error made when approximating the solution $\bf y$ from which the snapshots were collected. We will derive error bounds for the two ROMs: a) when only snapshots from $\bf y$  are used; b) when snapshots from both $\bf f$  and $\bf y$ are used. To this end, we need to make an estimate for $g(t)$ in (\ref{errnorm2}).

Let  $t_i, i=1, . . ., m \leq n, t_1=0, t_m=T$ be the time points at which the $2m$ snapshots  ${\bf y}(t_i)$  and  ${\bf f}({\bf y}(t_i), t_i)$ were calculated.  Denote $\Delta_i=t_{i+1}-t_i.$

\subsubsection{Method 1 -  Using only solution snapshots}
Let the snapshot matrix be ${\bf Y}=[{\bf y}(t_1) \, \vdots \,  . . . \, \vdots \, {\bf y}(t_m)]$, which, given $l\leq m$, generates the projection matrix $\tilde {\bf P}^Y$ and let the  solution of (\ref{apprsol}) be $ {\bf y}_{\tilde {\bf P}^Y}(t)$.

Let us take some value of $t\in [0,T]$ and let it belong to the interval $[t_i, t_{i+1}]$ for some $i$. Applying Lagrange interpolation, we get
\beq
\bes
&{\bf y}(t)={\bf y}(t_i) \frac{(t-t_{i+1})}{(t_i-t_{i+1})}+{\bf y}(t_{i+1}) \frac{(t-t_{i})}{(t_{i+1}-t_{i})} + R({\bf y}(t)), \text{ where }\\
& {\bf R}({\bf y}(t))=\frac{(t-t_i)(t-t_{i+1})}{2}\frac{d^2{\bf y}}{dt^2}(t)|_{t=\zeta}= \frac{(t-t_i)(t-t_{i+1})}{2}\frac{d {\bf f}(t)}{dt}|_{t=\zeta}
\end{split}
\eeq
where $\zeta\in (t_i, t_{i+1})$ and where it is assumed that ${\bf f}$ has first derivatives in all variables \cite{stoer}. Because of Assumption \ref{ass1}, $||\frac{d {\bf f}({\bf y}(t),t)}{dt}||_2$ is bounded on $[t_i, t_{i+1}]$, {\it i.e.} there exists a constant $\Psi_i$, such that  \beq \label{psi} ||\frac{d {\bf f}({\bf y}(t),t)}{dt}||_2 \leq \Psi_i, t\in [t_i, t_{i+1}].\eeq
Multiplying the above  expression for ${\bf y}(t)$ by $(\tilde {\bf P}^Y)^\perp$, we get 
\beq
(\tilde {\bf P}^Y)^\perp {\bf y}(t)= (\tilde {\bf P}^Y)^\perp {\bf y}(t_i)\frac{t-t_{i+1}}{t_i-t_{i+1}}+ (\tilde {\bf P}^Y)^\perp {\bf y}(t_{i+1})\frac{t-t_{i}}{t_{i+1}-t_{i}} + (\tilde {\bf P}^Y)^\perp \frac{(t-t_i)(t-t_{i+1})}{2}\frac{d {\bf f}(t)}{dt}|_{t=\zeta}
\eeq

\noindent Since $\max_{t\in[t_i, t_{i+1}]}|(t-t_i)(t-t_{i+1})|=\frac{\Delta_i}{4}$, it follows
\beq
\label{ppf}
||(\tilde {\bf P}^Y)^\perp {\bf y}(t)||_2\leq || (\tilde {\bf P}^Y)^\perp {\bf y}(t_i)||_2+|| (\tilde {\bf P}^Y)^\perp {\bf y}(t_{i+1})||_2+ \frac{\Delta_1^2}{8}|| (\tilde {\bf P}^Y)^\perp \frac{d {\bf f}(t)}{dt}|_{t=\zeta}||_2
\eeq
Using (\ref{rrq1}) we get from (\ref{ppf}):
\beq
\label{g1}
g(t)=||(\tilde {\bf P}^Y)^\perp {\bf y}(t)||_2\leq 2 \sigma^Y_{l+1} +||(\tilde {\bf P}^Y)^\perp ||_2\frac{\Delta_i^2}{8}\Psi_i. \, 
\eeq
Using (\ref{errnorm2}) we get:
\beq
\label{good1}
||{\bf e}^Y(t)||_2=||{\bf y}(t)-{\bf y}_{\tilde {\bf P}^Y}(t)||_2 \leq  \Bigl[2\sigma^Y_{l+1} + \Psi_i \frac{\Delta_i^2}{8} \Bigr]e^{\Lambda_{\tilde {\bf P}^Y} t}, 
\eeq
where $\sigma^Y_{l+1}$ is the $(l+1)$-th singular value of the snapshot matrix Y  and  $\Psi_i$ and $\Lambda_{\tilde {\bf P}^Y}$ are defined via (\ref{psi}) and (\ref{lambda}) (with $\tilde {\bf P}= \tilde {\bf P}^Y$).

\bigskip 

Thus, we obtain the following 
\begin{proposition}\label{ss2}
Let $\bf f$ satisfy Assumption 1 and let ${\bf y}(t)$ be a solution of (\ref{system1}). Let $t_j, j=1, . . ., m\leq n$ be a set of points in [0,T] and let ${\bf Y}$ be a matrix of solution snapshots, ${\bf Y}=[{\bf y}(t_1) \, \vdots \,  . . . \, \vdots \, {\bf y}(t_m)]$. Let $\tilde {\bf U}^Y$ be the matrix of the  truncated set of the first $l\leq m$ singular vectors ${\bf u}_k^Y, k=1, . . . , l$ of ${\bf Y}$  and $\tilde {\bf P}^Y$ be the corresponding projection matrix; $\tilde {\bf P}^Y=\tilde {\bf U}^Y (\tilde {\bf U}^Y)^T$. Let ${\bf y}_{\tilde {\bf P}^Y}(t)$ be a vector function whose projection in $\Span\{{\bf u}_k^Y, k=1, . . . , l\}$ solves the ROM (\ref{lower}). Then the 2-norm of the error ${\bf e}^Y(t)={\bf y}(t)-{\bf y}_{\tilde {\bf P}^Y}(t)$ satisfies the bound (\ref{good1}), where  $\Delta_i= t_{i+1}-t_i$ and $t\in [t_i, t_{i+1}]$. 
\end{proposition}

\medskip

\subsubsection{Method 2 -  Using snapshots of the solution and the time derivatives}
We use similar logic to the one in the previous section. 

Let the snapshot matrix be $Z=[{\bf y}(t_0) \, \vdots \,  . . . \, \vdots \, {\bf y}(t_m) \, \vdots \,  {\bf f}({\bf y}(t_0), t_0), \, \vdots \, . . . \, \vdots \, {\bf f}({\bf y}(t_m), t_m)]$, which, given $l$, generates the projection matrix $\tilde {\bf P}^Z$ and let the  solution of (\ref{apprsol}) be $ {\bf y}_{\tilde {\bf P}^Z}(t)$.

Let $t\in [0,T]$ belong to the interval $[t_i, t_{i+1}]$ for some $i$. We now apply Hermite interpolation to define a vector  polynomial ${\bf p}^i(t)$ such that ${\bf p}^i(t_k)={\bf y}(t_k), \frac{d{\bf p}^i}{dt}|_{t=t_k}={\bf f}({\bf y}(t_k), t_k), k=i, i+1$.

\beq
\label{hermit}
\bes
{\bf p}^i(t)=&{\bf y}(t_i)\Bigl[1+\frac{2(t-t_i)}{(t_{i+1}-t_{i})}\Bigr]\frac{(t-t_{i+1})^2}{(t_i-t_{i+1})^2}+{\bf y}(t_{i+1})\Bigl[1+\frac{2(t-t_{i+1})}{(t_{i}-t_{i+1})}\Bigr]\frac{(t-t_{i})^2}{(t_i-t_{i+1})^2}+\\
&{\bf f}({\bf y}(t_i), t_i)\frac{(t-t_i)(t-t_{i+1})^2}{(t_i-t_{i+1})^2}+{\bf f}({\bf y}(t_{i+1}), t_{i+1})\frac{(t-t_i)^2(t-t_{i+1})}{(t_i-t_{i+1})^2}
\end{split}
\eeq
and 

\beq
\label{r}
{\bf y}(t)={\bf p}^i(t)+\frac{1}{24}(t-t_i)^2(t-t_{i+1})^2 \frac{d^3 {\bf f}({\bf y}(t), t)}{dt^3}|_{t=\theta}
\eeq
where $\theta\in(t_i, t_{i+1})$ and where it is assumed that ${\bf f}$ is three times differentiable in $y$ and in $t$ on [0,T]. 

\begin{assumption} We further assume that $||\frac{d^3 {\bf f}(t)}{dt^3}||_2$ is continuous and thus, bounded on $[0,T]$. \end{assumption} 
 Thus, for each interval $[t_i, t_{i+1}]$, there exists a constant $\Phi_i$, such that  \beq \label{fi} ||\frac{d^3 {\bf f}(t)}{dt^3}||_2 \leq  \Phi_i, t\in [0,T]\eeq

We apply the following bounds
\beq
\bes
&\max_{t\in[t_i, t_{i+1}]}|2(t-t_i)(t-t_{i+1})^2|=\max_{t\in[t_i, t_{i+1}]}|2(t-t_i)^2(t-t_{i+1})|=\frac{8}{27}\Delta_i^3;\\
&\max_{t\in[t_i, t_{i+1}]}(t-t_i)^2=\max_{t\in[t_i, t_{i+1}]}(t-t_{i+1})^2=\frac{\Delta_i^2}{4},\\
&\max_{t\in[t_i, t_{i+1}]}|(t-t_i)^2(t-t_{i+1})^2|=\frac{\Delta_i^4}{16}
\end{split}
\eeq
and (\ref{rrq1}-\ref{rrq2}) in (\ref{hermit}) and (\ref{r}) to get the bound 
\beq
\label{g2}
\bes
g(t)=&||(\tilde {\bf P}^Z)^\perp {\bf y}(t)||_2\leq [2 \sigma^Z_{l+1} (\frac{1}{4}+\frac{8}{27}+\frac{4}{27}\Delta_i)+ \frac{\Delta_i^4}{16\cdot 24}\Phi_i];\\ 
\end{split}
\eeq
Using this upper bound, we obtain from (\ref{errnorm2}):
\beq
\label{good2}
||{\bf e}^Z(t)||_2 =||{\bf y}(t)-{\bf y}_{\tilde {\bf P}^Z}(t)||_2 \leq  \Bigl[ \sigma^Z_{l+1} (\frac{59}{54}+\frac{4}{27}\Delta_i)+ \frac{\Delta_i^4}{384}\Phi_i \Bigr]e^{\Lambda_{\tilde { P}^Z} t}, 
\eeq
where $\sigma^Z_{l+1}$ is the $l+1$-th g of the snapshot matrix $\bf Z$  and $\Phi_i$ and $\Lambda_{\tilde {\bf P}^Z}$ are defined via (\ref{fi}) and (\ref{lambda}).  $\Lambda_{\tilde {\bf P}^Z}$ is the bound for the Jacobian as defined in section \ref{deriv} and depending on the specific projection ${\bf P}^Z$.

Thus, we obtain the following 
\begin{proposition}
Let $\bf f$ satisfy Assumptions 1,  2 and let ${\bf y}(t)$ be a solution of (\ref{system1}). Let $t_j, j=0, . . ., m\leq n$ be a set of points in [0,T] and let  $Z=[{\bf y}(t_0) \, \vdots \,  . . . \, \vdots \, {\bf y}(t_m) \, \vdots \,  {\bf f}({\bf y}(t_0), t_0), \, \vdots \, . . . \, \vdots \, {\bf f}({\bf y}(t_m), t_m)]$, i.e., $Z$ is  a matrix of solution and time derivative snapshots. Let $\tilde {\bf U}^Z$ be the matrix of the  truncated set of the first $l\leq 2m$ singular vectors ${\bf u}_k^Z, k=1, . . . , l$ of $Z$  and $\tilde {\bf P}^Z$ be the corresponding projection matrix; $\tilde {\bf P}^Z=\tilde {\bf U}^Z (\tilde {\bf U}^Z)^T$. Let ${\bf y}_{\tilde {\bf P}^Z}(t)$ be a vector function whose projection in $\Span\{{\bf u}_k^Y, k=1, . . . , l\}$ solves the ROM (\ref{lower}). Then the 2-norm of the error ${\bf e}^Z(t)={\bf y}(t)-{\bf y}_{\tilde {\bf P}^Z}(t)$ satisfies the bound (\ref{good2}), where  $\Delta_i= t_{i+1}-t_i$ and $t\in [t_i, t_{i+1}]$.

\end{proposition}

\subsection{Remarks on the error bounds} \label{sno}
(a) To the best of our knowledge, the bounds (\ref{good1}) and (\ref{good2}) are the first derived bounds that include the time points at which the snapshots were taken. So far published bounds \cite{hinze, chatur1, iliescu, kv1, kv2, serban}, if including time information at all, assume equidistant snapshots and include the time step. Our bounds do not assume any specific distribution of the snapshot times. For equidistant snapshots $\Delta_i=\Delta$, these  bounds also contain local information by including derivative bounds in the intervals $[t_i, t_{i+1}]$. Including local information may be helpful for rational snapshot selection via error estimates. 

\bigskip

\noindent
(b) We note that the derived bounds are still far from being exact or being fully informative. For a linear system $\dot{\bf x}={\bf Ax}$, 
\beq
{\bf e}(t)=e^{{\bf \tilde P A}t}\Bigl[\tilde {\bf P}^\perp {\bf y}(0)+\int_0^t e^{- {\bf \tilde P A}} \tilde {\bf P}^\perp {\bf A} e^{{\bf A}s}  \tilde {\bf P}^\perp {\bf y}(0) ds\Bigr]
\eeq
{\it i.e.} the error depends on the eigenvalues of $ {\bf \tilde P A}$, which is neither reflected in the bounds derived, nor, in fact in any bounds or estimates published. 
Yet, the derived bounds are valuable because they provide insight about the relationship  between the basis truncation and the location of snapshots.  
 
\bigskip

\noindent
(c) For  linear systems of the form $\dot{\bf x}={\bf A}{\bf x}+{\bf b}, $ an explicit bound, not depending on ${\bf \tilde P}$, can be found. Indeed:
\beq 
||{\bf A}_{\tilde {P}}^*(t)||_2=||{\bf \tilde P A}||_2= \sigma_1({\bf \tilde PA}) \leq ||{\bf A}||_2= \sigma_1({\bf A}),
\eeq 
where $\sigma_1({\bf A})$ is  the largest singular value of {\bf A}, {\it i.e.} we can assume that  $\Lambda_{\tilde {P}}=\sigma_1({\bf A})$. Also, let $\theta_i=  \max_{t\in [t_i,t_{i+1}]}||{\bf y}(t)||_2$.
Then 
\beq \begin{split} 
&\Psi_i=\max_{t\in [t_i, t_{i+1}]}||{\bf A y}(t)||_2\leq ||{\bf A}||_2 \theta_i= \sigma_1({\bf A}) \theta_i,\\  
&\Phi_i=\max_{t\in [t_i,t_{i+1}]}||{\bf A^3 y}(t)||_2\leq ||{\bf A}||_2^3 \theta_i= (\sigma_1({\bf A}))^3 \theta_i
\end{split}
\eeq
and from (\ref{good1}) and (\ref{good2}):
\beq
\label{good11}
||{\bf e}^Y(t)||_2\leq \Bigl[2\sigma^Y_{l+1} + \sigma_1({\bf A}) \frac{\Delta_i^2}{8} \theta_i \Bigr]e^{\sigma_1({\bf A})t}
\eeq
and
\beq
\label{good21}
||{\bf e}^Z(t)||_2  \leq  \Bigl[\sigma^Z_{l+1} (\frac{118}{108}+\frac{4}{27}\Delta_i)+ (\sigma_1({\bf A}))^3 \frac{\Delta_i^4}{384} \theta_i \Bigr]e^{\sigma_1({\bf A}) t}.
\eeq

We next try to compare the bounds for the two ROMs considered. Above, the only quantities that depend on the method used are $\sigma^Y_{l+1} $ and $\sigma^Z_{l+1} $. Obviously, if the dimension of the RBS is taken to be equal to the number of snapshots used (i.e. $l=m$ in the first case and $l=2m$ in the second case and then $\sigma_{l+1}^Y=\sigma_{l+1}^Z=0$),  the reduced model by Method 2 will have twice larger dimension than from Method 1, but, for small $\Delta_i$,  its error might be much  smaller (possibly with two orders in $\Delta_i$ as predicted by the error bounds). Thus, for small $\sigma_{l+1}$ and $\Delta_i$, Method 2 may lead to a significantly smaller error.  Note that the terms "possibly" and "may" are used because these are upper bounds and not exact estimates. We check these predictions in the next section.

In the case where $\sigma_{l+1}$ is not sufficiently small, using both solution and derivative snapshots yields an error bound that is first order in $\Delta_i$, whereas using solution snapshots only yields an error bound that is zero-th order in $\Delta_i$. So improvements in the error for the solution snapshot - only case come primarily from increasing the dimension of the RBS, whereas the error in the the time derivative and solution snapshot case decreases when snapshots are sampled more frequently.  

These derivations  demonstrate the trade-off between error and dimension of the RBS for the two methods  and call for investigating the comparison between the distributions of the singular values  for the two methods. 
For example, for a fixed $l$ and equidistant snapshots, when $\Delta_i=\Delta$ is decreased,  $\sigma^Y_{l+1}$ and $\sigma^Z_{l+1}$ may increase (since adding more snapshots means adding more singular values).  Therefore, there may exist an optimal value of $\Delta$ such that further decrease of the distance between snapshots would not decrease the error as it will be dominated by the error caused by the  truncation of the snapshot-generated basis. We demonstrate this phenomenon in the next section. 

Now let us look at the case when the dimension $l$ of the RBS is fixed for both methods at $l=m$.  Then $\sigma^Y_{l+1}=0$ but $\sigma^Z_{l+1}>0$.  However, $\sigma^Y_{l+1}$ may be insignificant compared to the rest of the error.
If the bounds were good estimates of the error, the error of Method 2 would be smaller than the error of Method 1 in such cases. This possibility is explored and demonstrated on examples presented in the next section. 

Arguments similar to the ones presented above for linear systems  hold also for nonlinear systems for which the existence of a constant $\Lambda$ not dependent on $\tilde P$, such that $\Lambda\geq \Lambda_{\tilde {P}}$ can
 be proved. 

\bigskip

(d) As pointed above, investigating the distribution of the singular values and how it changes with decreasing $\Delta$ is important. It is not clear how dependent on the particular problem this distribution is. For the purpose of comparing methods using snapshots containing additional information about the problem (such as derivatives), it is important to understand how adding this information may change the distribution of singular values.  Some examples demonstrating these differences are considered in the next section.

\section{Numerical experiments}  
To validate the above bounds and to compare numerically the error from the two methods, 
we performed numerical experiments with systems of ODEs derived from a method-of-lines discretization of the FitzHugh-Nagumo system with diffusion (FHND)  \cite{keener}. This system has been used as test problem in various studies of model reduction methods, {\it e.g.} \cite{chatur1}. To provide some background, the FHND system is an approximation of the Hodgkin-Huxley system of equations, designed to describe the propagation and dynamics of an action potential (difference in external and intracellular voltage) generated along the nerve axon.   Since the dynamic behavior of the solutions to the FHN system is very sensitive to changes in some of the parameter values \cite{kostova}, and the solutions are characterized by a combination of fast and slow dynamics, it is often used as a test case for the accuracy of numerical methods. In the FHND system,  $V(x,t)$ is the membrane potential, {\it i.e.}, the difference between the extracellular and intracellular potentials and $w(t)$ is a "recovery variable".

Different texts ({\it e.g.}, \cite{olmos},  \cite{rinzel} , \cite{keener}) consider different forms of FHND. In general, the 1D version of the system is a system of two reaction-diffusion equations 
\beq
\bes
\label{eqform}
&v_t=D_1 v_{xx} + f_1(v,w)\\
&w_t=D_2 w_{xx}+ f_2(v,w)\\
\end{split}
\eeq
where  $x\in [0, X]$, representing a one dimensional axon with length $X$. 
In most texts, $f_1$ is cubic in $v$ and linear in $w$:  $f_1=\lambda \Bigl[ v(1-v)(v-a)-w\Bigr]$ and $f_2=cv-bw$. 

The initial conditions are set up so that initially the nerve membrane is at equilibrium: 

\beq \label{inifhn}
v(x,0)=0; w(x,0)=0.
\eeq 

Different boundary conditions (BC), depending on the problem, can be considered as outlined in \cite{keener}. The Neumann BCs $v_x(0,t)= -  {I}_0(t), v_x(X,t)=0$ correspond to applying current ${ I}_0$ at the "left" (at 0) end of the axon (where $r$ is a constant, depending on the internal and external resistances) and "sealing" the axon at the other end (no current).  

\bigskip 

In this paper, we solve numerically equations (\ref{eqform}) together with the following boundary conditions for $w(x,t)$:
\beq
\bes
\label{BCw}
&w(0,t)=w_0(t)\\ 
&w(X,T)=w_X(t),\\
\end{split}
\end{equation}

and  the BC for $v(x,t)$:

\beq
\bes
\label{BCa}
&v_x(0,t)= -  {I}_0(t)\\
&v_x(X,t)= -  {I}_X(t).\\
\end{split}
\eeq

Above, $w_0(t), w_X(t),  {I}_0(t),  {I}_X(t)$ are input functions.

Equations (\ref{eqform} -- \ref{BCw} -- \ref{BCa}) were semidiscretized using the method of lines with finite differences. 
We  denote: 
 $n= 2(L+1); \Delta x= X/L; x_j=j \Delta x, j=0, . . . , L;  v_j (t)=v(x_j, t), j=0, . . ., L $

Using the approximations 
\beq
\bes
&v_{xx}(0,t)\approx \frac{ \displaystyle{\frac{v_2(t)-v_1(t)}{\Delta x}} - v_x(0,t)}{\Delta x}; \, v_{xx}(x_j,t)\approx  \frac{v_{j-1}+v_{j+1}-2v_j}{(\Delta x)^2}, j>0; \\ & v_{xx}\approx \frac{v_x(X,t)- \displaystyle{\frac{v_L-v_{L-1}}{\Delta x} }}{\Delta x},
\end{split}
\eeq
we derive discretizations by the method of lines for the equations (\ref{eqform}) - (\ref{BCw}, \ref{BCa}) as follows. 

\beq
\label{fhn1}
\bes
&\frac{dv_0}{dt}=\frac{D_1}{(\Delta x)^2} [v_2-v_1 +\Delta x   {I}_0(t)]+f_1(v_0, w_0),\\
&\frac{dv_j}{dt}=\frac{D_1}{(\Delta x)^2} [v_{j+1}-2v_j+v_{j-1}]+ f_1(v_j, w_j), j=1, . . ., L-1\\
&\frac{dv_L}{dt}=\frac{D_1}{(\Delta x)^2}[v_{L-2}-v_{L-1}- \Delta x  {I}_1(t)]+ f_1(v_L, w_L)\\
&\frac{dw_j}{dt}=\frac{D_2}{(\Delta x)^2} [w_{j+1}-2w_{j}+w_{j-1}] + f_2(v_j, w_j), \, j=1, . . ., L-1\\
\end{split}
\eeq
where the initial conditions are $v_i(0)=0, \, w_i(0)=0$ and $w_0, w_L$  are defined in (\ref{BCw}).

\subsection{Experiments} 
The two POD ROM methods described above were implemented in a Matlab code. 
Specifically, we solved equations (\ref{BCw}), (\ref{BCa}), (\ref{fhn1}) and selected  $m$ equally spaced on the time interval $[0,T]$  snapshots from the solution.  The ODE solutions were obtained numerically using Matlab's routine $\texttt{odes15s}$, where the absolute and relative tolerances were sufficiently small (usually equal to $10^{-12}, 10^{-14}$) to ensure stable performance of the integrator.  The right hand side of (\ref{fhn1}) was calculated using the selected snapshot vectors to obtain the time derivative snapshots. The snapshot matrices $\bf Y$ and $\bf Z$ were formed and their SVD calculated using Matlab's $\texttt {svd}$ routine. The dimension $l$ of the ROM was either predefined or chosen so that $\sigma_{l+1}<\varepsilon\leq \sigma_{l}$ for a predefined $\varepsilon$. After determining the value of  $l$,  the respective  projection matrices for the two methods, $\tilde {\bf P}^Y$ and $\tilde {\bf P}^Z$ were calculated as described. The solutions of (\ref{apprsol}) with $\tilde {\bf P}=\tilde {\bf P}^Y$ and $\tilde {\bf P}=\tilde {\bf P}^Z$ were calculated respectively for the two models.  
 The distribution of the singular values in the two methods for different snapshot selections was also calculated. 

To be more specific, let us denote by ${\bf y}^Y=({\bf v}^Y, {\bf w}^Y)$ the  solution obtained by using only solution snapshots and by ${\bf y}^Z=({\bf v}^Z, {\bf w}^Z)$  the solution obtained by using both solution and derivative snapshots. We calculate  $||{\bf e}^Y(t)||_2 = ||{\bf y}(t)-{\bf y}^Y(t)||_2$ and $||{\bf e}^Z(t)||_2 =||{\bf y}(t)-{\bf y}^Z(t)||_2$. 

All experiments were done with $X=10, w_0=w_l=0$.

\subsubsection{Experiment A}
The two ROM methods were explored for a linear homogeneous system, obtained from  (\ref{fhn1}) with  $\lambda=0$, of 402 ({\it i.e.},  $L=200$) equations with constant coefficients $\Delta x= 0.05, D_1=15, D_2=10,  \mu=10, \gamma=5$ and ${I}_0(t)=1, {I}_1(t) =5$. ODE solutions were computed in Matlab using $\texttt{ode15s}$  with both absolute and relative tolerances set to $10^{-14}$. The dynamics of the system can be predicted theoretically. It has no equilibria, the eigenvalues are nonpositive and there are 2 zero eigenvalues. Due to the zero initial conditions, the solutions $v_j, w_j$ behave asymptotically like $te^{-gt}$, where $g>0$ is a constant, so the solutions eventually converge to 0.

The system was solved on the time interval $[0, 0.5]$. The calculated FOM solution is shown on Figure S1 (Supplement).

For both methods, three sets of equidistant snapshots were collected, with $\Delta$ = 0.01 (50 snapshots), 0.005 (100 snapshots), 0.0025 (200 snapshots) and the errors $||{\bf e}^Y(t)||_2 , ||{\bf e}^Z(t)||_2$ were calculated. The plots in Figure 1 demonstrate the size of these errors when $\Delta$ is reduced while keeping the first neglected singular value $\sigma_{l+1}$ constant, to compare with the predicted behavior of $||{\bf e}^Y(t)||_2$ and $||{\bf e}^Z(t)||_2$. 

The top two plots in Figure 1 correspond to $\sigma_{l+1}\leq 10^{-15}$. If the bounds (\ref{good1}) and (\ref{good2}) were exact predictors of the error, the distances between the curves for these cases  would be roughly equal to $\log_{10}||{\bf e}_{\Delta}^Y(t)||_2 - \log_{10}||{\bf e}_{\Delta/2}^Y(t)||_2= \log_{10}\frac{||{\bf e}_{\Delta}^Y(t)||_2}{||{\bf e}_{\Delta/2}^Y(t)||}\approx \log_{10} \frac{\Delta^2}{(\Delta/2)^2}=\log_{10}(2^2) \approx 0.602$ and  $\log_{10}||{\bf e}_{\Delta}^Z(t)||_2 - \log_{10}||{\bf e}_{\Delta/2}^Z(t)||_2  \approx \log_{10}(2^4)  \approx 1.204$, respectively (in the latter the sub-index $\Delta$ was added to denote the value of the spacing between snapshots used in the calculations). The actual distances appear to be twice as big, indicating that the error is probably of higher order in $\Delta$ than the error bounds. 

For larger $\sigma_{l+1}$ the distances between the curves decrease below these values (0.602 and 1.204). This behavior would be expected if the bounds were exact estimates of the error. Indeed,  if this was the case  $\log_{10}||{\bf e}_{\Delta}^Y(t)||_2 - \log_{10}||{\bf e}_{\Delta/2}^Y(t)||_2= \log_{10}\frac{||{\bf e}_{\Delta}^Y(t)||_2}{||{\bf e}_{\Delta/2}^Y(t)||} \approx  \log_{10}\frac{2\sigma_{l+1}+\Psi_i\Delta^2}{2\sigma_{l+1}+\Psi_i\Delta^2/4}$ which is decreasing function of $\sigma_{l+1}$. This behavior is  demonstrated on the  middle two ($\sigma_{l+1}\leq  \varepsilon = 10^{-9}$) and bottom two ($\sigma_{l+1} \leq  \varepsilon = 10^{-1}$)  plots. Increasing the value of  $\sigma_{l+1}$ in general decreases the distance between the curves corresponding to different decreasing values $\Delta$ and increases the error in all three cases since it is dominated by the value of $\sigma_{l+1}$.  Note that the error corresponding to smaller $\Delta$ is not always smaller (middle and bottom right plots) when the value of $\sigma_{l+1}$  is significant. 

Often in practice, the dimension of the ROM is predefined. Therefore, it is interesting to compare the performance of the two methods with fixed ROM basis dimension.  We present the results of an experiment with the same linear system where we compare the error from the two methods with the three different time steps and with a fixed dimension of the RBS in Figure 2. For the 402-variable FOM, equidistant snapshots in time of both the FOM solution and its time derivative were collected at spacings of $\Delta =0.01, 0.0025, 0.005$. After calculating the respective singular vectors, the first  $l=5, 10, . . . , 50$ vectors were used to calculate the ROM solution. Presented are 6 of the calculations, with ROM basis dimensions 5, 10, 15, 20, 35, 50 to illustrate the observed tendency (Figure S2). The plots in red in Figure 2 correspond to the three ROMs via Method 1 and the plots in blue correspond correspond to ROMs calculated by Method 2. 

To get understanding of the plots presented on Figure 2, we investigate the distributions of the base-10 logarithm of the singular values in the two ROMs  (Figure S3). 
 
For  $l=5$, all three types of snapshot spacings in Method 2 produce extremely large error ($>10^2$), while the error produced by Method 1 is of size ($10^{-1}$) comparable with the amplitude of the solution ($0-0.5$).  Concomitantly, the first neglected singular value,  $\sigma_6^Y$, for all 3 snapshot selections in Method 1 is less than than 1, while $\sigma_6^Z$, for all 3 snapshot selections in Method 2 is at least 10 times larger (Figure S2). In this case, for both methods, decreasing the distance between snapshots has insignificant effect on the value of the error, which is evidently dominated by the neglected singular value.  

When the ROM basis dimension is  increased to $l=10$, $||{\bf e}^Y(t)||_2$ and  $||{\bf e}^Z(t)||_2$ are comparable. For $l \geq15$, the error produced by Method 2 is considerably smaller than from Method 1. Having in mind the derived bounds, we explain this behavior with the steep decrease of $\sigma_{l+1}^Y$ and $\sigma_{l+1}^Z$ and consequent dominance of the $O(\Delta^2)$ and  $O(\Delta^4)$ in Methods 1 and 2, respectively. At $l=15$ the neglected singular value is less than $10^{-1}$, at $l=25$ it is less than $10^{-5}$, and at  $l=35$ it is less than $10^{-10}$ for all between-snapshot distances for both Method 1 and Method 2 (Figure S2). For this problem, for ROM dimension higher than 15, the error is mostly due to the terms that do not depend on the truncation of the basis, but depend on the error contributed by the spacing of the snapshots; therefore the method with higher order (O($\Delta^{4}$)) produces smaller error.

\begin{figure}
\label{A2}
\epsfig{file=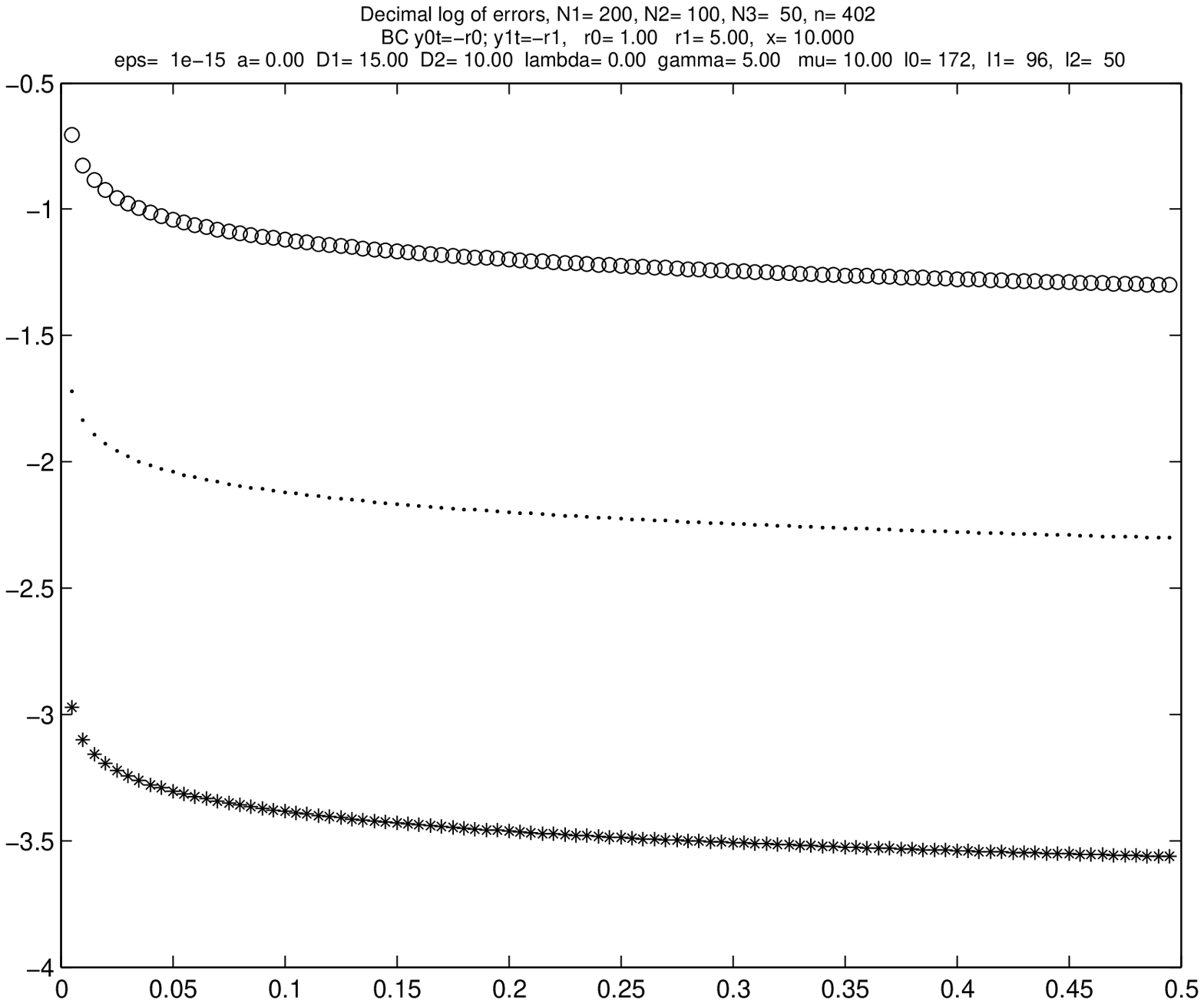, height=2.7in,width=3.2in}
\epsfig{file=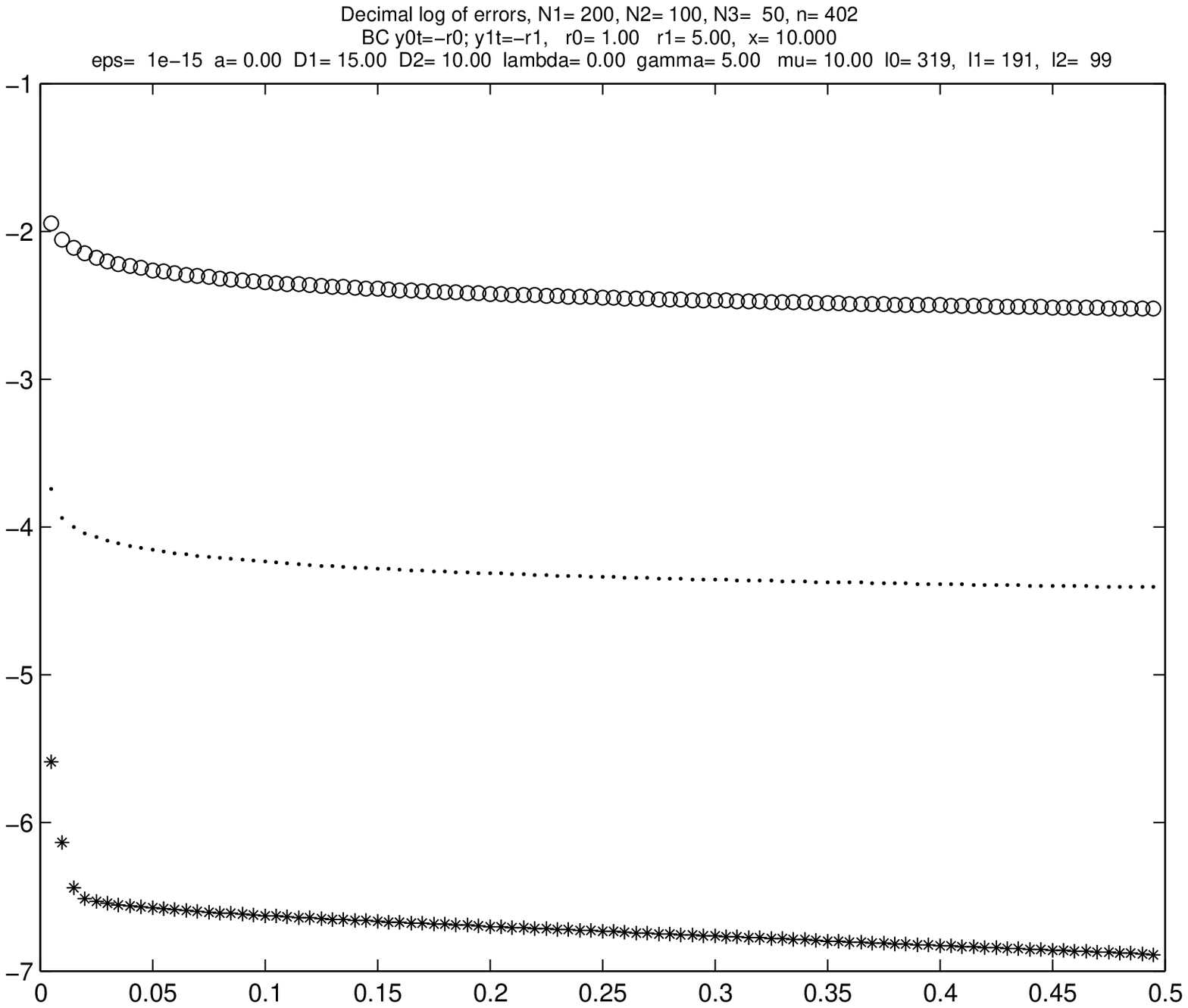, height=2.7in,width=3.2in}\\
\epsfig{file=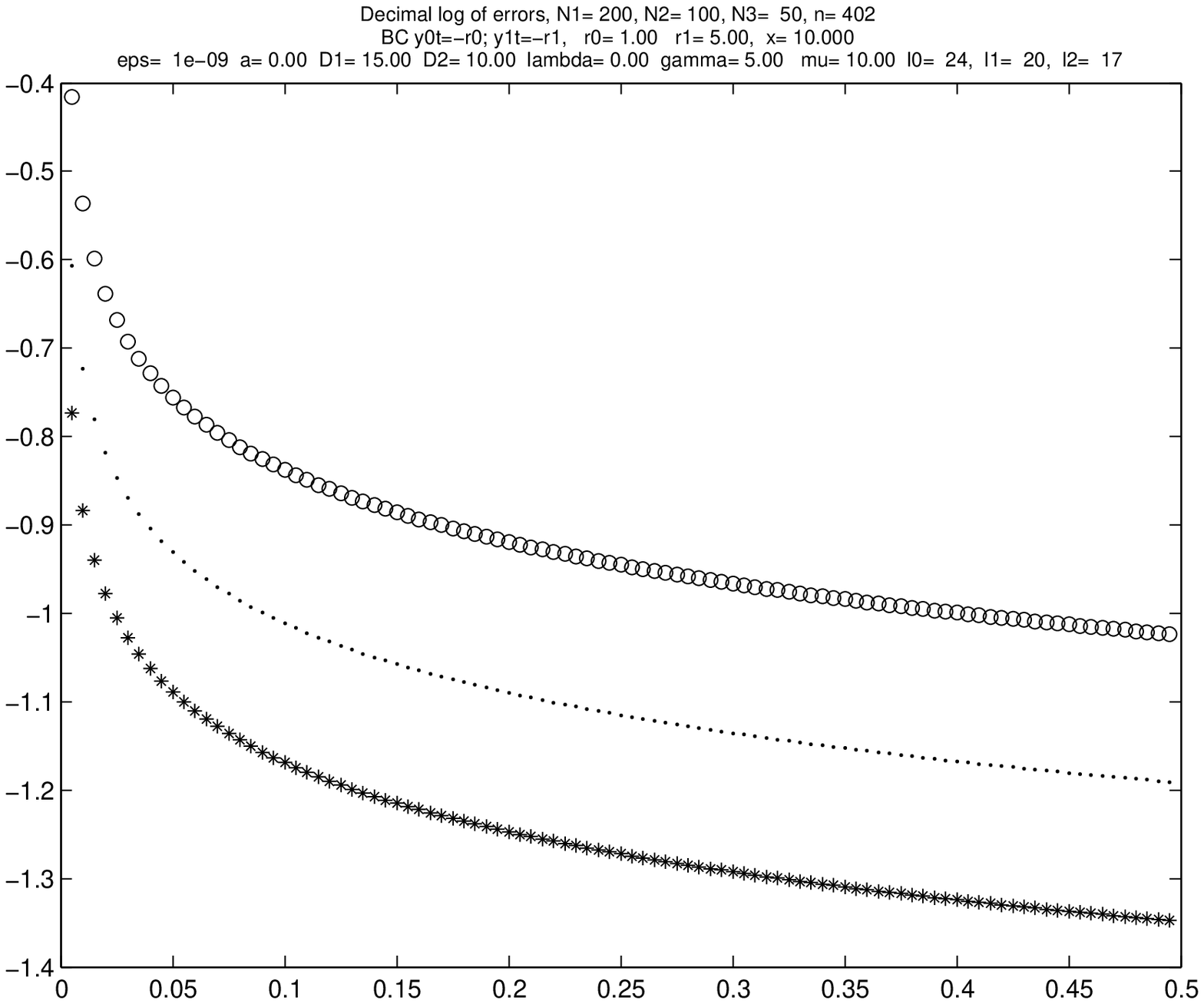, height=2.7in,width=3.2in}
\epsfig{file=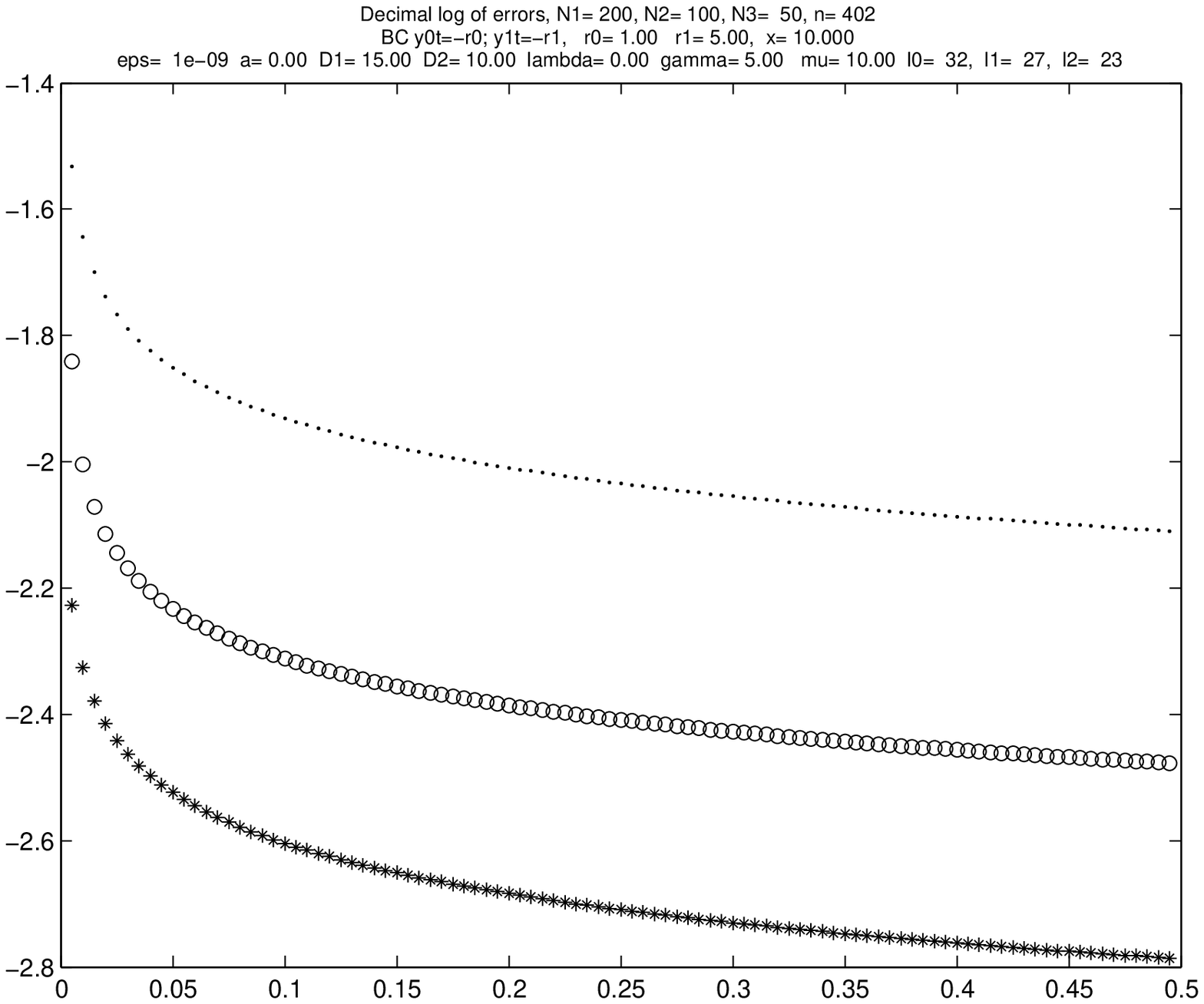, height=2.7in,width=3.2in}\\
\epsfig{file=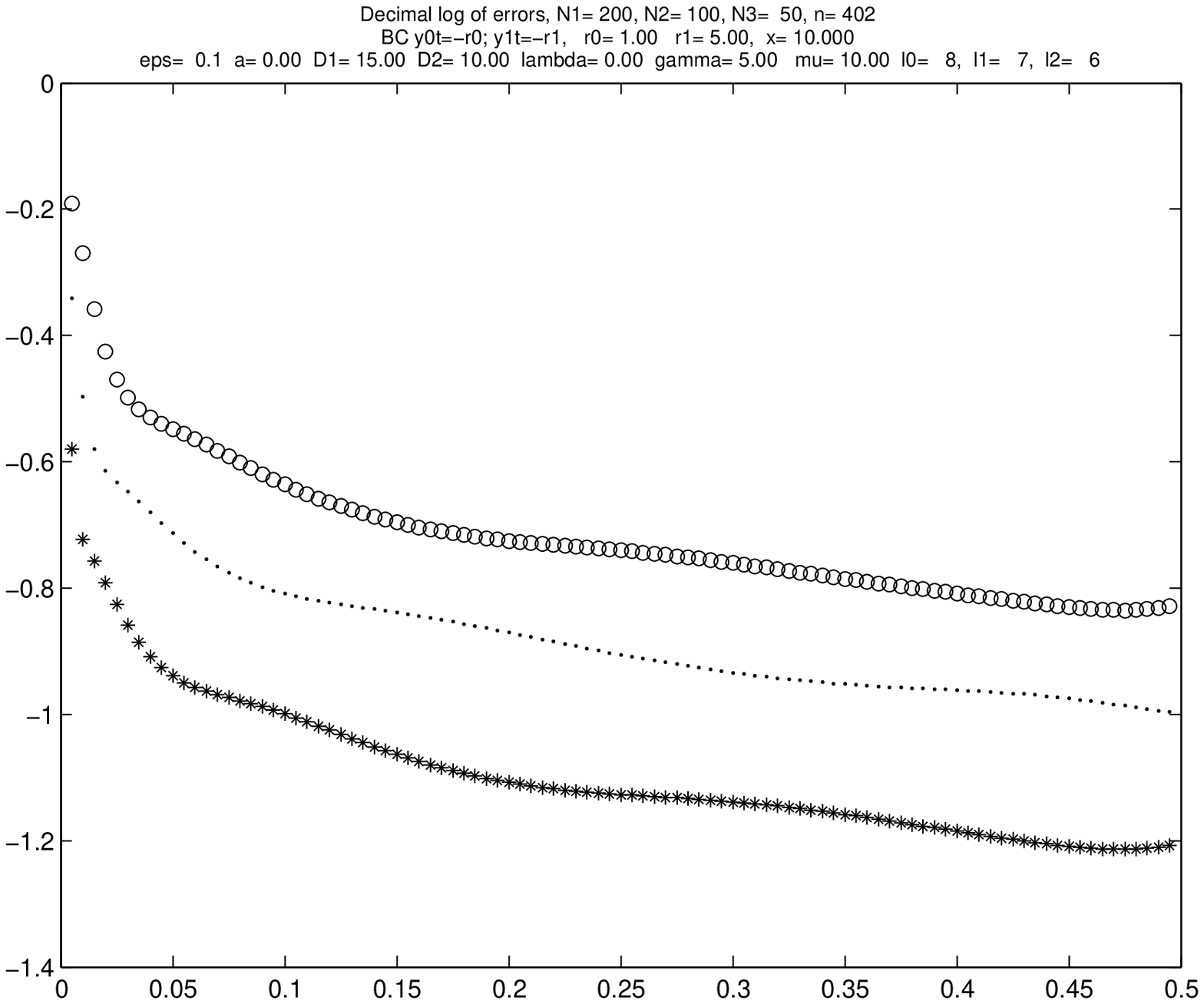, height=2.7in,width=3.2in}
\epsfig{file=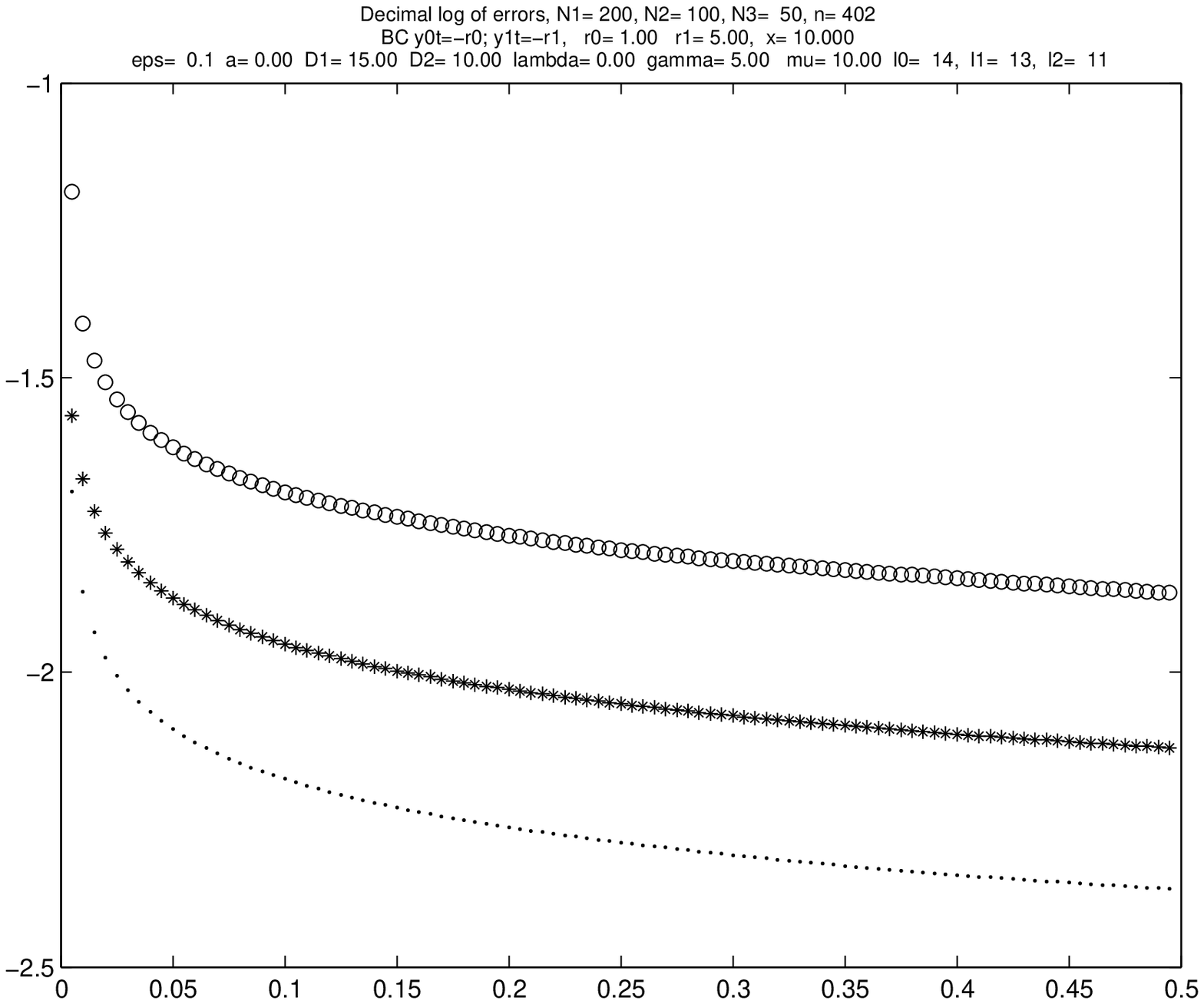, height=2.7in,width=3.2in}\\
\caption{Experiment A.  Error from the two methods at three different values of $\Delta =0.01, 0.0025, 0.005$ and different cutoff ($\varepsilon$) values. The x-axis is time $t$ and the y-axis is $\log_{10} (||{\bf e}^Y(t)||_2)$ (left) and $\log_{10} (||{\bf e}^Z(t)||_2)$ (right).  Circles correspond to $\Delta t=0.01$, dots - to $\Delta t=0.005$, and crosses to $\Delta t=0.0025$. The plots on the right correspond to error from Method 2 (solution and derivative snapshots) and plots on the left correspond to error from Method 1 (no derivative snapshots). }
\end{figure}

\begin{figure}[H]
\label{A3}
\epsfig{file=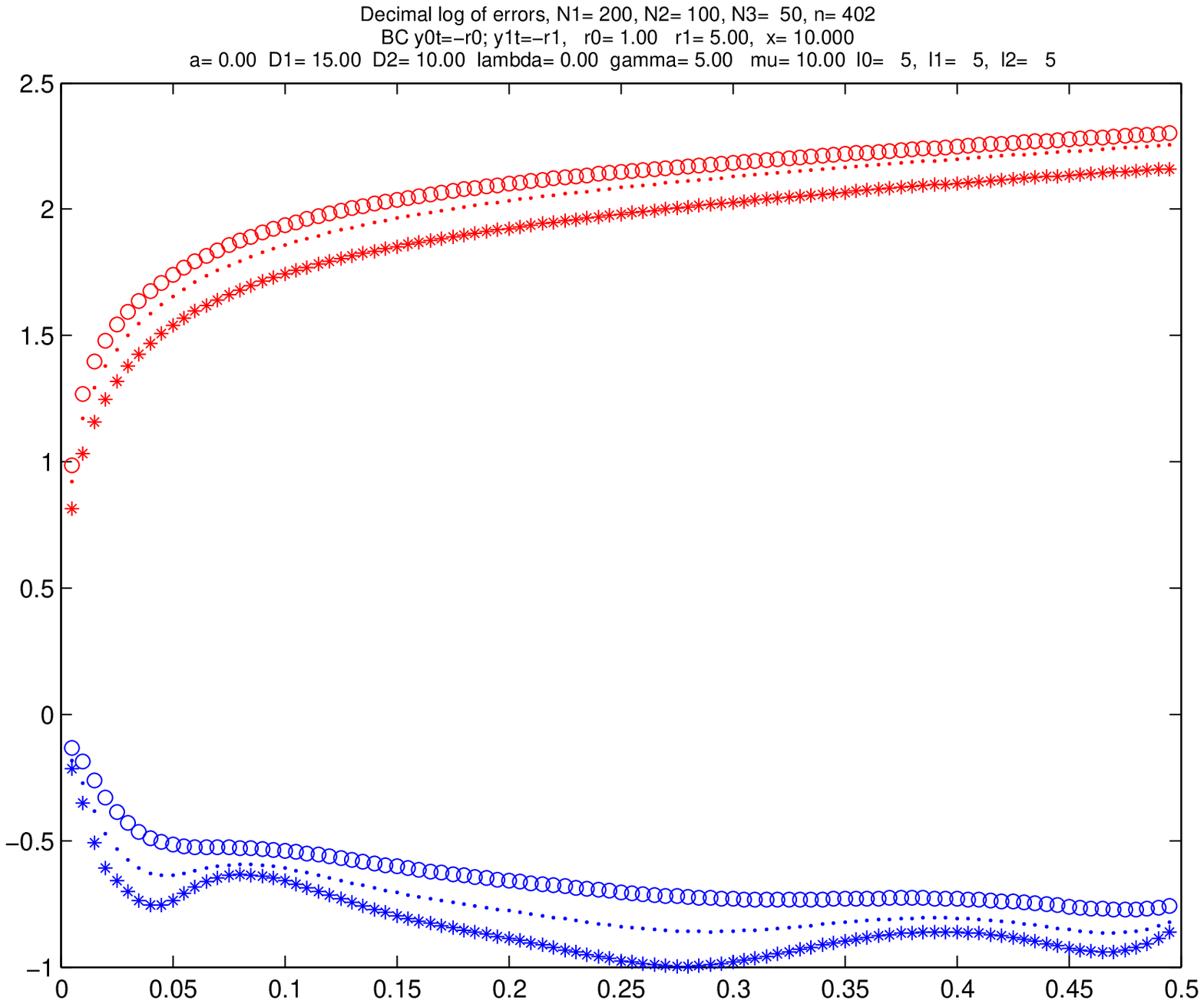, height=2.7in,width=3.2in}
\epsfig{file=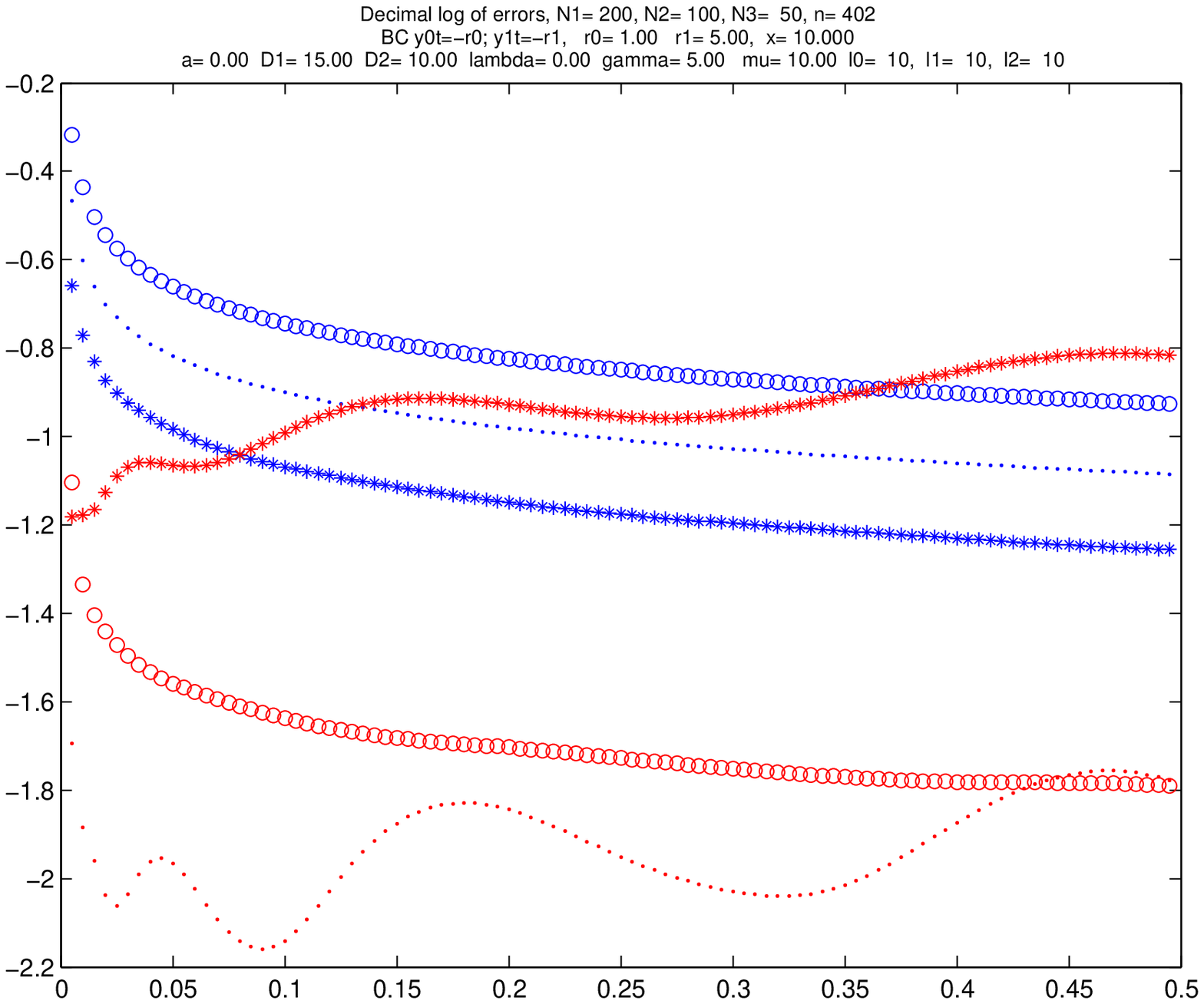, height=2.7in,width=3.2in}\\
\epsfig{file=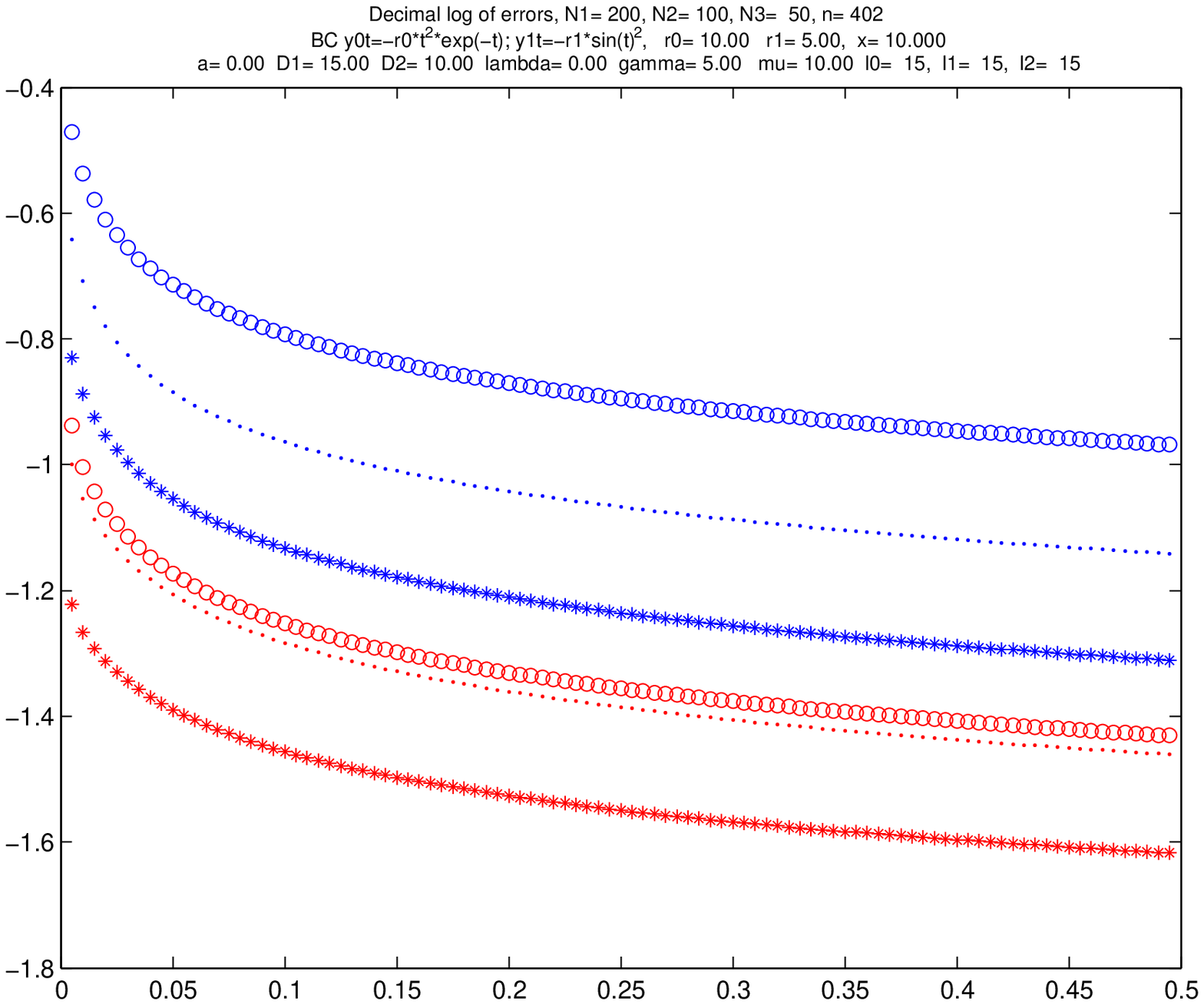, height=2.7in,width=3.2in}
\epsfig{file=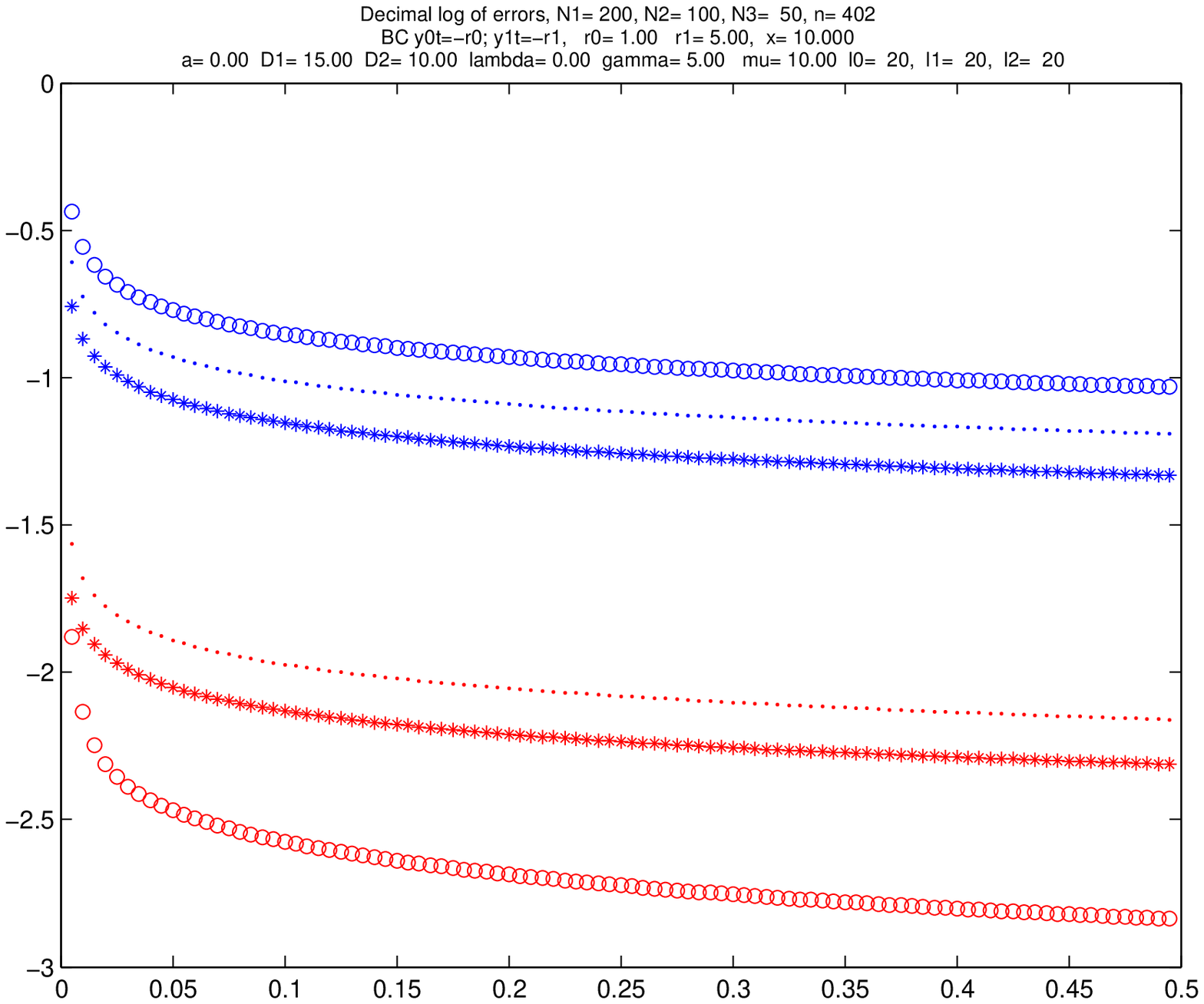, height=2.7in,width=3.2in}\\
\epsfig{file=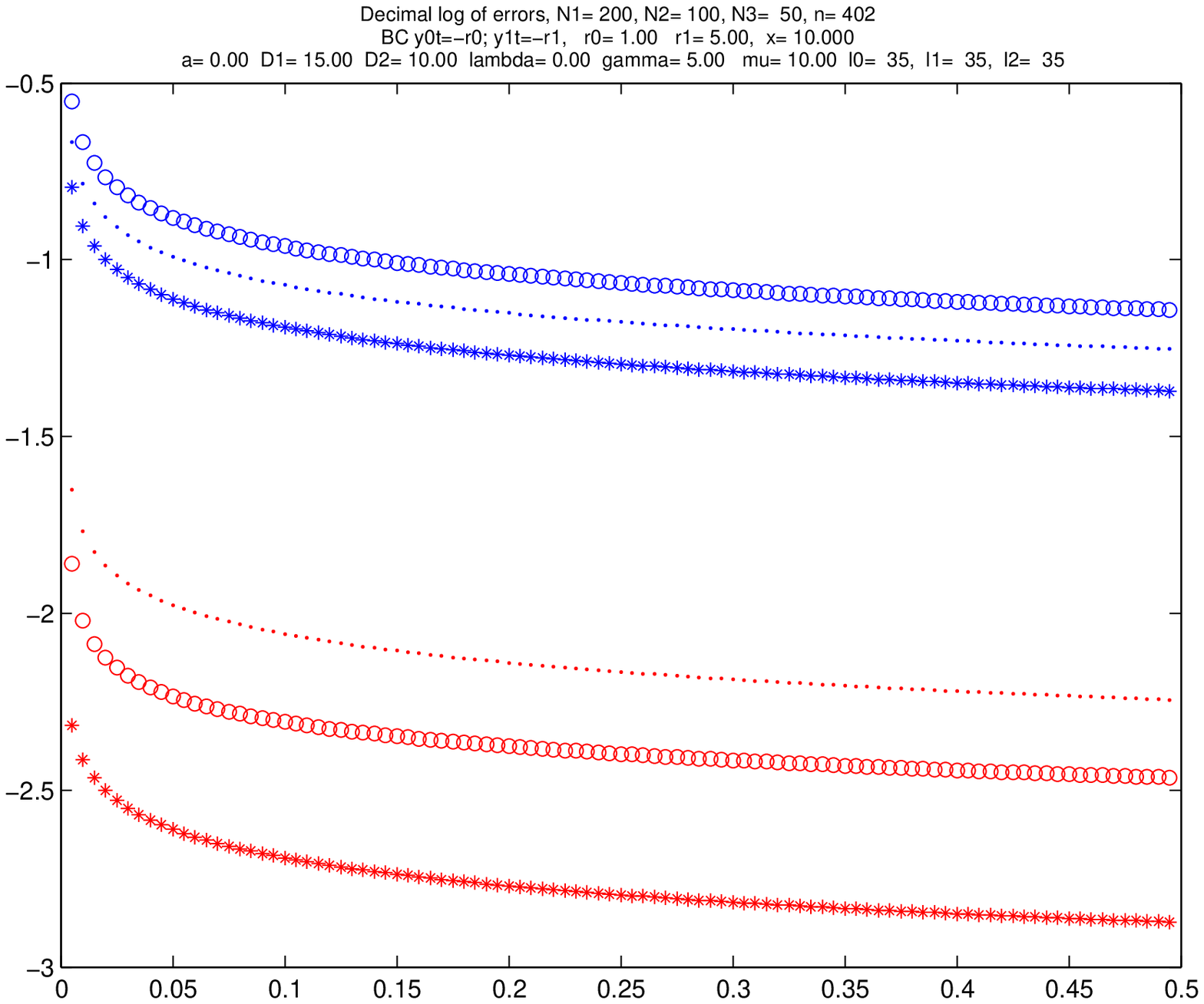, height=2.7in,width=3.2in}
\epsfig{file=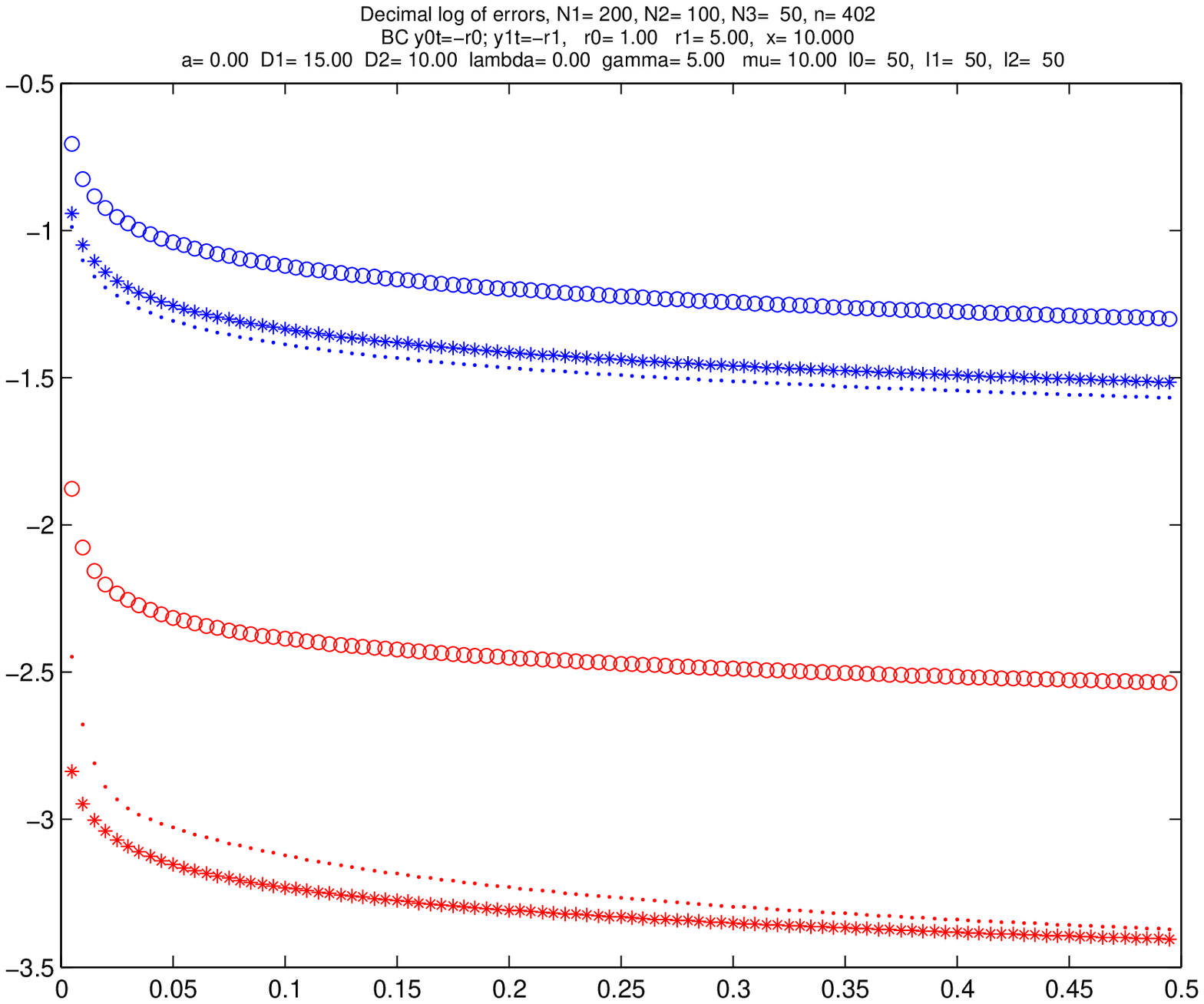, height=2.7in,width=3.2in}\\
\caption{Experiment A.  Error from the two methods at three different values of $\Delta =0.01, 0.0025, 0.005$ and different fixed RBS dimensions ($l=5, 10, 35, 50$). The x-axis is time $t$ and the y-axis is $\log_{10} (||{\bf e}^Y(t)||_2)$ (left) and $\log_{10} (||{\bf e}^Z(t)||_2)$ (right).  Circles correspond to $\Delta =0.01$, dots - to $\Delta =0.005$ and crosses to $\Delta =0.0025$. Red - Method 2, blue - Method 1. The y-axis is the decimal logarithm of the error.}
\end{figure}

We note again that we have only derived upper bounds and not exact estimates of the error. The phenomenon demonstrated above may not necessarily be valid in all cases. However, Experiment A demonstrates the potential for achieving better accuracy of approximation when using Method 2 instead of  Method 1 for relatively low-dimensional ROMs. 

\subsubsection{Experiment B}
In this experiment, the two ROM methods were explored for the nonlinear system  (\ref{fhn1}) with  $\lambda=2$, of 402 ({\it i.e.}  $L=200$) equations with constant coefficients $\Delta x= 0.05, D_1=5, D_2=1, \mu=1, \gamma=5, a=0.1$ and ${ I}_0(t)=1.5 (\sin t)^2, {I}_1(t) =0.5 (\sin t)^2$. Again, ODE solutions were computed in Matlab using $\texttt{ode15s}$  with both absolute and relative tolerances set to $10^{-14}$.

For the convenience of the reader, some of the text is identical with the previous section. 

 The system was solved on the time interval $[0, 2]$. The calculated FOM solution is shown on Figure S4 (Supplement). The figure shows the plots of the 402 variables $y_j(t)$ where time $t$ is on the x-axis. 

As in the previous example, for both methods, three sets of equidistant snapshots were collected, with $\Delta = 0.04, 0.02, 0.01$ and the errors from Method 1 and Method 2, $||{\bf e}^Y(t)||_2 , ||{\bf e}^Z(t)||_2$, were calculated. The plots on Figure 3 demonstrate the size of these errors when $\Delta$ is reduced, while keeping the first neglected singular value $\sigma_{l+1}$ below a constant threshold, to compare with the predicted, by the derived  bounds,  behavior of $||{\bf e}^Y(t)||_2$ and $||{\bf e}^Z(t)||_2$. 

The top two plots in Figure 3 correspond to $\sigma_{l+1}\leq 10^{-15}$. Even though the problem is nonlinear, the distances between the error curves suggest, similarly to the previous numerical experiment, that the error is of higher order in $\Delta$ than the error bounds indicate. The error incurred by applying a ROM with derivative snapshots is again at least one order smaller than the error incurred by Method 1.

The middle and bottom two plots demonstrate the dependence of the error on the size of the largest neglected singular value  ($\sigma_{l+1}\leq \varepsilon =10^{-7}$ and $\sigma_{l+1}\leq  \varepsilon = 10^{-4}$)  respectively. Similarly to the previous numerical experiment (A),  increasing the value of  $\sigma_{l+1}$ generally decreases the distance between the curves corresponding to different decreasing values of $\Delta$ and increases the error in all three cases, since it is dominated by the value of $\sigma_{l+1}$ which corresponds to the behavior of the error predicted by the error bounds. This behavior is clearly seen on the bottom plots of Figure 3, where the  plots indicate that reducing $\Delta$ has miniscule (right plot, Method 2) or zero (left plot, method 1) effect on the magnitude of the error. Note again that the error corresponding to smaller $\Delta$ is not always smaller (middle and bottom right plots) when the value of $\sigma_{l+1}$  is significant. 

Next we  present the results of an experiment with the same nonlinear system where we compare the error from the two methods with the three different snapshot spacings where the  dimension of the RBS is fixed (Figure 4). For the 402-variable FOM, equidistant snapshots in time of both the FOM solution and its time derivative were collected at spacings of $\Delta=0.01, 0.02$ and 0.04. Dimensions of $5, 10, 15, \ldots , 50$ were considered. Presented are 4 of the calculations, with ROM space dimensions =5, 20, 25, 50, to illustrate the observed tendency.  For a very low ROM dimension, {\it e.g.} $l=5$, Method 1 and Method 2 produce similar errors ($\approx 10^{-2}$) and for each method decreasing $\Delta$ does not have any effect on the error (Figure 4). At $l=5$, $\sigma_6\approx 10^{-2}$ for both methods and all three values of $\Delta$ (Figure S4, Supplement) and the magnitude of the error is determined by the basis truncation as suggested by our error bounds.  
We further focus on the case with $\Delta=0.04$ (50 collected snapshots, plots with circles). When the ROM basis has dimension $l = 20$  $||{\bf e}^Y(t)||_2$ becomes larger than   $||{\bf e}^Z(t)||_2$  and remains so  for all fixed RBS dimensions 25, . . ., 50. The same tendency of the error from Method 2 to become smaller than that of Method 1 when increasing the RBS dimension is observed in the case with $\Delta=0.02$ (100 snapshots).  More precisely, for the case when 100 solution snapshots were collected, at $l =25 $, $||{\bf e}^Y(t)||_2$ becomes larger than  $||{\bf e}^Z(t)||_2$  and remains so  for the larger RBS dimensions 30, 35, . . ., 50 (dotted plots). This tendency is likely  true for the third case $\Delta=0.01$ (200 snapshots) as well (calculations with RBS dimensions over 50 were not done).

The latter again demonstrates the potential for achieving better accuracy of approximation when using Method 2 compared to Method 1 for relatively low - dimensional ROMs. 

Similarly to the previous example, plotting the distributions of the singular values helps to get insight about the behavior of the error, observed on Figure 4. The plot of the distributions of the singular values in the two ROMs is presented on Figure S6 (Supplement). It is observed that for a given $\Delta$, there is  a dimension $l$ such that $\sigma_{l+1}$ is comparable to the magnitude of the error produced by both methods. Then Method 2 produces less error than Method 1 for ROM dimensions greater than $l$. 

\begin{figure} [H]
\label{B1}
\epsfig{file=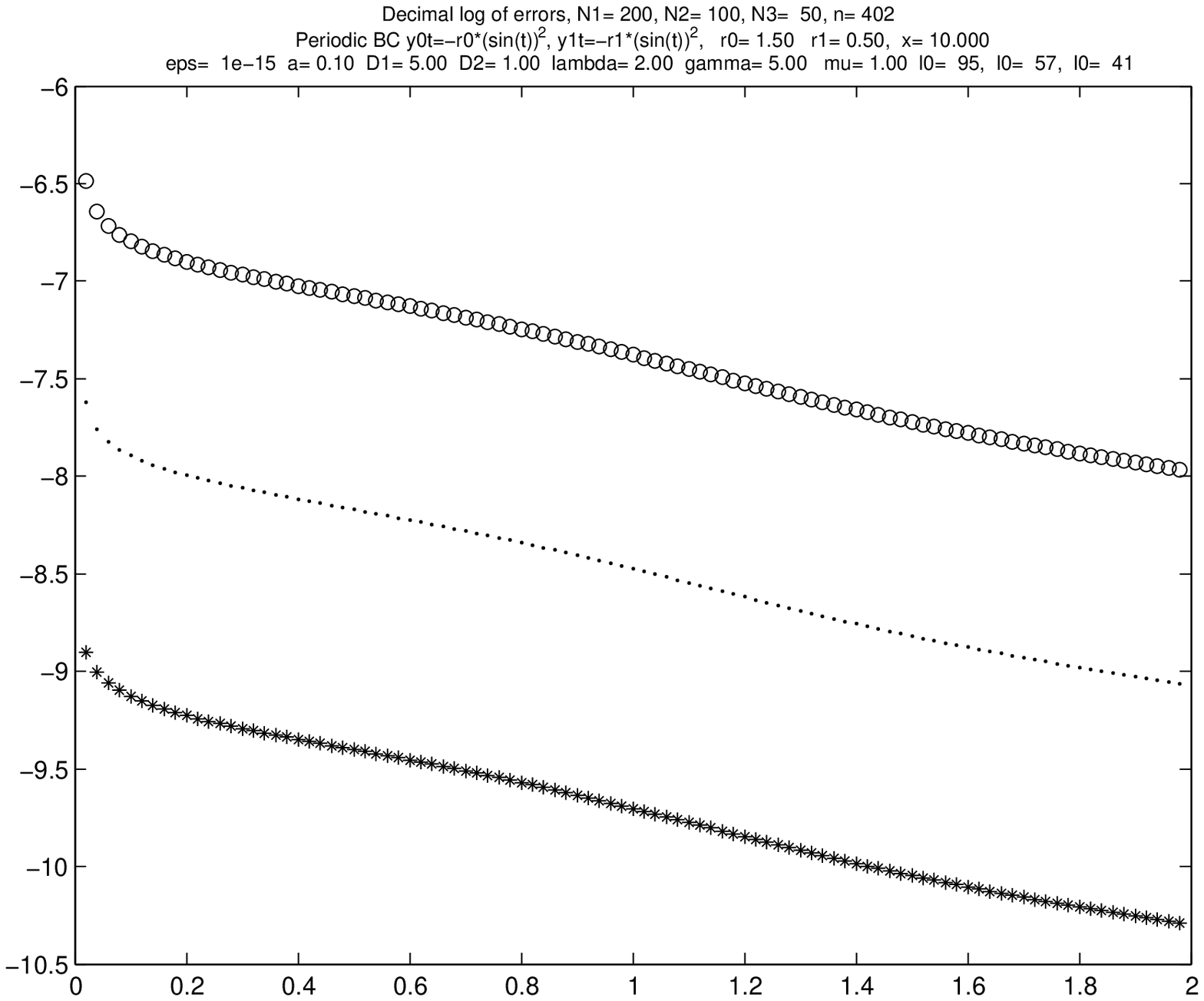, height=2.7in,width=3.2in}
\epsfig{file=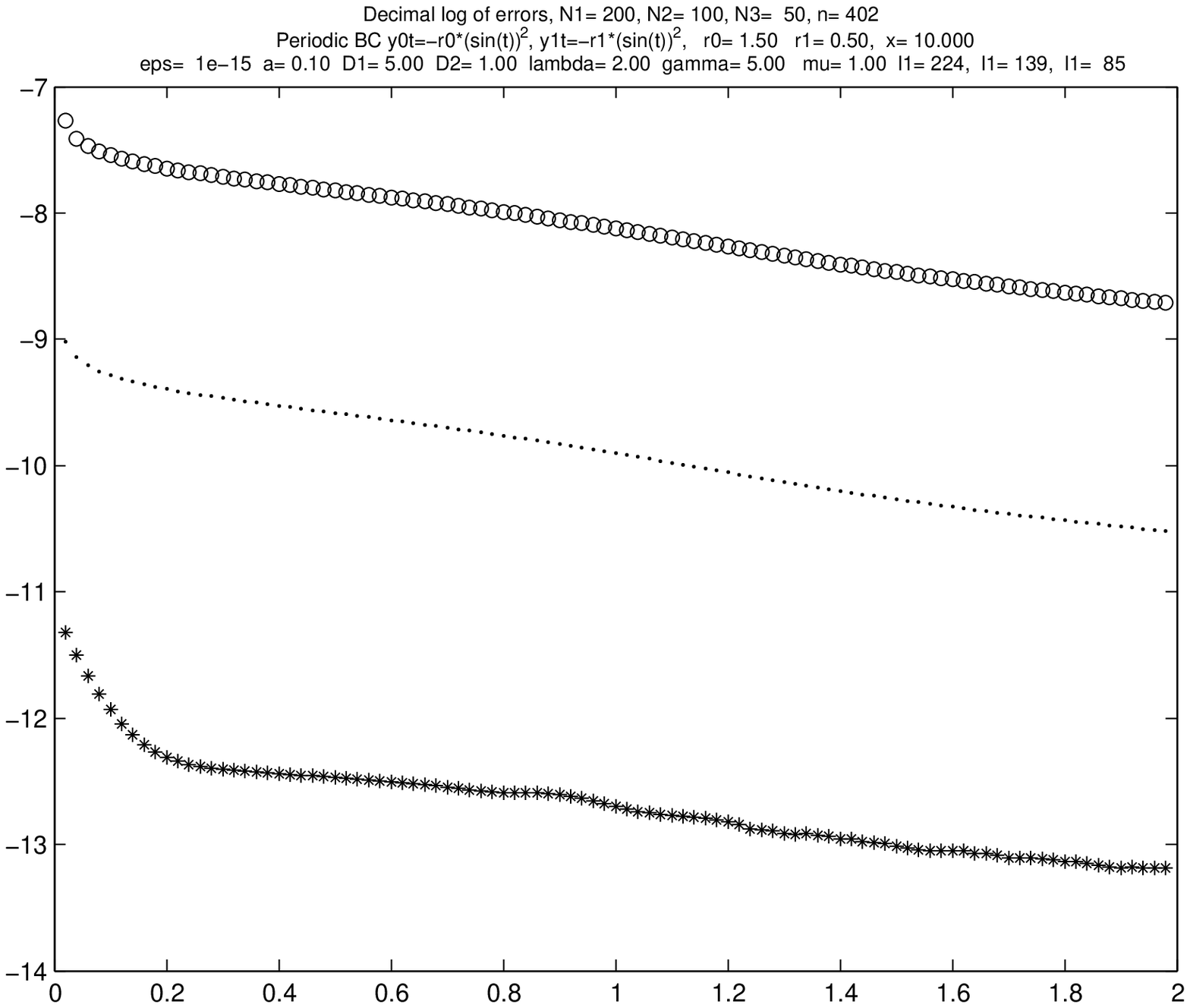, height=2.7in,width=3.2in}\\
\epsfig{file=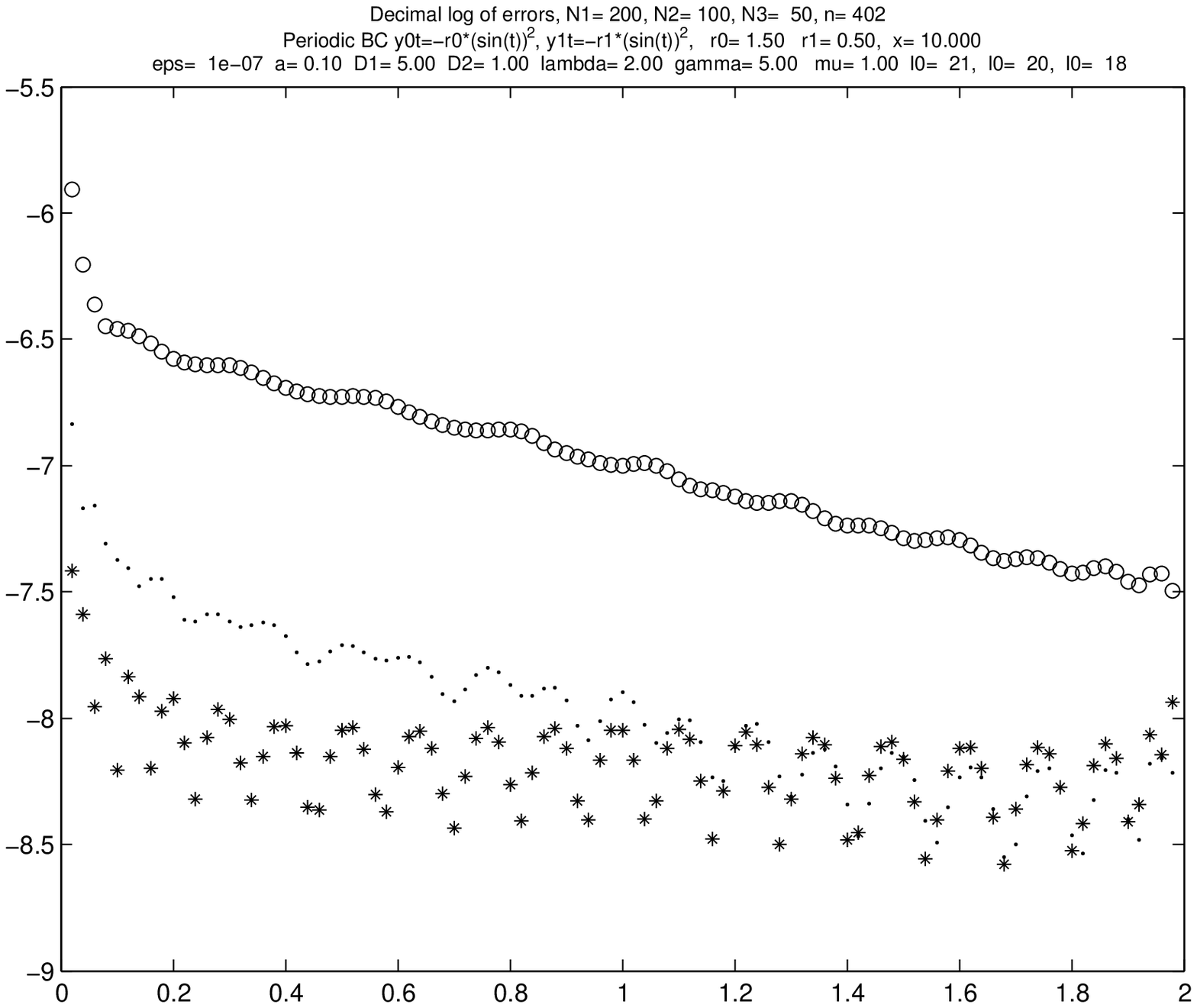, height=2.7in,width=3.2in}
\epsfig{file=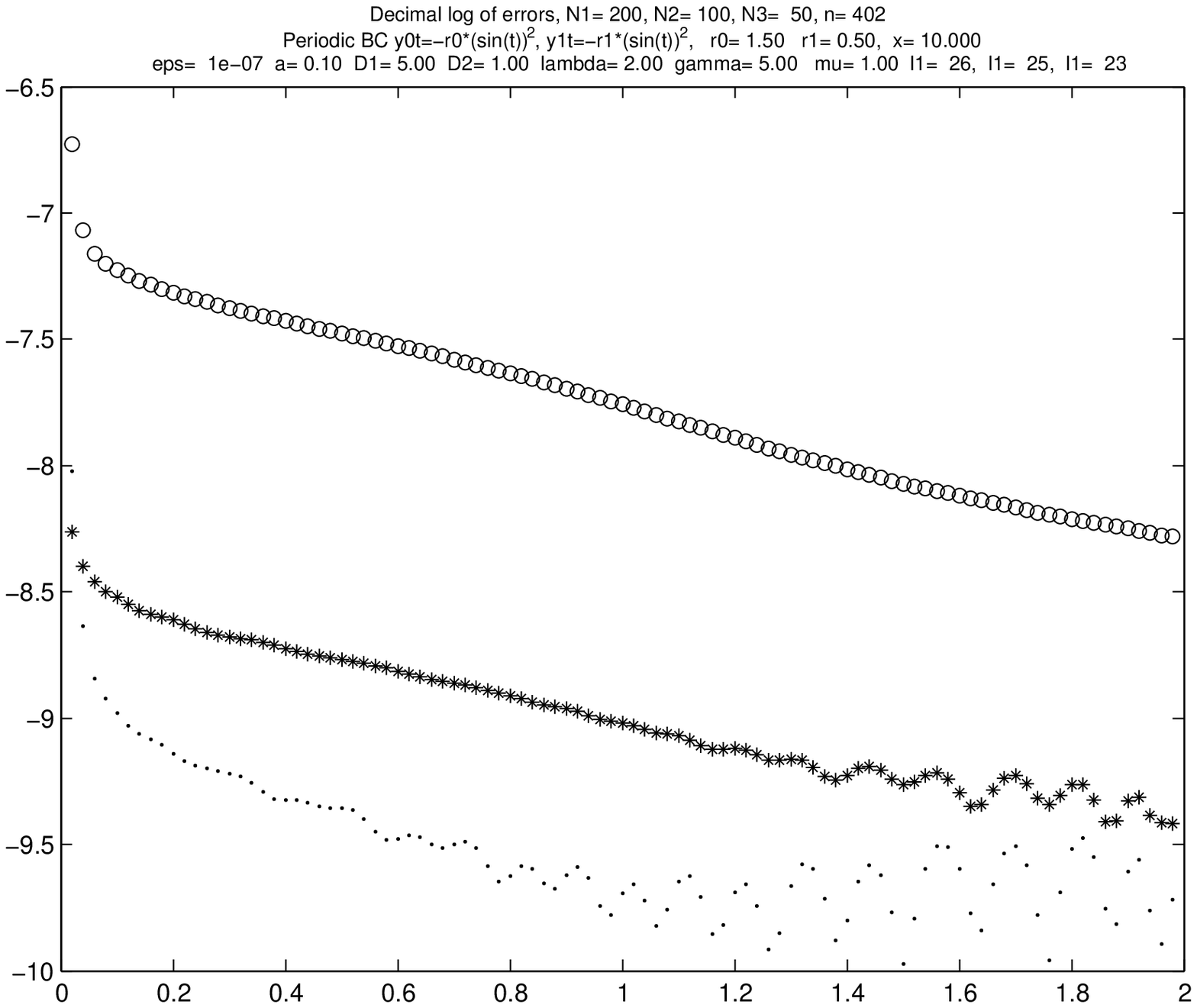, height=2.7in,width=3.2in}\\
\epsfig{file=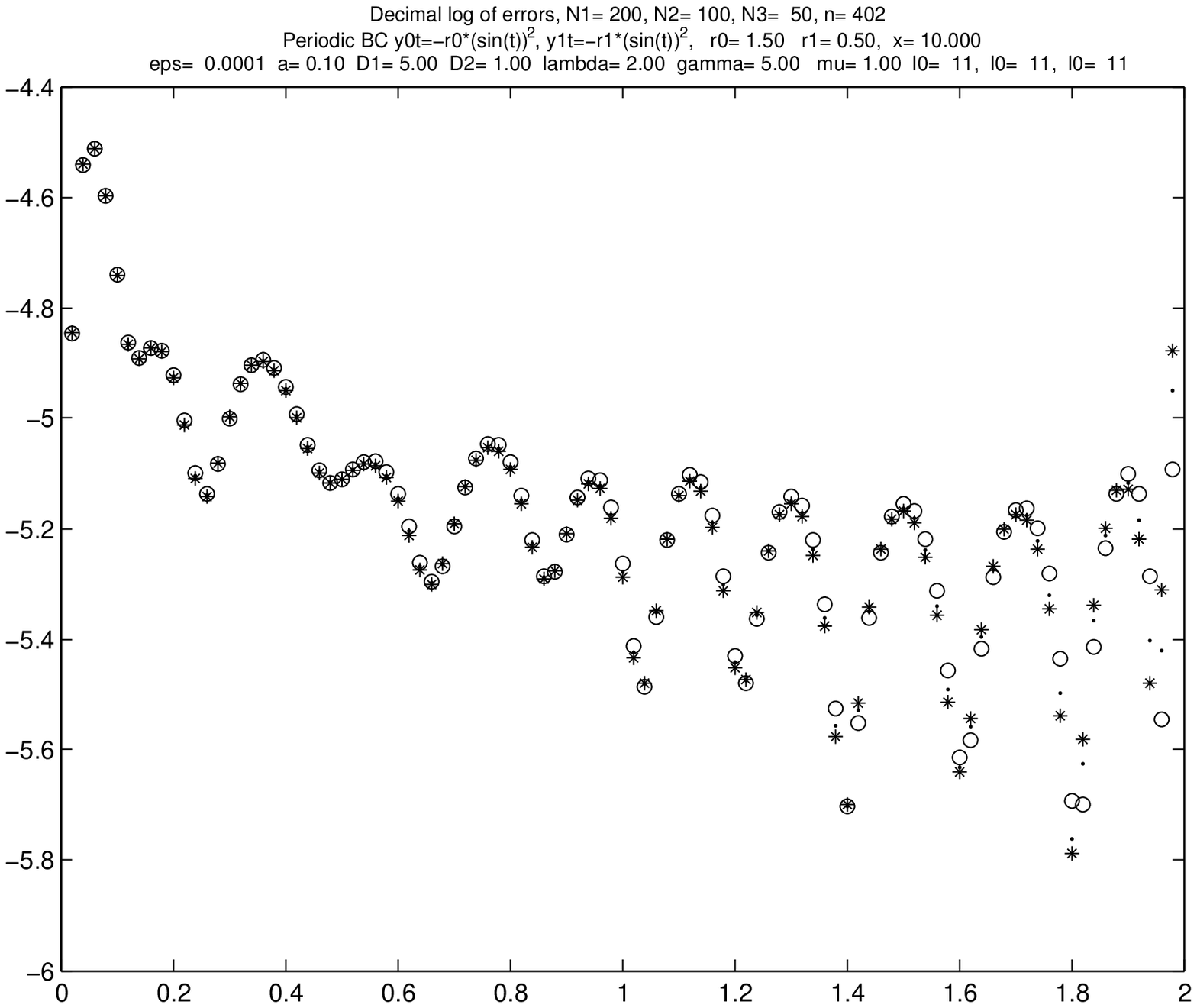, height=2.7in,width=3.2in}
\epsfig{file=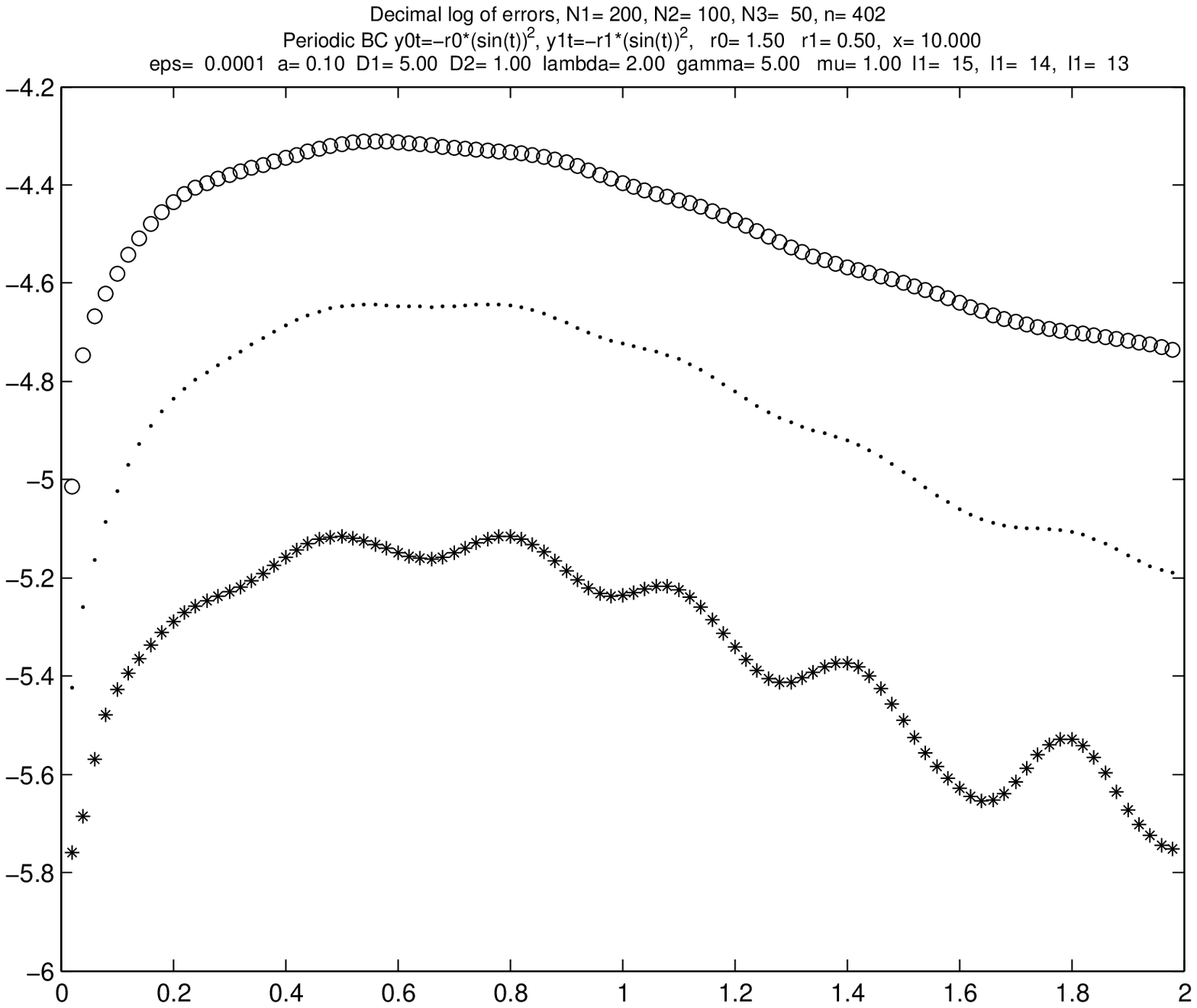, height=2.7in,width=3.2in}\\
\caption{Experiment B.  Error from the two methods at three different values of $\Delta t=0.01, 0.02, 0.04$ and different cutoff ($\varepsilon$) values. The x-axis is time $t$ and the y-axis is $\log_{10} (||{\bf e}^Y(t)||_2)$ (left) and $\log_{10} (||{\bf e}^Z(t)||_2)$ (right).  Circles correspond to $\Delta =0.04$, dots - to $\Delta =0.02$ and crosses to $\Delta =0.01$.  The plots on the right correspond to error from Method 2 (solution and derivative snapshots) and plots on the left correspond to error from Method 1 (no derivative snapshots). }
\end{figure}

\begin{figure}[H]
\label{B2}
\epsfig{file=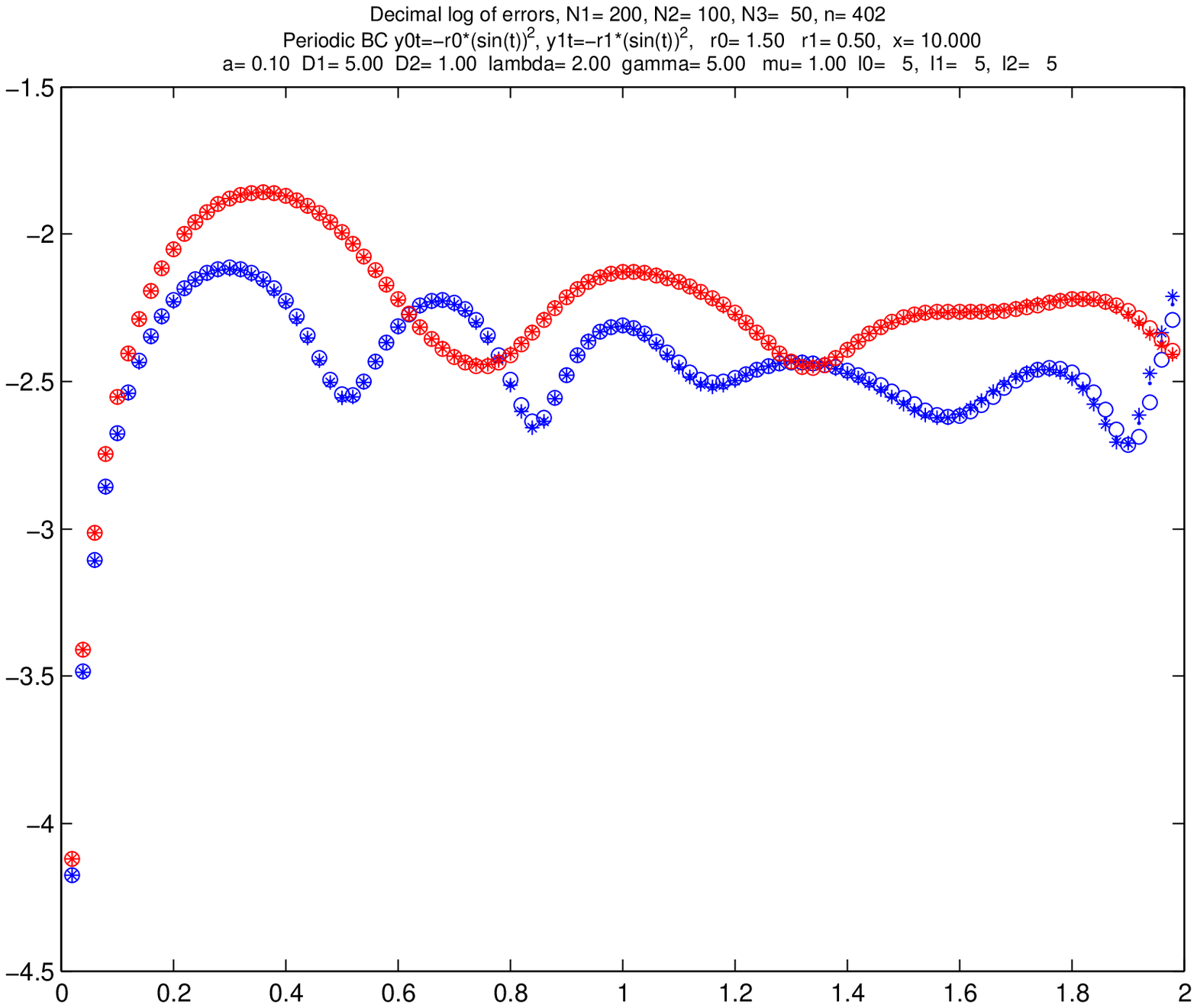, height=2.7in,width=3.2in}
\epsfig{file=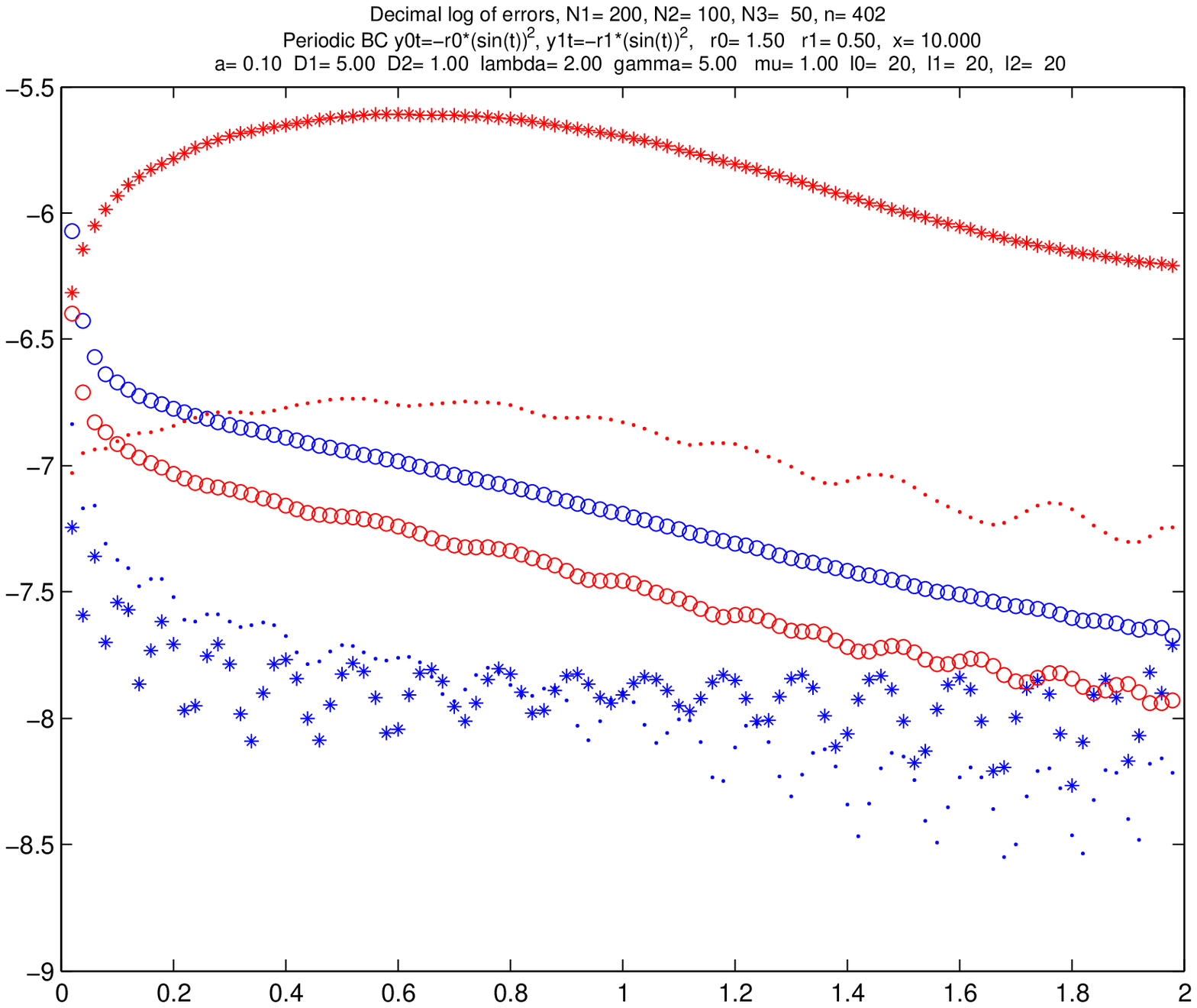, height=2.7in,width=3.2in}\\
\epsfig{file=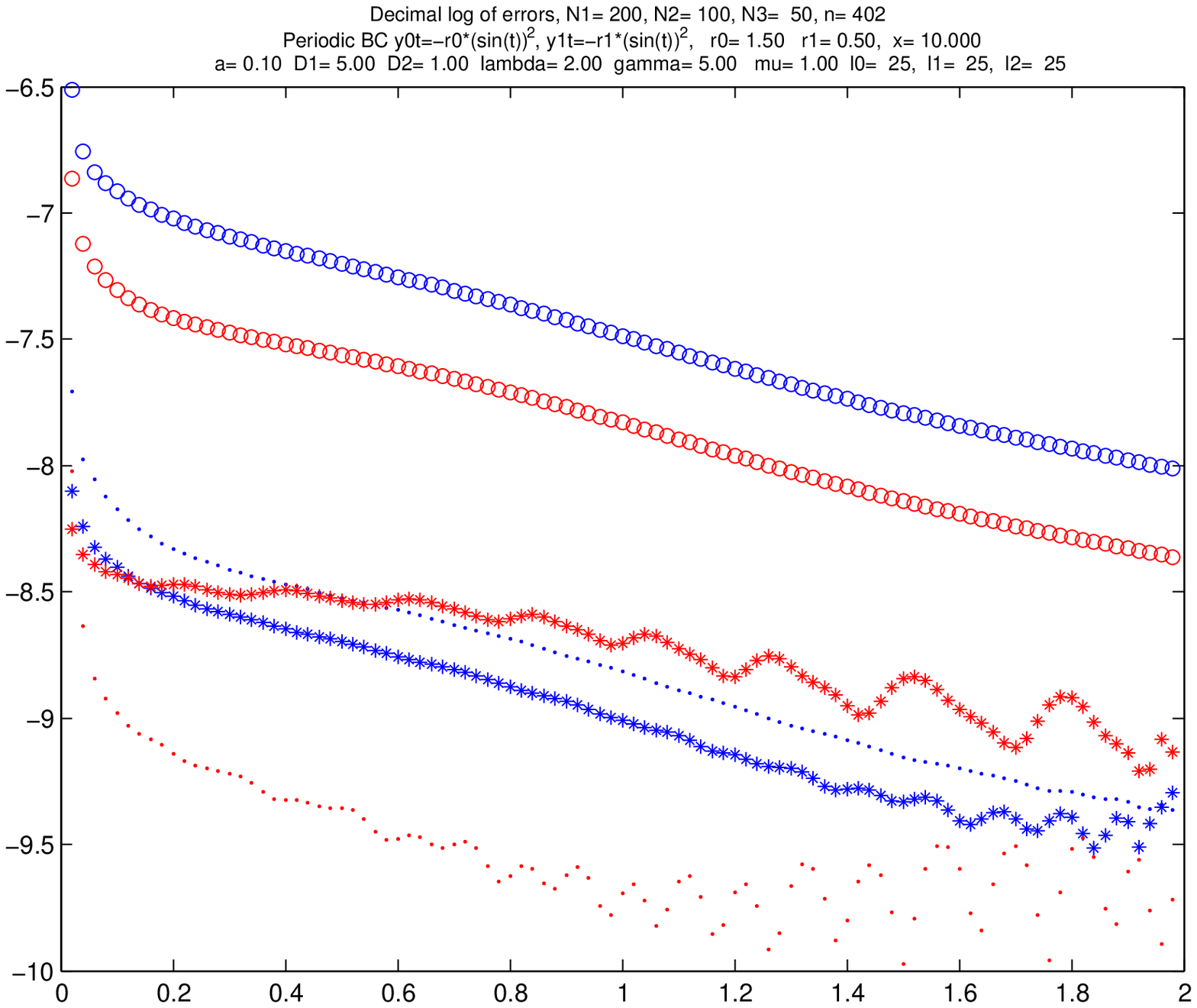, height=2.7in,width=3.2in}
\epsfig{file=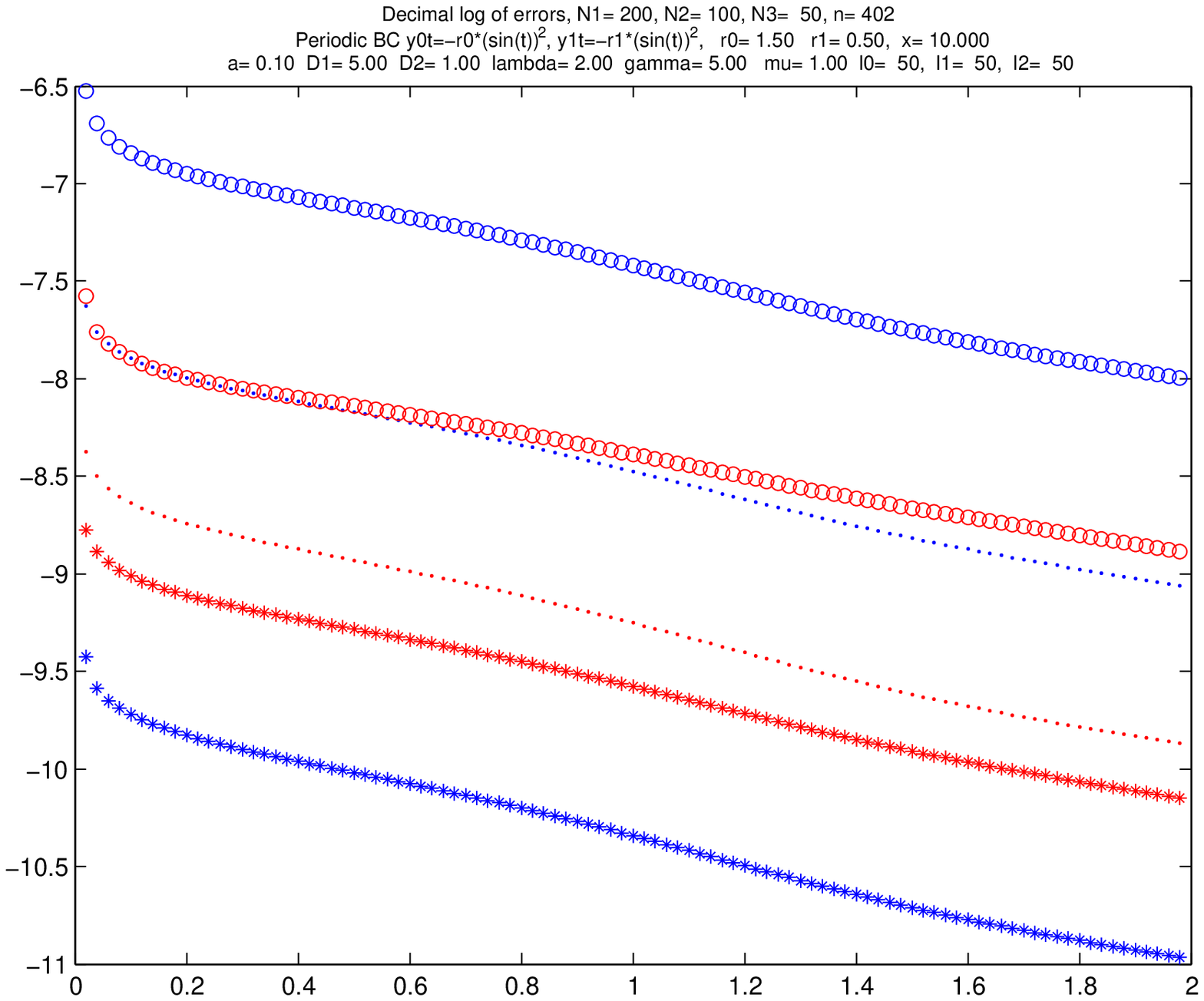, height=2.7in,width=3.2in}\\
\caption{Experiment B.  Error from the two methods at three different values of $\Delta t=0.04, 0.02, 0.01$ and different fixed RBS dimensions ($l=5, 20, 25, 50$ and for the problem in Experiment B. Circles correspond to $\Delta =0.04$, dots - to $\Delta = 0.02$ and crosses to $\Delta =0.01$.  Blue - Method 1, red - Method 2. The y-axis is the decimal logarithm of the error.}
\end{figure}

For $\Delta=0.04$ and $l\geq 25$, $\sigma_{l+1}<10^{-8}$ for  for both Method 1 and Method 2 ROMs and  the maximum of error is of similar order ($10^{-6.5}- 10^{-7.5}$. As we saw, at $l\geq 20$, Method 2 yields more accurate approximation of the FOM than Method 1. As the bounds we derived suggest, this phenomenon is a trade-off between the value of $\sigma_{l+1}$ which at some point (here $l=20$) becomes sufficiently small so that the error is dominated by the terms containing $\Delta^2$ and $\Delta^4$ in both methods. The latter trade-off leads to Method 2 becoming more accurate than Method 1, because of its higher order of approximation $O(\Delta^4)$.

\subsubsection{Experiment C}
In this experiment we explore what happens with the error if the distance between the times at which the snapshots were collected is not small. We consider the same nonlinear system of equations as in experiment B, which was, however, solved numerically on the time interval [0,20]. For the ROMs 40, 20 and 10 equidistant snapshots were collected with $\Delta=0.5, 1, 2$ respectively.   The solution of the system is shown on Figure S5 (Supplement). 

We consider the same type of scenarios as in the previous examples. Having in mind the bounds we derived, we still expect to see decrease in the error when $\sigma_{l+1}$ is very small. This is indeed demonstrated on Figure 5, top plots ($\sigma_{l+1}<10^{-15}$). Again, as in the previous examples, the difference between the plots is larger than expected. Similarly to the previous experiments as well, for larger $\sigma_{l+1}$, the errors become of similar magnitude when $\Delta$ is decreased, indicating that the value of  $\sigma_{l+1}$ dominates the error. Interestingly, when $t$ grows,  the error decreases initially and then stabilizes around some value. For this specific problem, the error produced by Method 1 is larger than the error from Method 2, even for  $\Delta \geq 1$. The error is remarkably small given the relatively small number of snapshots taken. 

Figure 6 presents the results of an experiment with the same nonlinear system where we compare the error from the two methods with the three different time steps and where the  dimension of the RBS is fixed.  Presented are plots of the error for 8 cases where the RBS dimensions are held fixed at 5, 10, 15, 20, 25, 30, 35 and 40. When the dimension of the RBS is 5 or 10, the magnitude of the error is similar for both methods and all distances between the snapshots (top two plots on Figure 6) and also comparable to the magnitude of the first neglected singular value ($\sigma_6, \sigma_{11}$ respectively).   For larger RBS dimensions, the error produced by Method 2 becomes smaller than that produced by Method 1. Specifically, since one cannot construct a POD ROM of higher dimension than the number of snapshots collected, the approximation from Method 2 becomes better than that from Method 1 when increasing the the RBS dimension above 10 (circle symbol plots). 

\begin{figure}[H]
\label{C1}
\epsfig{file=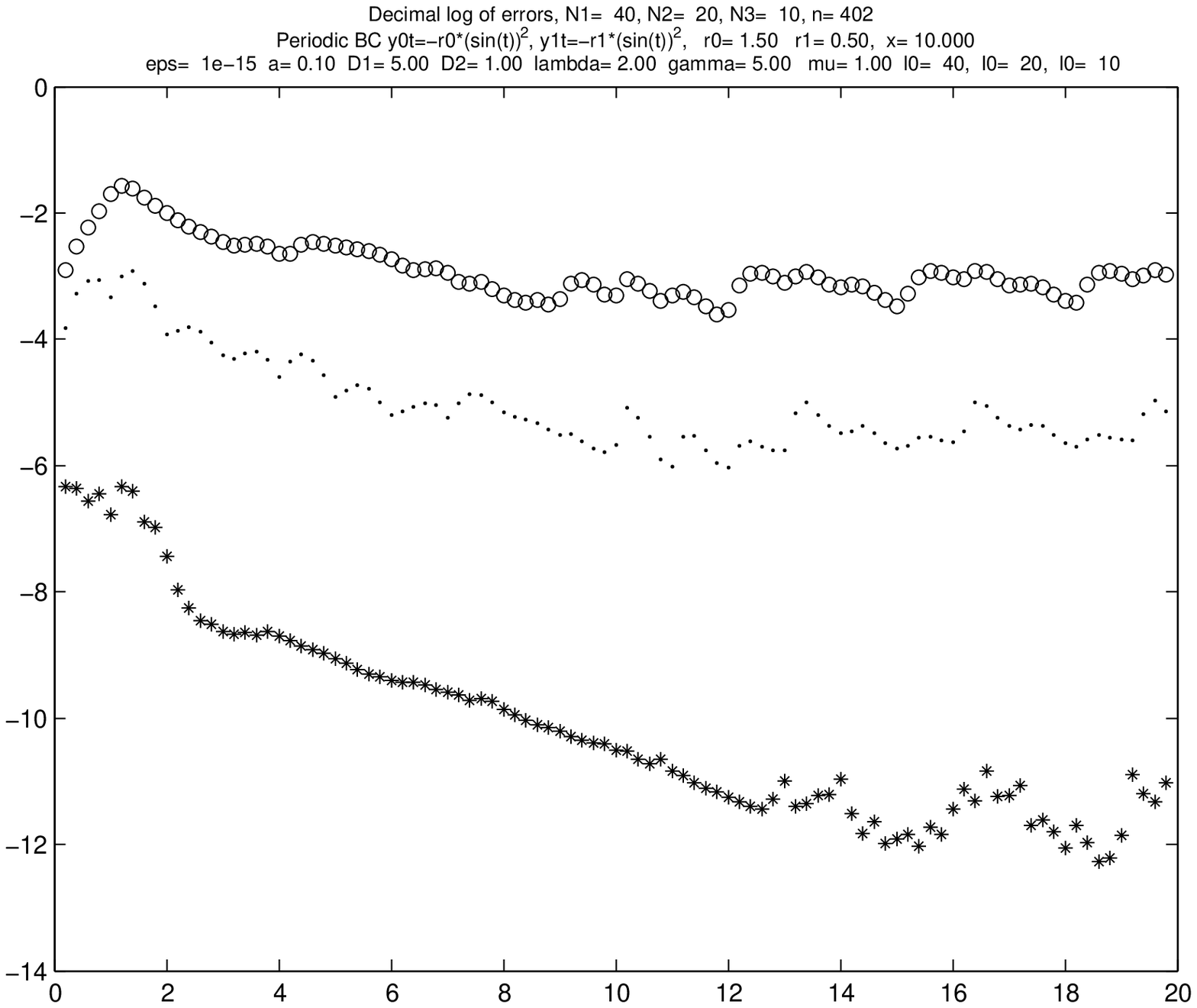, height=2.7in,width=3.2in}
\epsfig{file=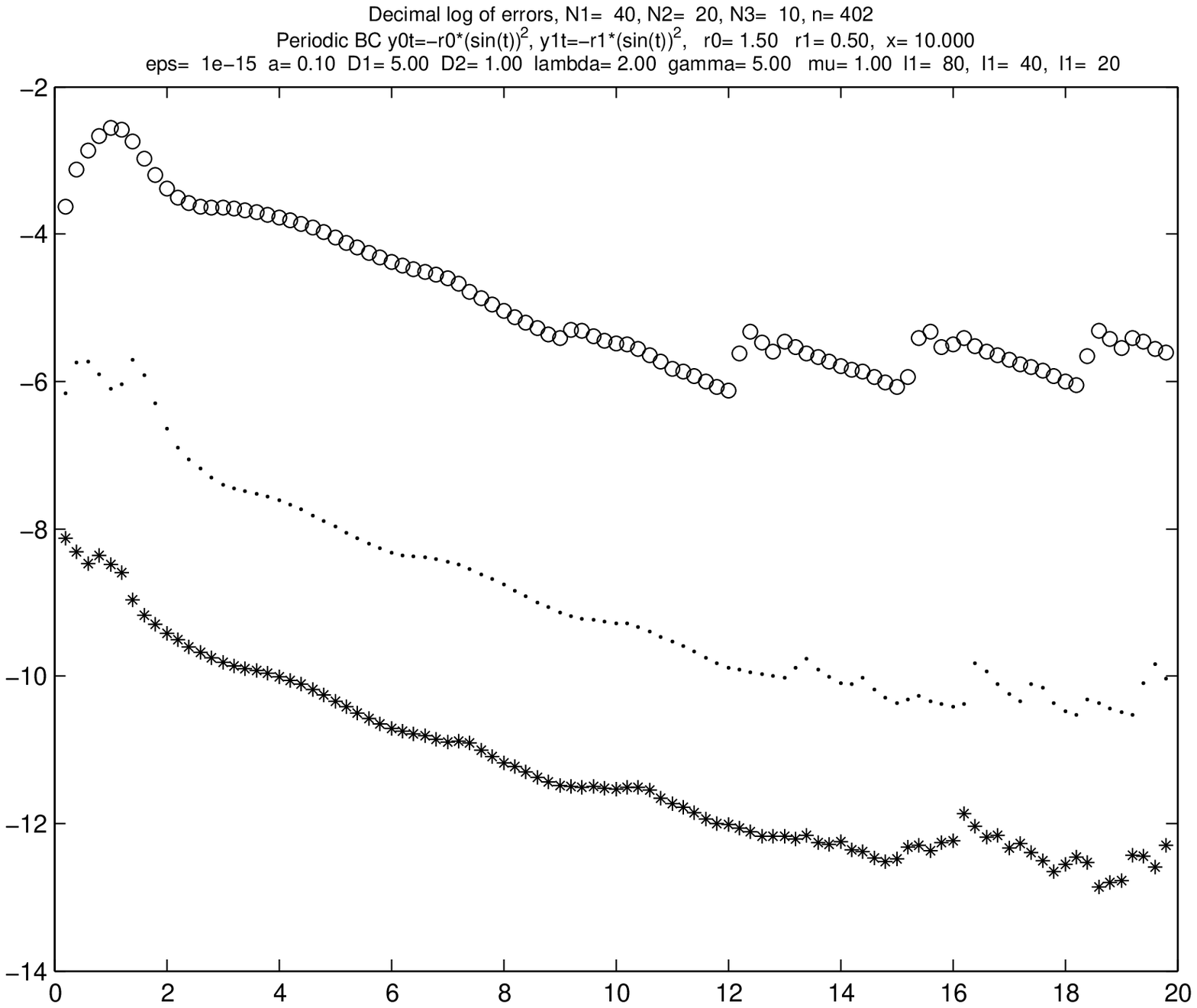, height=2.7in,width=3.2in}\\
\epsfig{file=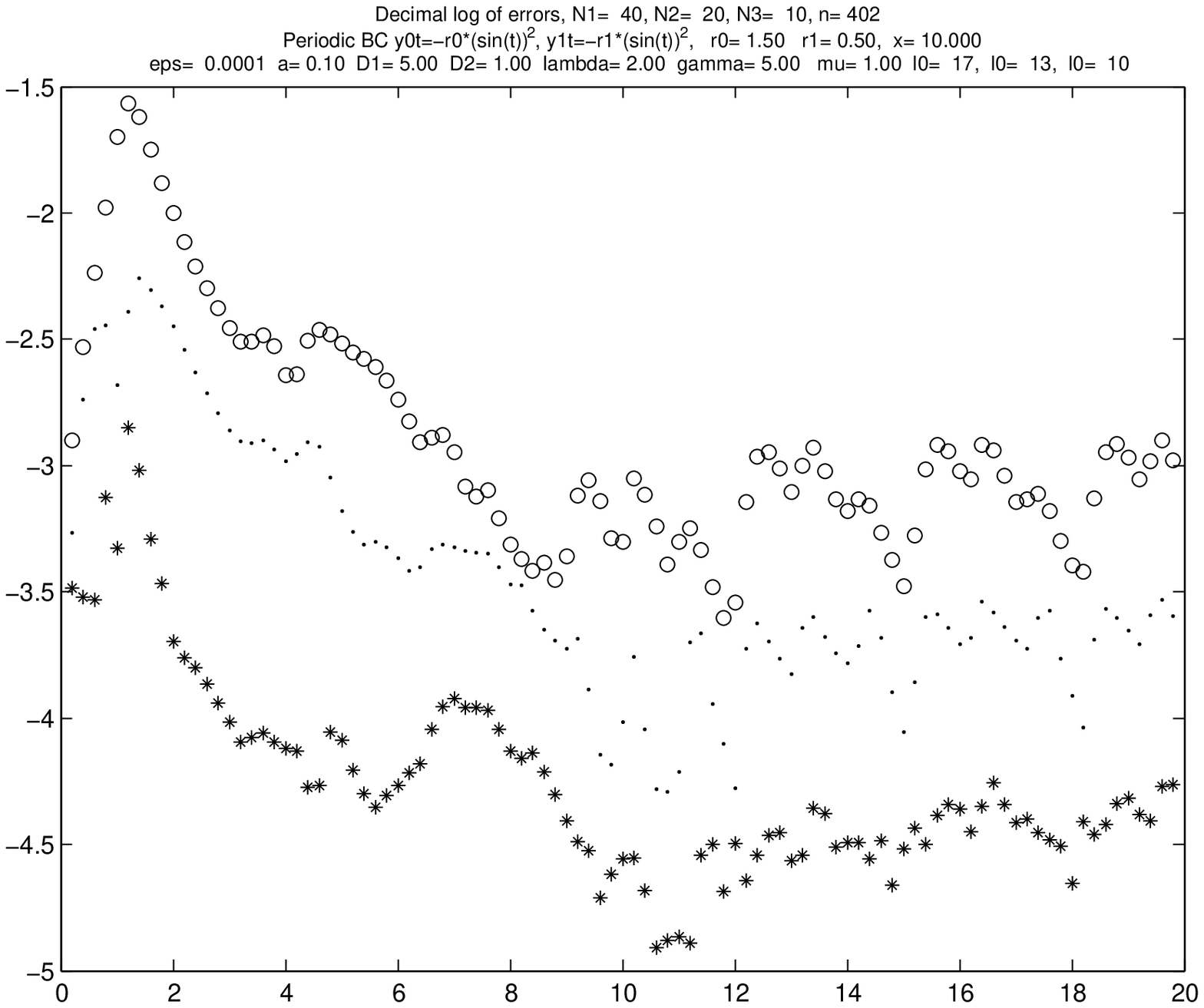, height=2.7in,width=3.2in}
\epsfig{file=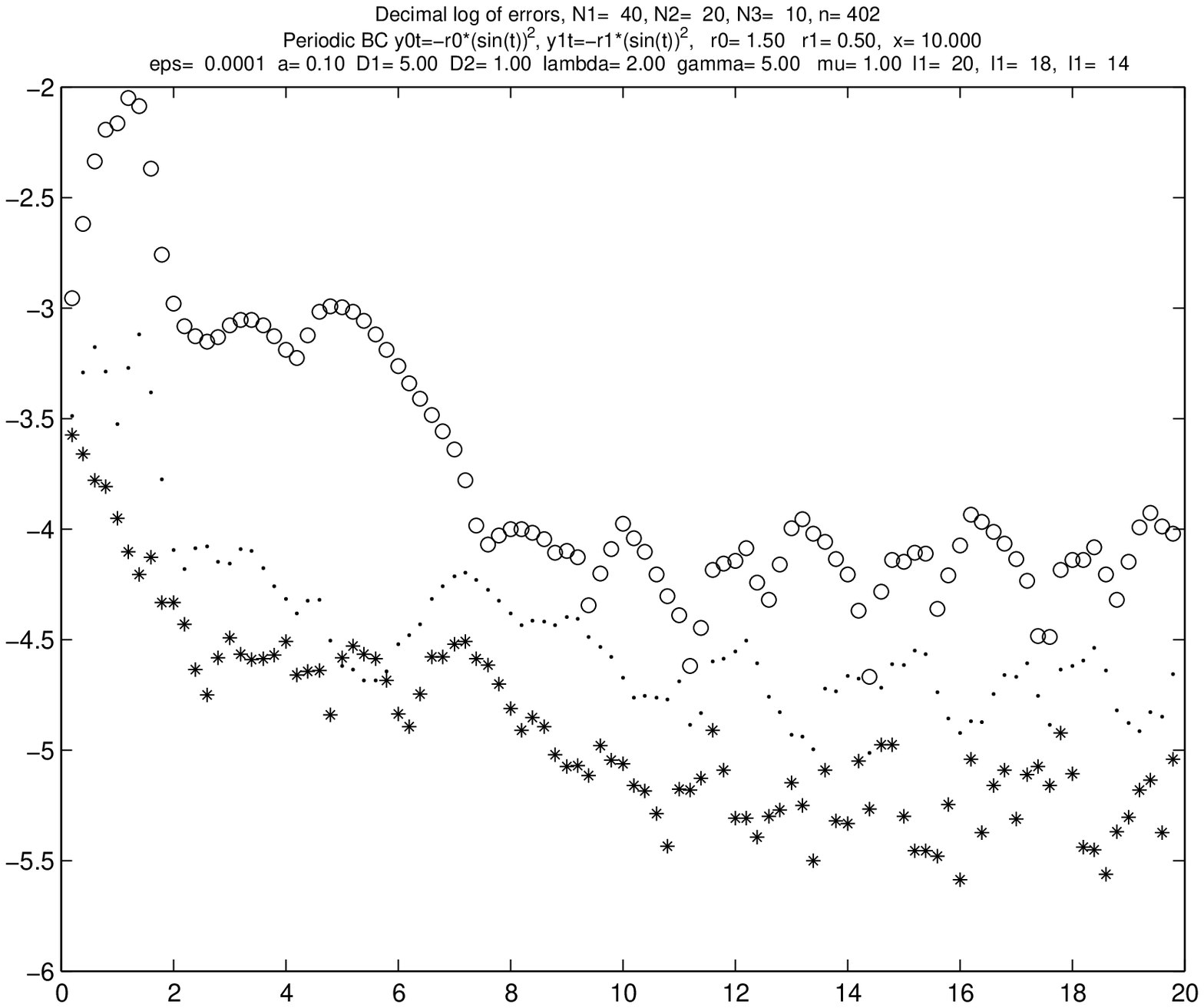, height=2.7in,width=3.2in}\\
\caption{Experiment C.  Error from the two methods at three different values of $\Delta t=0.5, 1, 2$ and different cutoff ($\varepsilon$) values. Circles correspond to $\Delta=0.5$, dots - to $\Delta=1$ and crosses to $\Delta = 2$. Left plot  - Method 1; right  - Method2.}
\end{figure}

\begin{figure}[H]
\label{C2}
\epsfig{file=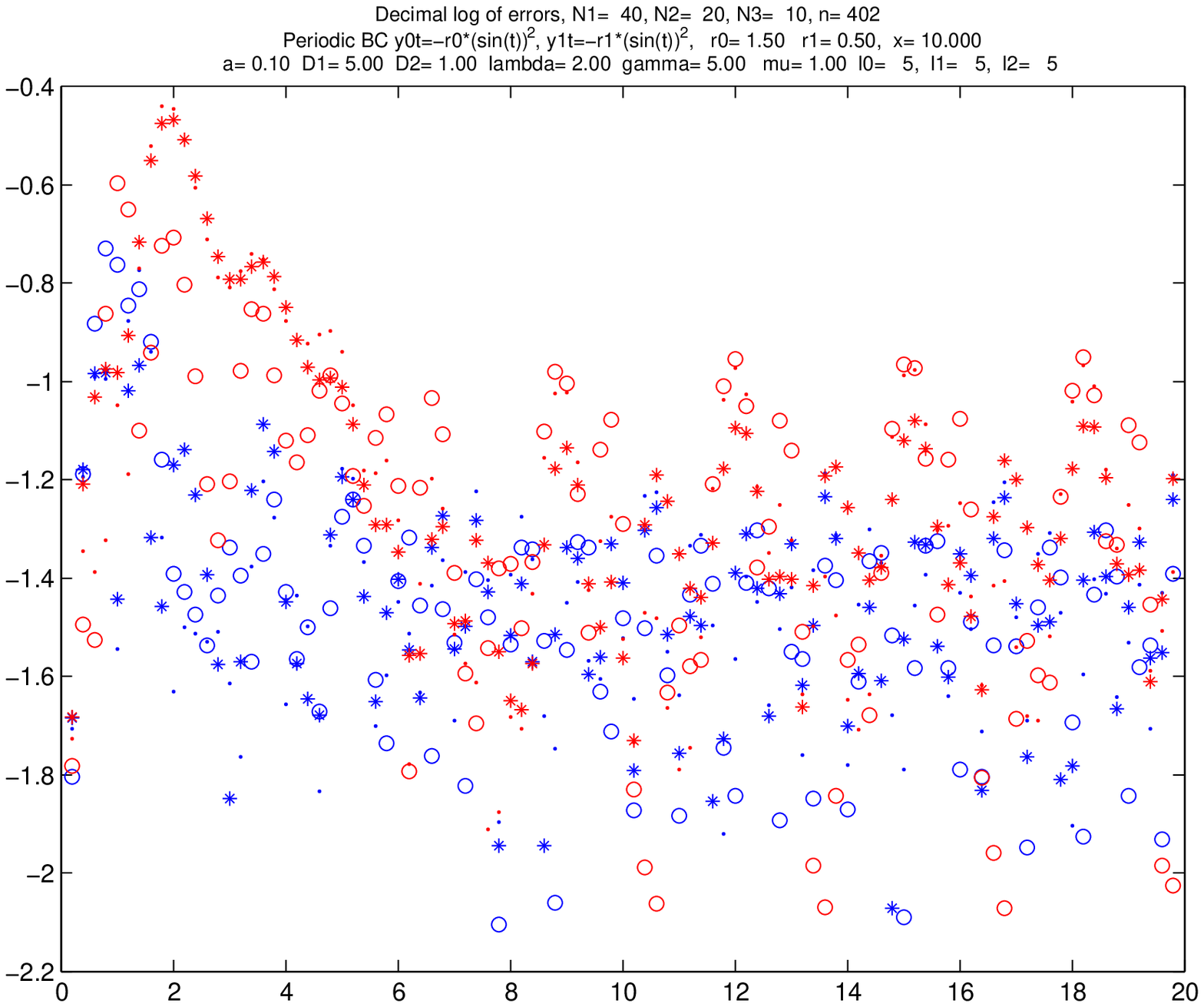, height=2.7in,width=3.2in}
\epsfig{file=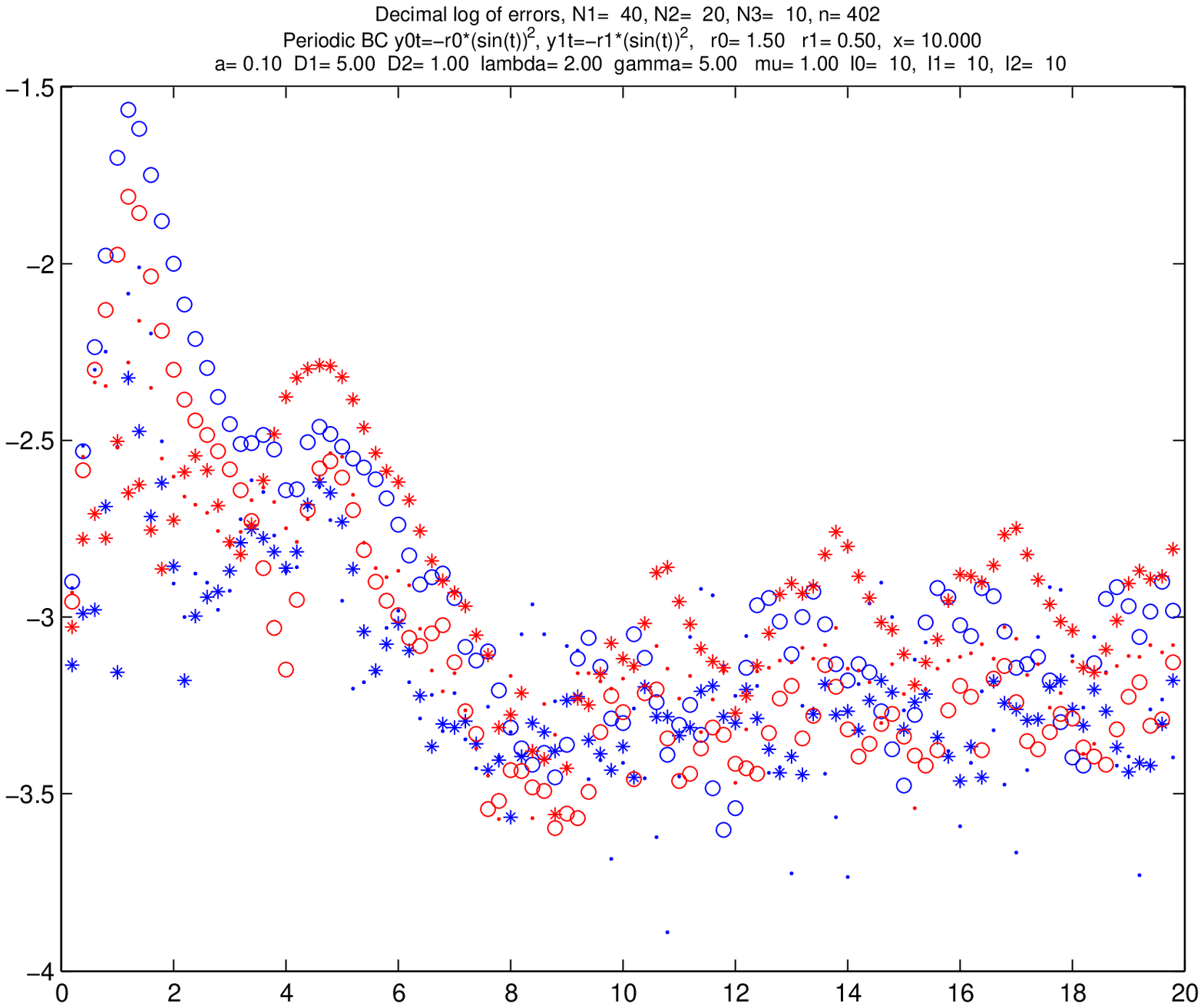, height=2.7in,width=3.2in}\\
\epsfig{file=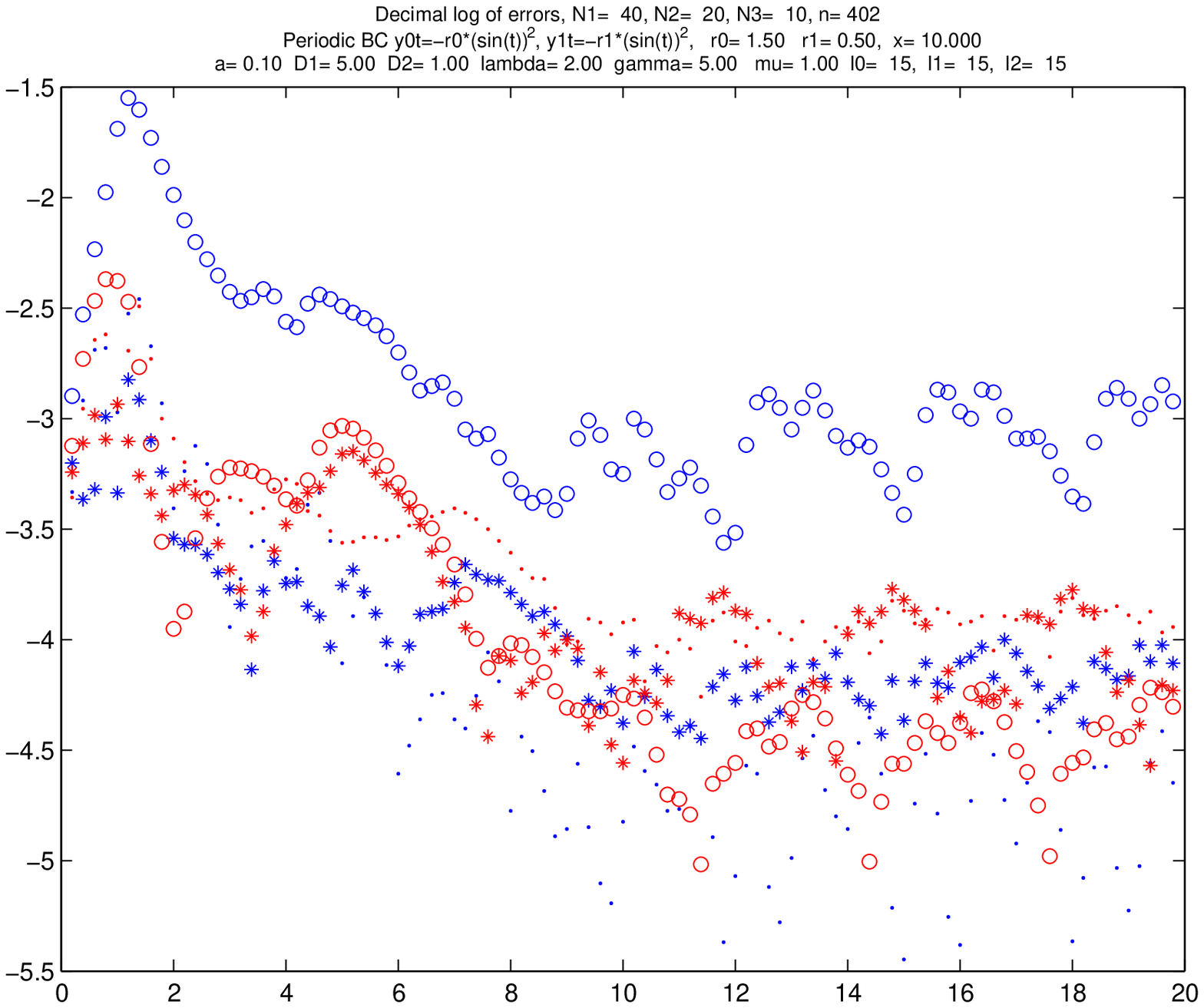, height=2.7in,width=3.2in}
\epsfig{file=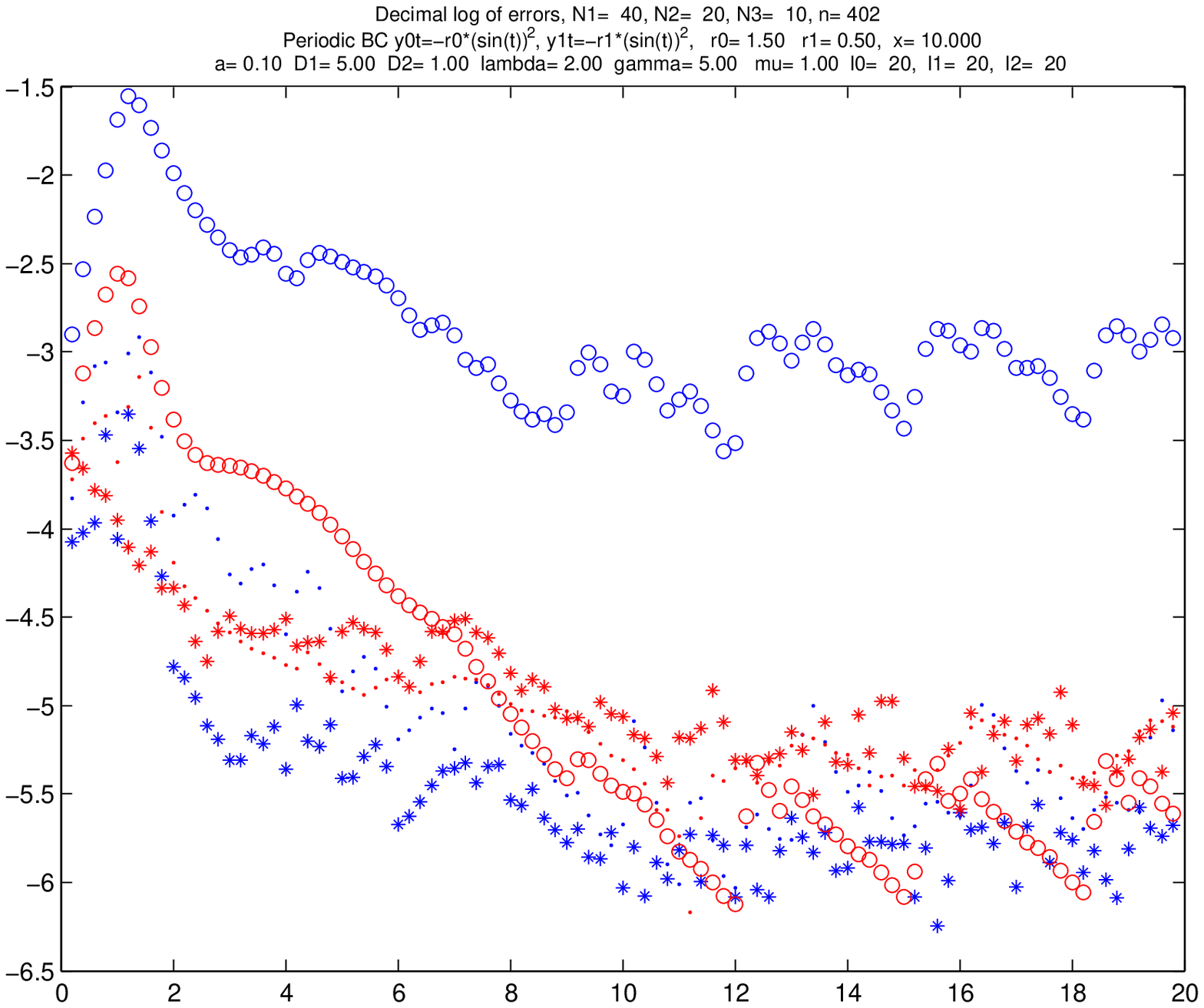, height=2.7in,width=3.2in}\\
\epsfig{file=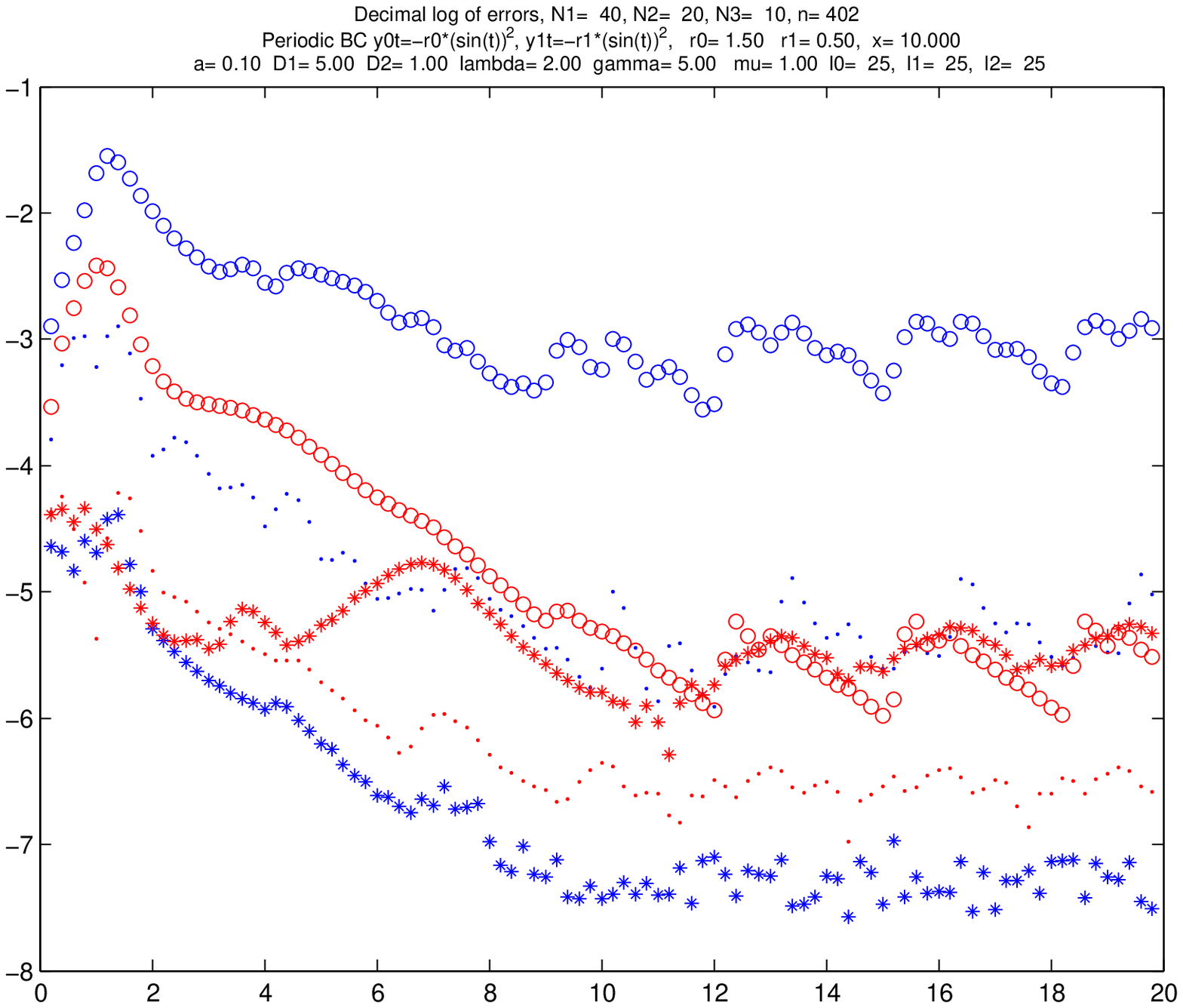, height=2.7in,width=3.2in}
\epsfig{file=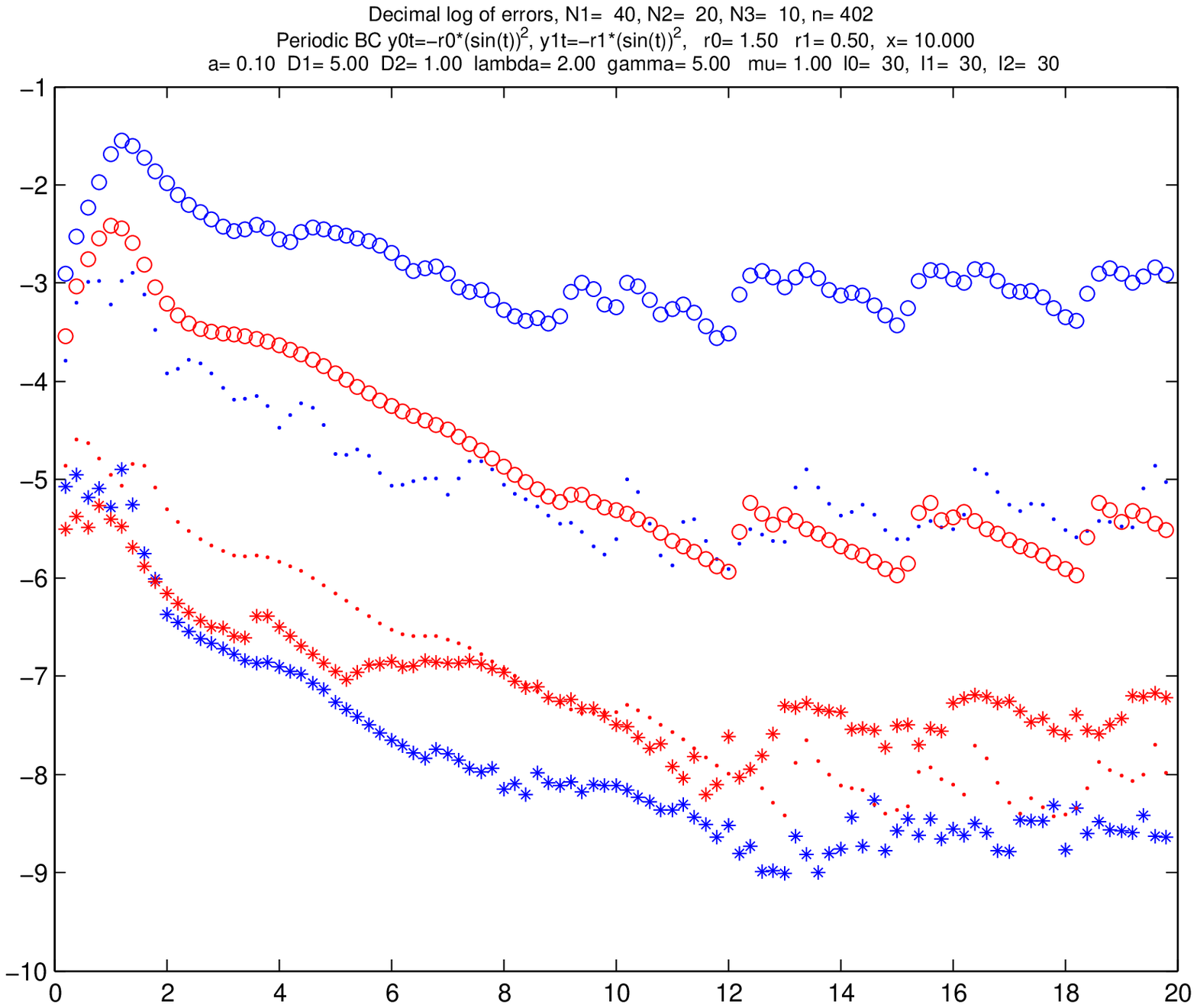, height=2.7in,width=3.2in}\\
\caption{Experiment C. Error from the two methods at three different values of $\Delta =0.5, 1, 2$ and different fixed RBS dimensions ($l=5, 10, 15, 20, 25, 30, 35, 40$). Circles correspond to $\Delta =0.5$, dots - to $\Delta =1$ and crosses to $\Delta =2$. Red - method 2, blue - method 1. The y-axis is the decimal logarithm of the error.}
\end{figure}

Similarly, the error from Method 1 when using 20 snapshots (blue dots, Figure 6) is larger or equal to the error from Method 2 (red dots) when the dimension of the RBS is greater than 20 (third row of plots). 
In Figure S6 (Supplement) we present the distributions of the singular values in the 6 calculations (Method 1 and Method 2, three valued of $\Delta$).

These results demonstrate the possibility of reducing the error of approximation using ROMs constructed via Method 2 compared to ROMs via Method 1 while keeping a low dimension of the ROM.

\section{Discussion}
Time derivatives include important information about the behavior of solutions of time dependent problems. It is reasonable to expect that POD ROMs including  this additional information would approximate more accurately the solution of the FOM. It can be argued that difference quotients, used by some authors as discussed in the Introduction, are good approximations of the time-derivatives, which is true when the time intervals between snapshots are sufficiently small. DQs are not approximations if the solution snapshots are taken at relatively large intervals. Further, for dynamical systems of the form (\ref{system1})  calculating the derivatives at the times at which solution snapshots were calculated, come at almost no additional  computational cost. Finally, unlike DQs, derivative snapshots do not belong to the reduced space generated by the solution snapshots. These arguments justify a comparative study of the accuracy of the POD ROM based on two types of snapshot selection - with and without using time derivatives as snapshots.

We have derived error bounds that suggest that using time derivative snapshots may decrease the approximation error if the first neglected singular value is not too large. We also demonstrate by numerical examples that the method with time derivative snapshots can yield significantly smaller (error of approximation. Specifically, we show that when we take a small number of snapshots and use all of them to define the ROMs, Method 2 produces considerably smaller approximation error than Method 1.

The general method we use to derive the error bounds in this paper has been used by other authors: deriving an equation for the error, directly integrating it and applying the Gronwall inequality by using assumptions for sufficient smoothness of the right-hand side of the system and its solutions and the respective Lipschitz constants. The innovative part of the method used here is the application of interpolation methods (Lagrange and Hermite interpolations in the two cases considered). This approach enabled the derivation of error bounds containing the time intervals $\Delta_i$ between the snapshots and the orders of approximation expressed in terms of  $\Delta_i$. 

This analysis of the error is the first emphasizing the relative significance of the two sources of the error -  one coming from the reduced dimension  $l$, via the size of the first neglected singular value of the snapshot matrix, and the other from the term $O(\Delta^\alpha)$, where $\alpha$ has different values (4 and 2)  for the two methods (with and without time derivative snapshots). We believe that the study presented here contributes to understanding the error in POD ROMs. It is notable that even though we derive upper bounds for the error, and not estimates, these bounds are quite informative about these sources of error and the numerical experiments support what is expected from the bounds. 

Based on the numerical experiments performed it seems that the order $\alpha$ of the error may be higher than predicted by the error bounds we derive. More accurate bounds result from using the logarithmic norm as in \cite{rathinam} (which would primarily affect the exponential term in the error bounds) or other methods of error estimation that yet need to be defined.  Our further goal is to devise methods using these or improved bounds for rational selection of the time moments of the snapshots. 

\section{Acknowledgements}

This work was performed under the auspices of the U.S. Department
of Energy by Lawrence Livermore National Laboratory under contract
DE-AC52-07NA27344. 

Early discussions with William Henshaw are gratefully acknowledged.
All authors were supported by LLNL-LDRD 13-ERD-031 grant.

\vfill\eject

\voffset= -0.5in
\section{Supplement}

\beginsupplement
\begin{figure}[H]
\epsfig{file=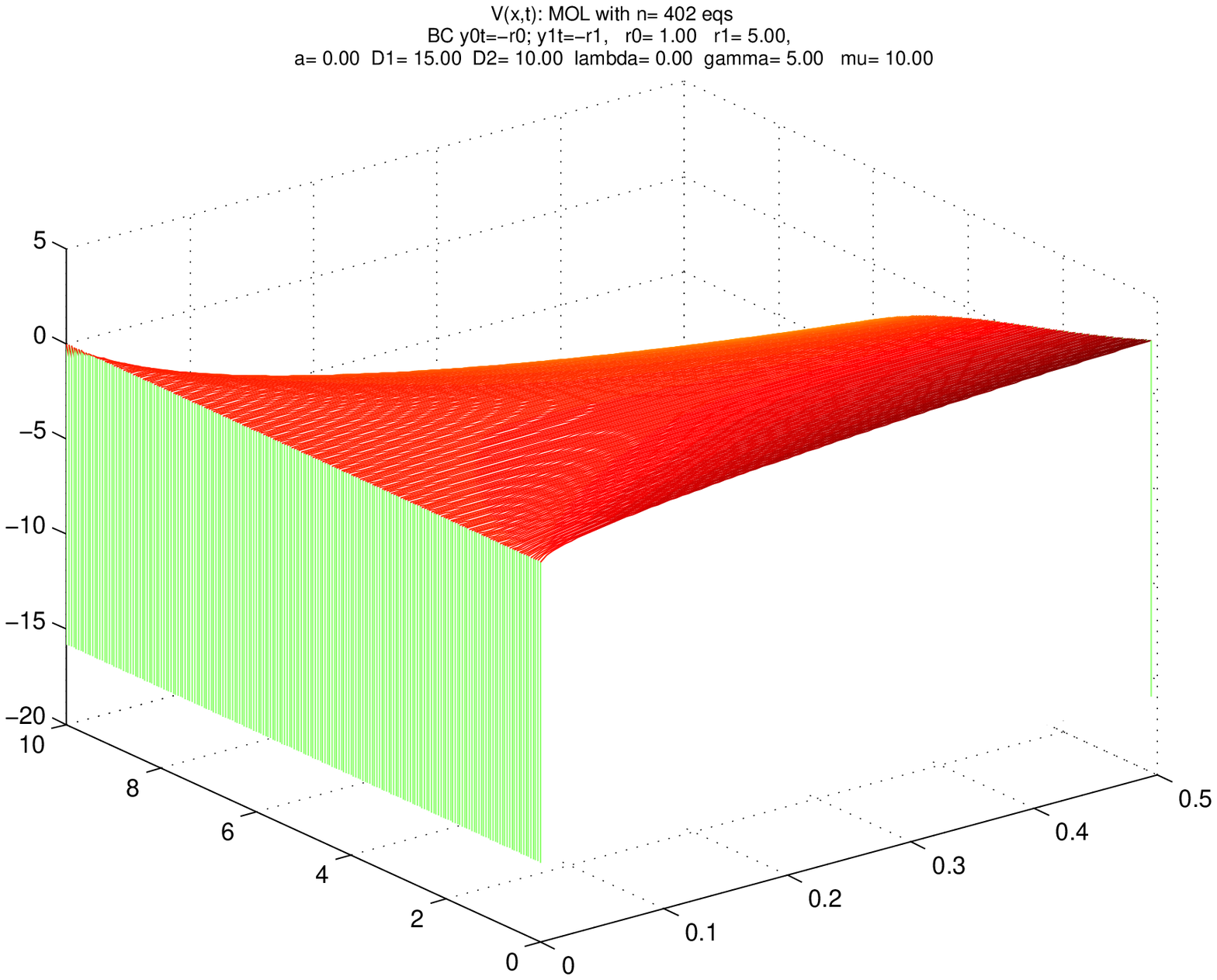, height=3.2in,width=4.5in}
\caption{Experiment A. Solution of the full system. x-axis - time, y-axis - space. The plot represents a set of 201 plots of the calculated solution $v_i(t)$, i=0, 200.} 
\end{figure}

\begin{figure}[H]
\label{A3}
\epsfig{file=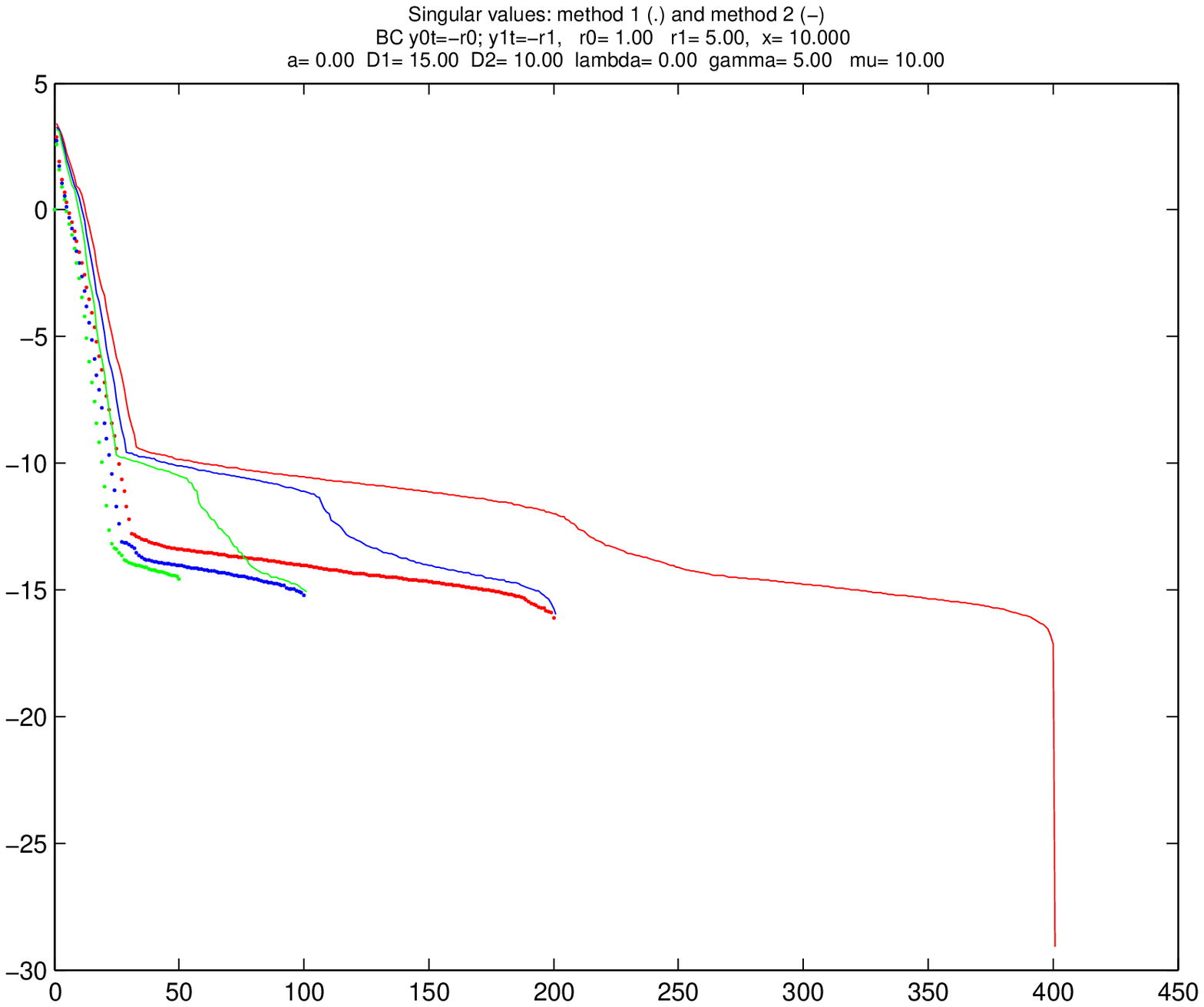, height=3.2in,width=4.5in}
\caption{Experiment A.  Distribution of the singular values  for the two ROMs and three different snapshot spacings. Dots: Method 1; continuous line: Method 2; Red: $\Delta=0.0025$; Blue: $\Delta=0.005$; Green: $\Delta=0.01$. }
\end{figure}

\begin{figure}[H]
\epsfig{file=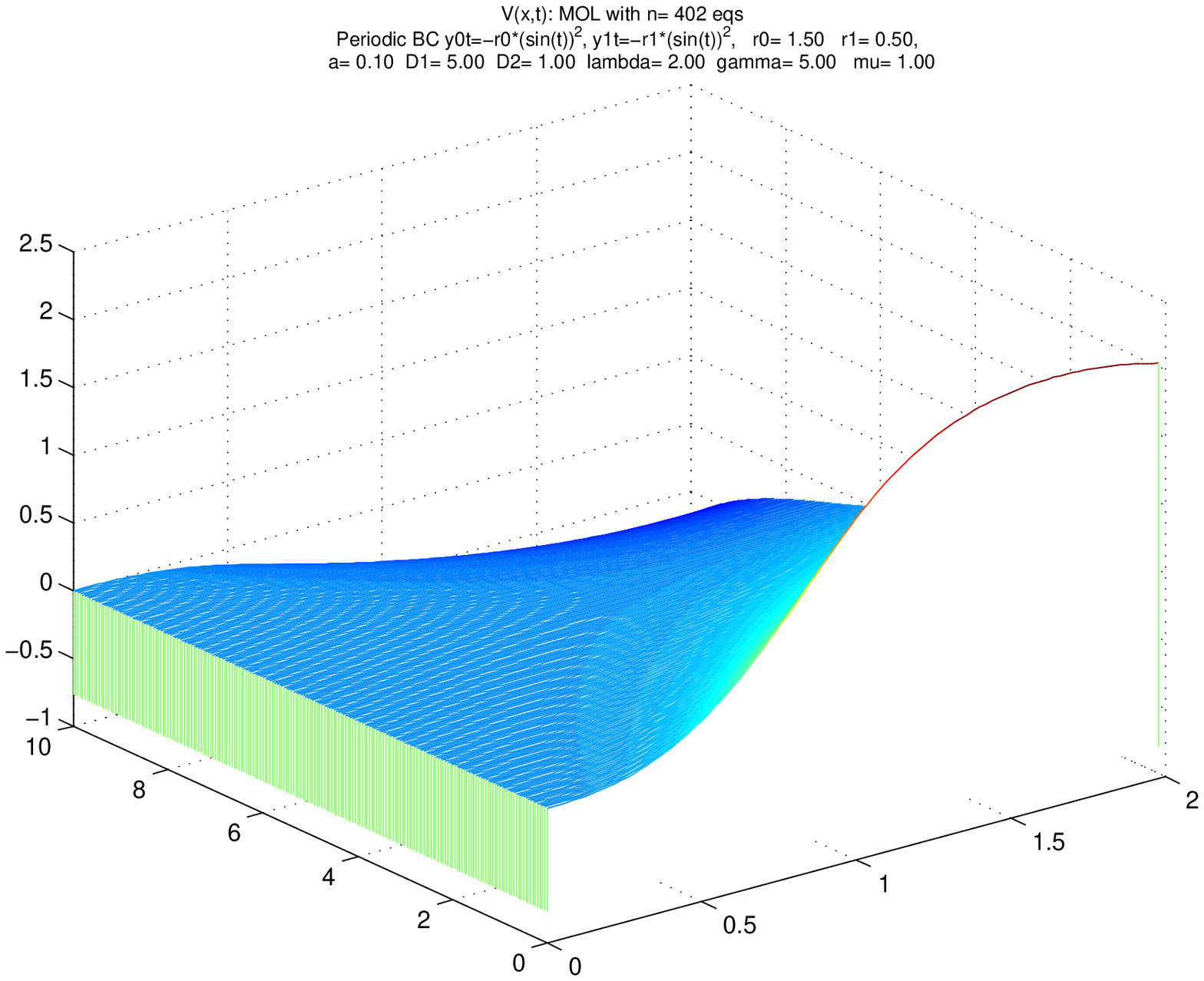, height=3.2in,width=4.5in}
\caption{Experiment B. Solution of the full system. x-axis - time, y-axis - space. The plot represents a set of 201 plots of the calculated solution $v_j(t)$, i=0, 200.} 
\end{figure}

\begin{figure}[H]
\label{B3}
\epsfig{file=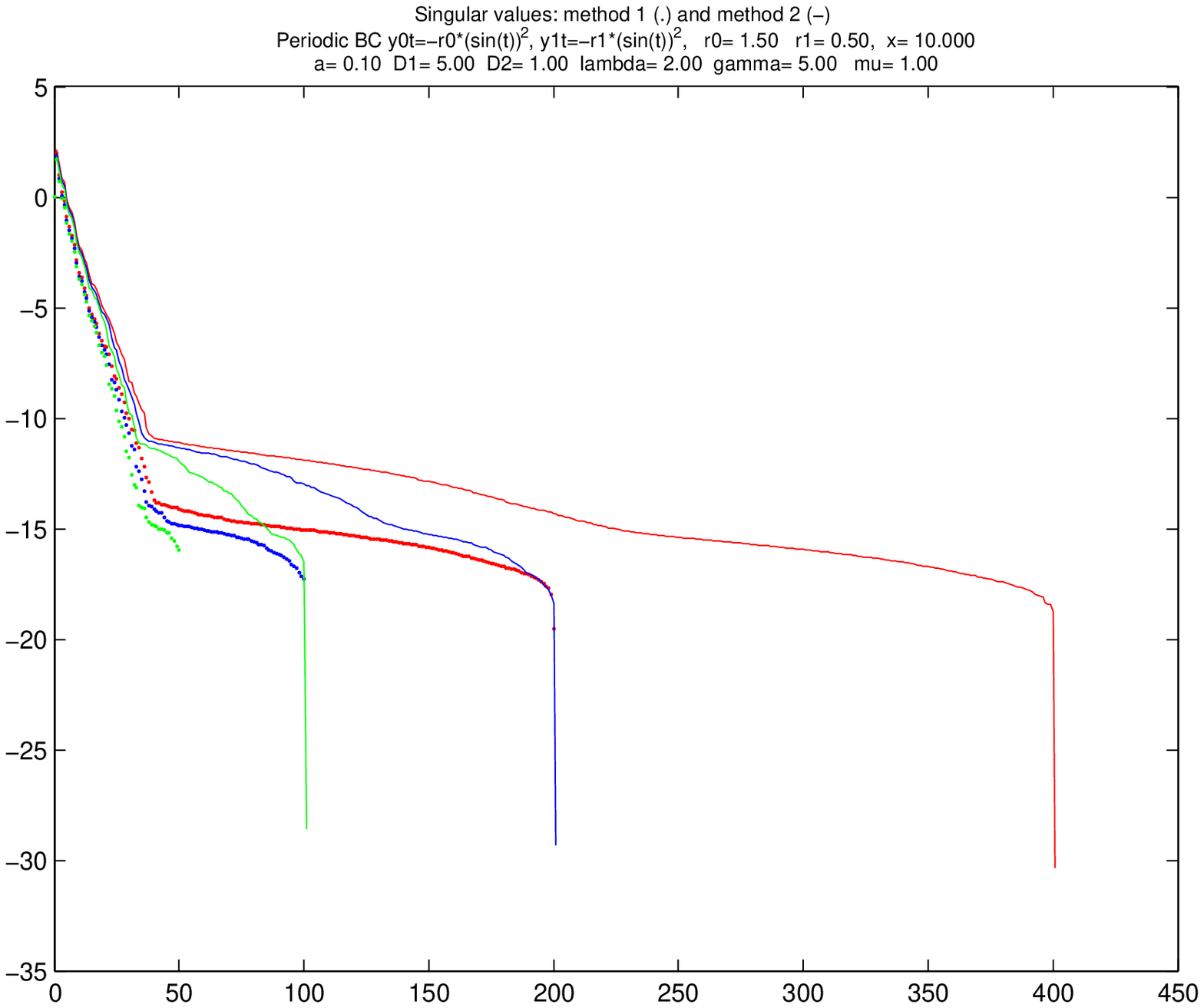, height=3.2in,width=4.5in}
\caption{Experiment B.  Distribution of the singular values  for the two ROMs and three different snapshot spacings. Dots: Method 1; continuous line: Method 2; Red: $\Delta=0.01$; Blue: $\Delta=0.02$; Green: $\Delta=0.04$. }
\end{figure}

\begin{figure}[H]
\epsfig{file=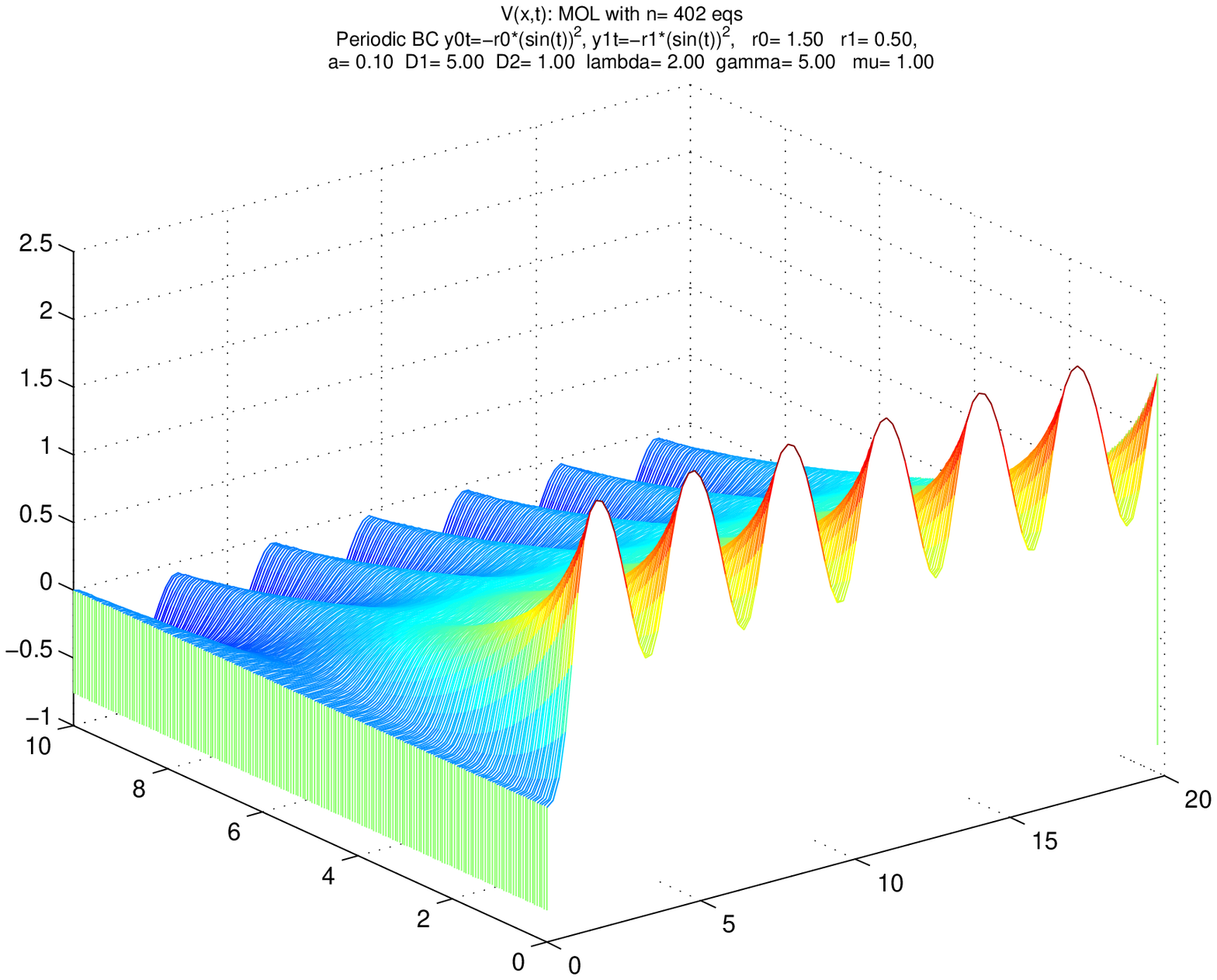, height=3.2in,width=4.5in}
\caption{Experiment C. Solution of the full system. x-axis - time, y-axis - space. The plot represents a set of 201 plots of the calculated solution $v_i(t)$, i=0, 200.} 
\end{figure}

\begin{figure}[H]
\label{C3}
\epsfig{file=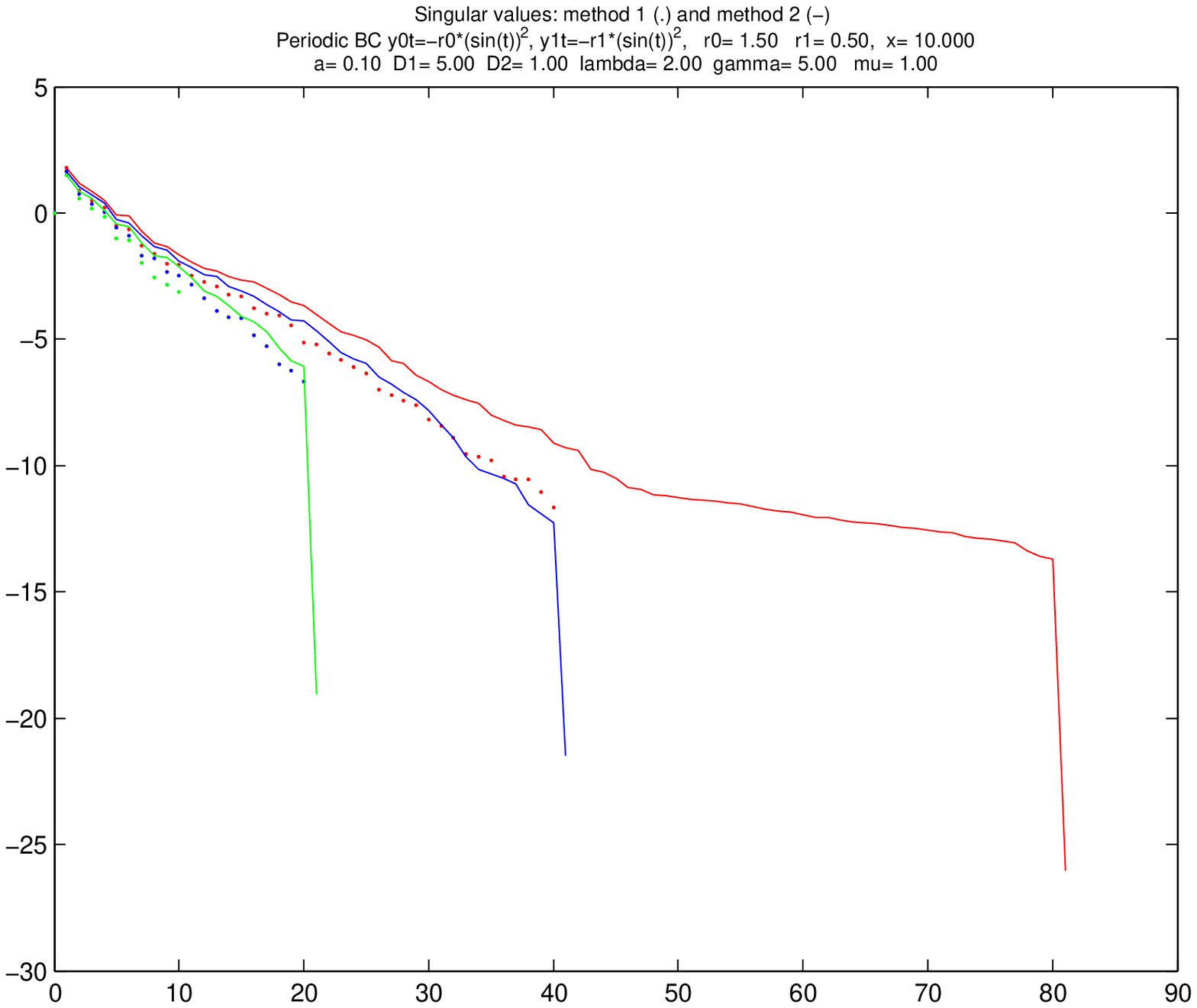, height=3.2in,width=4.5in}
\caption{Experiment C.  Distribution of the singular values  for the two ROMs and three different snapshot spacings. Dots: Method 1; continuous line: Method 2; Red: $\Delta=0.0025$; Blue: $\Delta=0.005$; Green: $\Delta=0.01$. }
\end{figure}


\begin{thebibliography}{99}
\bibitem{antoulas} A.C. ANTOULAS, {\it Approximation of large-scale Dynamical Systems}, SIAM Philadelphia (2005)

\bibitem{brunner}
H. BRUNNER, P.J. VAN DER HOUWEN, {\it The numerical solution of Volterra equations}, CWI Monographs, North-Holland (1986)

\bibitem{burkardt} J. BURKARDT, M. GUNZBURGER, L. H.-C. LEE, {\it POD and CVT-based reduced-order modeling of Navier-Stokes flows},
Computer Methods in Applied Mechanics and Engineering,
Volume 196, Number 1-3, 1, pages 337-355 (2006)

\bibitem{chapelle} D. Chapelle, A. Gariah, and J. Sainte-Marie, Galerkin approximation with proper orthogonal
decomposition: New error estimates and illustrative examples, ESAIM Math. Model.
Numer. Anal., 46, pp. 731–757 (2012)

\bibitem{chatur1}
Chaturantabut, S and Sorensen, DC, Nonlinear Model Reduction via Discrete Empirical Interpolation, SIAM Journal on Scientific Computing, 32, 2737--2764 (2010)

\bibitem{chatur2}
S. CHATURANTABUT, D.C. SORENSEN, A state space error estimate for POD --DEIM nonlinear model reduction, SIAM Journal on Numerical Analysis, 50, 1, pp. 46--63, (2012)

\bibitem{foias}
C. FOIAS, G.R. SELL, R. TEMAM, {\it Inertial manifolds for nonlinear evolutionary equations}, Journal of Differential Equations, Volume 73, Issue 2,  Pages 309--353, (1988)

\bibitem{guck}
J. GUCKENHEIMER, P. HOLMES, {\it Nonlinear Oscillations, Dynamical Systems and Bifurcations of vector Fields}, Springer (1983)

\bibitem{hinze}M. HINZE, S. VOLKWEIN, {\it Proper orthogonal decomposition surrogate models for nonlinear dynamical systems: error estimates and suboptimal control. In: Dimension Reduction of Large-Scale Systems}, Lecture Notes
in Computational Science and Engineering,
Springer 45, 261--306 (2006)

\bibitem{homescu}
C. HOMESCU, L.R. PETZOLD, R. SERBAN, {\it Error estimation for reduced-order models of dynamical systems}, SIAM Review, 
Vol. 49, No. 2, pp. 277--299 (2007)

\bibitem{iliescu}
T. ILIESCU, Z. WANG, {\it Are the snapshot difference quotients needed in the proper orthogonal decomposition?}
SIAM J. Sci.Comput. 2014, Vol. 36, No. 3, pp. A1221--A1250 (2014)

\bibitem{ito} 
K. ITO, S.S. RAVINDRAN, {\it A reduced-order method for simulation
and control of fluid flows}, J. Comp. Phys. 143, 403--425 (1998) 

\bibitem{keener}
J. KEENER, J. SNEYD, {\it Mathematical physiology}, Springer-Verlag,  New York (1998)

\bibitem{kostova}
T. KOSTOVA,R. RAVINDRAN, M. SCHONBEK, {\it FitzHugh-Nagumo revisted}, International J.Bifurcation \& Chaos, v 14, no. 3: 913-925 (2004)

\bibitem{kv1}
Kunisch K,  Volkwein S, {\it Galerkin proper orthogonal decomposition methods for parabolic
problems}, Numer. Math., 90, pp. 117--148  (2001).

\bibitem{kv2} K. KUNISCH, S. VOLKWEIN, {\it Galerkin proper orthogonal decomposition methods for a general
equation in fluid dynamics}, SIAM J. Numer. Anal., 40, pp. 492--515 (2002).

\bibitem{olmos} 
D. OLMOS, B.D. SHIZGAL, {\it Pseudospectral method of solution of the FitzHugh -- Nagumo equation}, Mathematics and Computers in Simulation 79 2258--2278 (2009) 

\bibitem{rathinam} M. RATHINAM, L. PETZOLD, {\it A new look at proper orthogonal decomposition}, SIAM J. Numer. Anal., 41, No. 5, pp. 1893--1925 (2003) 
 
\bibitem{rega}
G. REGA, H. TROGER, {\it Dimension reduction of dynamical systems: methods,models, applications}, Nonlinear Dynamics 41: 1--15  (2005) 

\bibitem{rinzel} J. RINZEL, D. TERMAN, {\it Propagation phenomena in a bistable reaction-diffusion system},  SIAM Journal on Applied Mathematics (1982) 

\bibitem{rudin} W. RUDIN, {\it Principles of mathematical analysis}, McGraw-Hill (1976)

\bibitem{serban} R. SERBAN R, C. HOMESCU, L. PETZOLD, {\it The effect of problem perturbations on nonlinear dynamical systems and their reduced order models}, SIAM J Sci Comput, Society for Industrial and Applied Mathematics Vol. 29, No. 6, pp. 2621--2643  (2007) 

\bibitem{singler} J.R. SINGLER, {\it New POD error expressions, error bounds, and asymptotic results for reduced order models of parabolic PDEs}, SIAM J. Numer. Anal., 52, pp. 852--876 (2014).

\bibitem{stoer}J. STOER, R. BULIRSCH, {\it Introduction to numerical analysis} (Texts in Applied Mathematics),  Springer-Verlag,  New York (1980)

\bibitem{sirovich} 
L. SIROVICH, {\it Turbulence and the dynamics of coherent structures}, parts I--III., Quart. Appl. Math. XLV: 561--590 (1987)

\end{thebibliography}
\end{document}